\documentclass[12pt,a4paper]{article}

\usepackage{latexsym}
\usepackage{tikz} 
\usepackage[new]{old-arrows}
\usepackage{multirow}
\usepackage{amsfonts,amsmath,amsthm}
\usepackage{amsthm}
\usepackage{amsmath}
\usepackage{amssymb}
\usepackage{mathrsfs}
\usepackage[top=1.5cm, bottom=1.8cm, left=1.5cm, right=1.5cm]{geometry}

\newenvironment{changemargin}[2]{\begin{list}{}{%
\setlength{\topsep}{0pt}%
\setlength{\leftmargin}{0pt}%
\setlength{\rightmargin}{0pt}%
\setlength{\listparindent}{\parindent}%
\setlength{\itemindent}{\parindent}%
\setlength{\parsep}{0pt plus 1pt}%
\addtolength{\leftmargin}{#1}%
\addtolength{\rightmargin}{#2}%
}\item }{\end{list}}

\newtheorem{theoreme}{Theorem}[section]
\newtheorem{proposition}[theoreme]{Proposition}%[section]
\newtheorem{corollaire}[theoreme]{Corollary}%[section]
\newtheorem{definition}[theoreme]{Definition}%[section]
\newtheorem{lemme}[theoreme]{Lemma}%[section]

\theoremstyle{remark}

\newtheorem*{exemple}{Example}

\newcommand{\finDemo}{{\hfill $\Box$} \vspace{0.5em}}
\newcommand{\debutDemo}{\noindent \textit{Proof: } \hspace{2pt}}

\newcommand{\fleche}[4]{                     % fonction
            \begin{array}{rcl} #1 & \rightarrow & #2 \\   %
                         #3 &\mapsto & #4          %
            \end{array}}

\newcommand{\fonc}[5]{                     % fonction
            \begin{array}{crll}#1 :& #2 & \rightarrow & #3 \\   %
                         &#4 &\mapsto & #5          %
            \end{array}}

\DeclareMathOperator{\End}{End}

\title{\bf Projective representations of mapping class groups in combinatorial quantization}
\author{Matthieu \textsc{Faitg}}
\date{}

\begin{document}

\maketitle

\vspace{-1em}
\noindent IMAG, Univ Montpellier, CNRS, Montpellier, France.
\\E-mail address: \texttt{matthieu.faitg@umontpellier.fr}
\vspace{2em}
%\hrulefill
%\tableofcontents
%\hrulefill
\begin{changemargin}{2cm}{2cm}
{\small
\noindent \textsc{Abstract}. \hspace{2pt} Let $\Sigma_{g,n}$ be a compact oriented surface of genus $g$ with $n$ open disks removed. The graph algebra $\mathcal{L}_{g,n}(H)$ was introduced by Alekseev--Grosse--Schomerus and Buffenoir--Roche and is a combinatorial quantization of the moduli space of flat connections on $\Sigma_{g,n}$. We construct a projective representation of the mapping class group of $\Sigma_{g,n}$ using $\mathcal{L}_{g,n}(H)$ and its subalgebra of invariant elements. Here we assume that the gauge Hopf algebra $H$ is finite-dimensional, factorizable and ribbon, but not necessarily semi-simple. We also give explicit formulas for the representation of the Dehn twists generating the mapping class group; in particular, we show that it is equivalent to a representation constructed by V. Lyubashenko using categorical methods. %We first recall the definition of $\mathcal{L}_{g,n}(H)$ and we give its main properties under these assumptions on $H$. Although the representation of the mapping class group is constructed for $n=0$, we explain in detail how to generalize the result to the case $n > 0$.
%\noindent \textsc{Abstract}. \hspace{2pt} Let $\Sigma_{g,n}$ be a compact oriented surface of genus $g$ with $n$ open disks removed. We construct a projective representation of the mapping class group of $\Sigma_{g,n}$ using the graph algebras $\mathcal{L}_{g,n}(H)$. These algebras have been introduced by Alekseev--Grosse--Schomerus and Buffenoir--Roche and are a combinatorial quantization of the moduli space of flat connections on $\Sigma_{g,n}$. Here we assume that the gauge Hopf algebra $H$ is finite-dimensional, factorizable and ribbon, but not necessarily semi-simple. We also give explicit formulas for the representation of the Dehn twists generating the mapping class group; in particular, we show that the representation is equivalent to that constructed by V. Lyubashenko. %We first recall the definition of $\mathcal{L}_{g,n}(H)$ and we give its main properties under these assumptions on $H$. Although the representation of the mapping class group is constructed for $n=0$, we explain in detail how to generalize the result to the case $n > 0$.
}
\end{changemargin}
\vspace{0.5em}

\section{Introduction}
\begin{figure}[h]
\begin{tikzpicture}[scale=0.87]
\draw [shift={(13.816832178429763,5.99)},line width=0.8pt]  plot[domain=0:3.141592653589793,variable=\t]({1*1.7075037526311627*cos(\t r)+0*1.7075037526311627*sin(\t r)},{0*1.7075037526311627*cos(\t r)+1*1.7075037526311627*sin(\t r)});
\draw [shift={(13.795787711044765,5.99)},line width=0.8pt]  plot[domain=0:3.141592653589793,variable=\t]({1*1.309575195744216*cos(\t r)+0*1.309575195744216*sin(\t r)},{0*1.309575195744216*cos(\t r)+1*1.309575195744216*sin(\t r)});
\draw [shift={(15.296936063235531,5.99)},line width=0.8pt]  plot[domain=0:1.8118338027760237,variable=\t]({1*1.4030639367644646*cos(\t r)+0*1.4030639367644646*sin(\t r)},{0*1.4030639367644646*cos(\t r)+1*1.4030639367644646*sin(\t r)});
\draw [shift={(15.296936063235531,5.99)},line width=0.8pt]  plot[domain=0:1.9894438844649243,variable=\t]({1*1.776642756206547*cos(\t r)+0*1.776642756206547*sin(\t r)},{0*1.776642756206547*cos(\t r)+1*1.776642756206547*sin(\t r)});
\draw [shift={(15.296936063235531,5.99)},line width=0.8pt]  plot[domain=2.2704846801922587:3.141592653589793,variable=\t]({1*1.3953711574291767*cos(\t r)+0*1.3953711574291767*sin(\t r)},{0*1.3953711574291767*cos(\t r)+1*1.3953711574291767*sin(\t r)});
\draw [shift={(15.296936063235531,5.99)},line width=0.8pt]  plot[domain=2.3877511379579315:3.141592653589793,variable=\t]({1*1.789237366642512*cos(\t r)+0*1.789237366642512*sin(\t r)},{0*1.789237366642512*cos(\t r)+1*1.789237366642512*sin(\t r)});
\draw [shift={(13.802618659443407,5.996993354832272)},line width=0.8pt]  plot[domain=0:3.141592653589793,variable=\t]({1*1.4991378524730479*cos(\t r)+0*1.4991378524730479*sin(\t r)},{0*1.4991378524730479*cos(\t r)+1*1.4991378524730479*sin(\t r)});
\draw [line width=0.8pt] (11.7,4.99)-- (17.494708954246594,4.99203777401749);
\draw [line width=0.8pt] (11.7,5.99)-- (12.1093284257986,5.99);
\draw [line width=0.8pt] (12.486212515300549,5.99)-- (13.508028251752195,5.99);
\draw [line width=0.8pt] (13.900815462007142,5.99)-- (15.105362906788981,5.99);
\draw [line width=0.8pt] (15.524335931060925,5.99)-- (16.7,5.99);
\draw [line width=0.8pt] (17.073578819442076,5.99)-- (17.507869270959233,5.99);
\draw [shift={(15.298192456594988,7.366022459255682)},line width=0.8pt]  plot[domain=3.8538008041679896:5.569368440196791,variable=\t]({1*2.1055979038964163*cos(\t r)+0*2.1055979038964163*sin(\t r)},{0*2.1055979038964163*cos(\t r)+1*2.1055979038964163*sin(\t r)});
\draw [shift={(15.293880815436204,5.974130858810398)},line width=0.8pt]  plot[domain=2.2962968369369396:3.1316089960333495,variable=\t]({1*1.6018682159422226*cos(\t r)+0*1.6018682159422226*sin(\t r)},{0*1.6018682159422226*cos(\t r)+1*1.6018682159422226*sin(\t r)});
\draw [shift={(15.293880815436204,5.974130858810398)},line width=0.8pt]  plot[domain=0.009962038173572503:1.9080540500790064,variable=\t]({1*1.5929876394337077*cos(\t r)+0*1.5929876394337077*sin(\t r)},{0*1.5929876394337077*cos(\t r)+1*1.5929876394337077*sin(\t r)});
\draw [shift={(13.806612855565357,7.123074732398353)},line width=0.8pt]  plot[domain=3.7832753073946166:5.638361073009996,variable=\t]({1*1.8763553623268683*cos(\t r)+0*1.8763553623268683*sin(\t r)},{0*1.8763553623268683*cos(\t r)+1*1.8763553623268683*sin(\t r)});
\draw [shift={(19.24576578255941,5.99)},line width=0.8pt]  plot[domain=0:3.141592653589793,variable=\t]({1*1.7542342174405903*cos(\t r)+0*1.7542342174405903*sin(\t r)},{0*1.7542342174405903*cos(\t r)+1*1.7542342174405903*sin(\t r)});
\draw [shift={(19.245382891279704,5.99)},line width=0.8pt]  plot[domain=0:3.141592653589793,variable=\t]({1*1.5646171087202987*cos(\t r)+0*1.5646171087202987*sin(\t r)},{0*1.5646171087202987*cos(\t r)+1*1.5646171087202987*sin(\t r)});
\draw [shift={(19.245,5.99)},line width=0.8pt]  plot[domain=0:3.141592653589793,variable=\t]({1*1.375*cos(\t r)+0*1.375*sin(\t r)},{0*1.375*cos(\t r)+1*1.375*sin(\t r)});
\draw [line width=0.8pt] (17.494708954246594,4.99203777401749)-- (21,4.99);
\draw [line width=0.8pt] (17.87,5.99)-- (20.62,5.99);
\draw [shift={(23.883555184850728,5.997624070343269)},line width=0.8pt]  plot[domain=0:3.141592653589793,variable=\t]({1*1.7542342174405903*cos(\t r)+0*1.7542342174405903*sin(\t r)},{0*1.7542342174405903*cos(\t r)+1*1.7542342174405903*sin(\t r)});
\draw [shift={(23.883172293571022,5.997624070343269)},line width=0.8pt]  plot[domain=0:3.141592653589793,variable=\t]({1*1.5646171087202987*cos(\t r)+0*1.5646171087202987*sin(\t r)},{0*1.5646171087202987*cos(\t r)+1*1.5646171087202987*sin(\t r)});
\draw [shift={(23.882789402291316,5.997624070343269)},line width=0.8pt]  plot[domain=0:3.141592653589793,variable=\t]({1*1.375*cos(\t r)+0*1.375*sin(\t r)},{0*1.375*cos(\t r)+1*1.375*sin(\t r)});
\draw [line width=0.8pt] (22.132498356537912,4.999661844360759)-- (25.637789402291318,4.997624070343269);
\draw [line width=0.8pt] (22.507789402291316,5.997624070343269)-- (25.25778940229132,5.997624070343269);
\draw [line width=0.8pt] (21,5.99)-- (21.247247909739123,5.989302899659494);
\draw [line width=0.8pt] (21.882036837396466,5.997431483584242)-- (22.129320967410138,5.997624070343269);
\draw [line width=0.8pt] (22.132498356537912,4.999661844360759)-- (21.886596832121036,4.997979206429364);
\draw [line width=0.8pt] (21.246207548250812,4.987839750755572)-- (20.999994186578796,4.990000003379588);
\draw [shift={(7.554556060708619,6.001511148703084)},line width=0.8pt]  plot[domain=0:3.141592653589793,variable=\t]({1*1.7075037526311627*cos(\t r)+0*1.7075037526311627*sin(\t r)},{0*1.7075037526311627*cos(\t r)+1*1.7075037526311627*sin(\t r)});
\draw [shift={(7.533511593323622,6.001511148703084)},line width=0.8pt]  plot[domain=0:3.141592653589793,variable=\t]({1*1.309575195744216*cos(\t r)+0*1.309575195744216*sin(\t r)},{0*1.309575195744216*cos(\t r)+1*1.309575195744216*sin(\t r)});
\draw [shift={(9.034659945514388,6.001511148703084)},line width=0.8pt]  plot[domain=0:1.8118338027760237,variable=\t]({1*1.4030639367644646*cos(\t r)+0*1.4030639367644646*sin(\t r)},{0*1.4030639367644646*cos(\t r)+1*1.4030639367644646*sin(\t r)});
\draw [shift={(9.034659945514388,6.001511148703084)},line width=0.8pt]  plot[domain=0:1.9894438844649243,variable=\t]({1*1.776642756206547*cos(\t r)+0*1.776642756206547*sin(\t r)},{0*1.776642756206547*cos(\t r)+1*1.776642756206547*sin(\t r)});
\draw [shift={(9.034659945514388,6.001511148703084)},line width=0.8pt]  plot[domain=2.2704846801922587:3.141592653589793,variable=\t]({1*1.3953711574291767*cos(\t r)+0*1.3953711574291767*sin(\t r)},{0*1.3953711574291767*cos(\t r)+1*1.3953711574291767*sin(\t r)});
\draw [shift={(9.034659945514388,6.001511148703084)},line width=0.8pt]  plot[domain=2.3877511379579315:3.141592653589793,variable=\t]({1*1.789237366642512*cos(\t r)+0*1.789237366642512*sin(\t r)},{0*1.789237366642512*cos(\t r)+1*1.789237366642512*sin(\t r)});
\draw [shift={(7.545148591741332,6.001511148703084)},line width=0.8pt]  plot[domain=0:3.141592653589793,variable=\t]({1*1.5074247094624793*cos(\t r)+0*1.5074247094624793*sin(\t r)},{0*1.5074247094624793*cos(\t r)+1*1.5074247094624793*sin(\t r)});
\draw [shift={(7.542140865125529,7.129716438011397)},line width=0.8pt]  plot[domain=3.78504805100958:5.641644250362013,variable=\t]({1*1.8804567618269643*cos(\t r)+0*1.8804567618269643*sin(\t r)},{0*1.8804567618269643*cos(\t r)+1*1.8804567618269643*sin(\t r)});
\draw [shift={(9.032426867588411,7.322376985857134)},line width=0.8pt]  plot[domain=3.8347075821896603:5.590627966512442,variable=\t]({1*2.0672882303154765*cos(\t r)+0*2.0672882303154765*sin(\t r)},{0*2.0672882303154765*cos(\t r)+1*2.0672882303154765*sin(\t r)});
\draw [shift={(8.938496818166326,6.041050601785358)},line width=0.8pt]  plot[domain=2.2945520657206804:3.168010420806889,variable=\t]({1*1.4981243856202084*cos(\t r)+0*1.4981243856202084*sin(\t r)},{0*1.4981243856202084*cos(\t r)+1*1.4981243856202084*sin(\t r)});
\draw [shift={(9.037723882278852,6.001511148703084)},line width=0.8pt]  plot[domain=0.0005593210374478744:1.9083098013248945,variable=\t]({1*1.5857206133786363*cos(\t r)+0*1.5857206133786363*sin(\t r)},{0*1.5857206133786363*cos(\t r)+1*1.5857206133786363*sin(\t r)});
\draw [line width=0.8pt] (5.437723882278853,6.001511148703084)-- (5.437723882278853,5.001511148703084);
\draw [line width=0.8pt] (5.437723882278853,5.001511148703084)-- (11.23243283652545,5.003548922720574);
\draw [line width=0.8pt] (5.437723882278853,6.001511148703084)-- (5.847052308077457,6.001511148703084);
\draw [line width=0.8pt] (6.223936397579406,6.001511148703084)-- (7.2457521340310524,6.001511148703084);
\draw [line width=0.8pt] (7.6385393442859995,6.001511148703084)-- (8.843086789067838,6.001511148703084);
\draw [line width=0.8pt] (9.262059813339782,6.001511148703084)-- (10.437723882278853,6.001511148703084);
\draw [line width=0.8pt] (10.811302701720935,6.001511148703084)-- (11.245593153238092,6.001511148703084);
\draw [line width=0.8pt] (25.637789402291318,4.997624070343269)-- (26,5);
\draw [line width=0.8pt] (25.637789402291318,5.997624070343269)-- (26,6);
\draw [shift={(19.24322702226514,7.300822577637479)},line width=0.8pt]  plot[domain=3.8396349332512414:5.586499614921624,variable=\t]({1*2.0394952698221545*cos(\t r)+0*2.0394952698221545*sin(\t r)},{0*2.0394952698221545*cos(\t r)+1*2.0394952698221545*sin(\t r)});
\draw [shift={(23.88841170623931,7.3222559223474075)},line width=0.8pt]  plot[domain=3.8400078207328283:5.579860283760005,variable=\t]({1*2.049771503067818*cos(\t r)+0*2.049771503067818*sin(\t r)},{0*2.049771503067818*cos(\t r)+1*2.049771503067818*sin(\t r)});
\draw [line width=0.8pt] (26,6)-- (26,5);
\draw[color=black] (14.065464814556085,8.225228746713648) node {$\overset{I}{B}(g)$};
\draw[color=black] (16.506783062916657,8.049594310731285) node {$\overset{I}{A}(g)$};
\draw [fill=black,shift={(13.783362211564688,7.496007527524644)},rotate=270] (0,0) ++(0 pt,3pt) -- ++(2.598076211353316pt,-4.5pt)--++(-5.196152422706632pt,0 pt) -- ++(2.598076211353316pt,4.5pt);
\draw [fill=black] (11.229765344432948,5.481433903743638) circle (1pt);
\draw [fill=black] (11.728945614264426,5.491417509140221) circle (1pt);
\draw [fill=black] (11.479355479348687,5.4864257064419295) circle (1pt);
\draw[color=black] (19.572197705594967,8.280296075623776) node {$\overset{I}{M}(g+1)$};
\draw[color=black] (24.197853334067627,8.280296075623776) node {$\overset{I}{M}(g+n)$};
\draw [fill=black] (21.290988745789484,5.53370059891846) circle (1pt);
\draw [fill=black] (21.790988745789484,5.53370059891846) circle (1pt);
\draw [fill=black] (21.540988745789484,5.53370059891846) circle (1pt);
\draw[color=black] (7.806145095075223,8.243584523017024) node {$\overset{I}{B}(1)$};
\draw[color=black] (10.247463343435793,8.023315207376505) node {$\overset{I}{A}(1)$};
\draw [fill=black,shift={(7.499751368533525,7.5082521174243375)},rotate=270] (0,0) ++(0 pt,3pt) -- ++(2.598076211353316pt,-4.5pt)--++(-5.196152422706632pt,0 pt) -- ++(2.598076211353316pt,4.5pt);
\draw [fill=black,shift={(9.024083427459018,7.587173093029878)},rotate=270] (0,0) ++(0 pt,3pt) -- ++(2.598076211353316pt,-4.5pt)--++(-5.196152422706632pt,0 pt) -- ++(2.598076211353316pt,4.5pt);
\draw [fill=black,shift={(15.281283481269876,7.567068687653846)},rotate=270] (0,0) ++(0 pt,3pt) -- ++(2.598076211353316pt,-4.5pt)--++(-5.196152422706632pt,0 pt) -- ++(2.598076211353316pt,4.5pt);
\draw [fill=black,shift={(19.21714251123375,7.554362227182415)},rotate=270] (0,0) ++(0 pt,3pt) -- ++(2.598076211353316pt,-4.5pt)--++(-5.196152422706632pt,0 pt) -- ++(2.598076211353316pt,4.5pt);
\draw [fill=black,shift={(23.82801754289718,7.561268739178991)},rotate=270] (0,0) ++(0 pt,3pt) -- ++(2.598076211353316pt,-4.5pt)--++(-5.196152422706632pt,0 pt) -- ++(2.598076211353316pt,4.5pt);
\draw [fill=black] (17.25437211155386,4.991953256886888) circle (2.5pt);
\end{tikzpicture}
\caption{To each generating loop in $\pi_1(\Sigma_{g,n}\!\!\setminus\! D, \bullet)$, we associate a matrix.}
\label{figureIntro}
\end{figure}
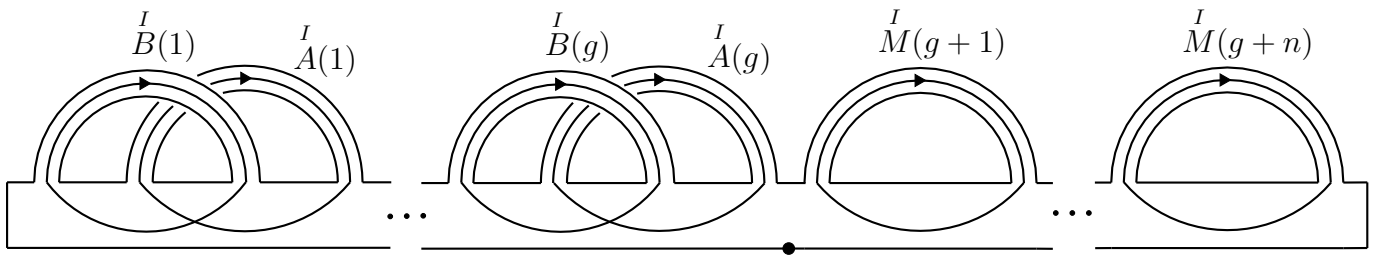

\indent Let %\let\thefootnote\relax
%\footnote{2010 Mathematics Subject Classification: Primary 16T05 ; Secondary 81R05, 81R50.\\
%\indent\hspace{0.45em} Keywords: combinatorial quantization, mapping class group, factorizable Hopf algebra.}
%\addtocounter{footnote}{-1}\let\thefootnote\svthefootnote 
$\Sigma_{g,n}$ be a compact oriented surface of genus $g$ with $n$ open disks removed. % and let $\Gamma = (V,E)$ be a filling graph. To simplify, we choose the canonical graph $\Gamma = \left(\{\bullet\}, \{b_1, a_1, \ldots, b_g, a_g, m_{g+1}, \ldots, m_{g+n}\}\right)$ with vertex $\bullet$ and edges the generators of the fundamental group. 
It is readily seen that $\Sigma_{g,n} \setminus D$ (where $D$ is an open disk) is homeomorphic to the tubular neighborhood of the graph $\Gamma$ whose edges are the generators of the fundamental group of the surface (see Figure \ref{figureIntro}); we will denote $\Sigma_{g,n}^{\mathrm{o}} = \Sigma_{g,n} \setminus D$. This particular choice of graph is not a loss of generality.
\smallskip\\
\indent Let $G$ be an algebraic Lie group (generally assumed connected and simply-connected, \textit{e.g.} $G = \mathrm{SL}_2(\mathbb{C})$). A lattice gauge field theory on $\Gamma$ is a discretization of the moduli space of flat $G$-connections on $\Sigma_{g,n}^{\mathrm{o}}$. It consists of a set of discrete connections $\mathcal{A} = G^{2g+n}$, a gauge group $\mathcal{G} = G$ and an algebra of functions $\mathbb{C}[\mathcal{A}] = \mathbb{C}[G]^{\otimes 2g+n}$ (see \textit{e.g.} \cite[2.3]{W}, \cite[Chap. 2]{labourie} for the general definitions). There is also a notion of discrete holonomy defined in a natural way. The gauge group acts on $\mathcal{A}$ (and dually on $\mathbb{C}[\mathcal{A}]$ on the right) by conjugation; the invariant functions are called classical observables.
\smallskip\\
\indent Lattice gauge field theory on $\Gamma$ is another description of the character variety of $\Sigma_{g,n}^{\mathrm{o}}$. More precisely, the discrete holonomy is a bijection between the set $\mathcal{A}/\mathcal{G}$ of discrete $G$-connections up to gauge equivalence and $\mathrm{Hom}\!\left(\pi_1(\Sigma_{g,n}^{\mathrm{o}}), G\right)/G$. The space $\mathcal{A}$ is endowed with a Poisson structure defined by Fock and Rosly \cite{FockRosly}. This Poisson structure descends to $\mathcal{A}/\mathcal{G}$ and moreover, $\mathbb{C}[\mathcal{A}/\mathcal{G}] = \mathbb{C}[\mathcal{A}]^{\mathcal{G}}$ is isomorphic to $\mathbb{C}\!\left[ \mathrm{Hom}\!\left(\pi_1(\Sigma_{g,n}^{\mathrm{o}}), G\right) \right]^G$, namely the space of functions on the character variety. Under this isomorphism, the Fock--Rosly Poisson structure corresponds to that given by the Goldman bracket, or equivalently, by the Atiyah--Bott symplectic form.
\smallskip\\
\indent The previous remarks apply to the original surface $\Sigma_{g,n}$ if we consider the subset of discrete flat connections $\mathcal{A}_f$ instead of $\mathcal{A}$. These are the discrete connections whose holonomy along the boundary of the unique face of the graph $\Gamma$ is trivial. 

\smallskip

\indent It is worthwhile to describe the algebra of functions $\mathbb{C}[\mathcal{A}]$ in terms of matrix coefficients $\overset{I}{T}{^i_j} \in \mathbb{C}[G]$ (where $I$ is a finite-dimensional $G$-module), since they linearly span $\mathbb{C}[G]$. For instance, we can construct a function $\overset{I}{A}(k){^i_j} \in \mathbb{C}[\mathcal{A}]$ by putting the function $\overset{I}{T}{^i_j}$ over the edge $a_k$ and the trivial function $1$ on the other edges. In particular, we get a matrix $\overset{I}{A}(k)$ with coefficients in $\mathbb{C}[\mathcal{A}]$, see Figure \ref{figureIntro}. The coefficients of such matrices span $\mathbb{C}[\mathcal{A}]$ as an algebra. The action of the gauge group is by conjugation, for instance $\overset{I}{A} \cdot g = \overset{I}{g} \overset{I}{A} \overset{I}{g}{^{-1}}$, where $\overset{I}{g} = \overset{I}{T}(g)$ is the representation of $g$ on $I$.
%To each edge $a_k$ (resp. $b_k$, $m_{\ell}$) in $\Gamma$, we associate a function $\overset{I}{A}(k)^i_j$ (resp. $\overset{I}{B}(k)^i_j$, $\overset{I}{M}(\ell)^i_j$) which is $\overset{I}{T}{^i_j}$ in $\mathbb{C}[G]_{a_k}$ (resp. on $\mathbb{C}[G]_{b_k}$, $\mathbb{C}[G]_{m_{\ell}}$) and the trivial function $1$ elsewhere (see \ref{figureIntro}). Then $\mathbb{C}[\mathcal{A}] = \mathbb{C}\langle \overset{I}{A}(k)^i_j,  \overset{I}{B}(k)^i_j, \overset{I}{M}(\ell)^i_j\rangle_{I, i, j, k, \ell}$. The (right) action of the gauge group $G$ is given by $\overset{I}{U}(k) \cdot h = \overset{I}{h} \overset{I}{U}(k) \overset{I}{h}{^{-1}}$, where $U$ is $A, B, M$ and $\overset{I}{h} = \overset{I}{T}(h)$ is the representation of $h \in G$ on $I$.

\smallskip

\indent In the works of Alekseev \cite{alekseev}, Alekseev--Grosse--Schomerus \cite{AGS, AGS2} and Buffenoir--Roche \cite{BR, BR2}, the Lie group $G$ is replaced by a quantum group $U_q(\mathfrak{g})$, with  $\mathfrak{g} = \mathrm{Lie}(G)$. The notions described above can be generalized in this setting. Then the graph algebra $\mathcal{L}_{g,n}(U_q(\mathfrak{g}))$ is a quantization of the Fock-Rosly Poisson structure on $\mathcal{A}$. It is an associative (non-commutative) deformation of $\mathbb{C}[\mathcal{A}]$, defined by means of equalities involving the matrices $\overset{I}{A}(k),  \overset{I}{B}(k), \overset{I}{M}(l)$ and the $R$-matrix of $U_q(\mathfrak{g})$ ($I$ now runs in the set of finite-dimensional $U_q(\mathfrak{g})$-modules). This algebra is endowed with an action of $U_q(\mathfrak{g})$, analogous to the action of the gauge group $G$ on $\mathbb{C}[\mathcal{A}]$. The multiplication in $\mathcal{L}_{g,n}(U_q(\mathfrak{g}))$ is designed so that it is an $U_q(\mathfrak{g})$-module-algebra with respect to this action. In particular, we have a subalgebra of invariant elements $\mathcal{L}^{\mathrm{inv}}_{g,n}(U_q(\mathfrak{g}))$, which is a quantized analogue of the algebra of classical observables of the initial lattice gauge field theory. 
%%%%%%%%%%%%%%%%%%%%%%%%%%%%%%%%%%%%%%%%%%%%%%%%%%%%%%
%\begin{comment}

%\end{comment}
%%%%%%%%%%%%%%%%%%%%%%%%%%%%%%%%%%%%%%%%%%%%%%%%%%
\indent These quantized algebras of functions and their generalizations appear in various works, \textit{e.g.} \cite{BFKB2, BNR, MW, BZBJ, AGPS}.
\smallskip\\
\indent The definition of the algebras $\mathcal{L}_{g,n}(U_q(\mathfrak{g}))$ is purely algebraic and we can replace the quantum group $U_q(\mathfrak{g})$ by any ribbon Hopf algebra $H$. The representation theory of $\mathcal{L}_{g,n}(H)$ and of its subalgebra of invariant elements is investigated in \cite{alekseev} when $H$ is the quantum group $U_q(\mathfrak{g})$ for $q$ generic and in \cite{AS} when $H$ is finite-dimensional and semi-simple, or a semisimple truncation of quantum group at a root of unity (the latter being defined in the setting of quasi-Hopf algebras).
\smallskip\\
\indent Moreover, in \cite{AS, AS2}, a projective representation of the mapping class group of $\Sigma_{g,n}$ based on $\mathcal{L}_{g,n}(H)$ is described. This representation is an analogue in the quantized setting of the obvious representation of the mapping class groups on $\mathbb{C}[\mathcal{A}]$ and $\mathbb{C}[\mathcal{A}_f].$%$\mathbb{C}\!\left[ \mathrm{Hom}\!\left(\pi_1(\Sigma_{g,n}), G\right) \right]^G$.
\medskip\\
\indent In this paper, we consider the algebras $\mathcal{L}_{g,n}(H)$ from a purely algebraic viewpoint, under the general assumption that the gauge algebra $H$ is finite-dimensional, factorizable and ribbon, but not necessarily semi-simple. The algebras $\mathcal{L}_{0,1}(H)$ and $\mathcal{L}_{1,0}(H)$, which are the building blocks of the theory (see Definiton \ref{definitionLgn}), and the associated projective representation of $\mathrm{SL}_2(\mathbb{Z})$, have already been studied under these assumptions in \cite{Fai18}.
\smallskip\\
\indent In section \ref{sectionLgnAlekseev}, we first quickly recall the definition and main properties of $\mathcal{L}_{0,1}(H)$ and $\mathcal{L}_{1,0}(H)$. Then we recall the definition of $\mathcal{L}_{g,n}(H)$, and we show that the Alekseev isomorphism \cite{alekseev}, which is a fundamental tool to construct representations of $\mathcal{L}_{g,n}(H)$, holds under our assumptions. In particular, when $n=0$, the Alekseev isomorphism implies that $\mathcal{L}_{g,0}(H)$ is isomorphic to a matrix algebra (because the Heisenberg double is a matrix algebra, see subsection \ref{heisenbergDouble} and \eqref{isoPsi}) and that the only indecomposable (and simple) representation of $\mathcal{L}_{g,0}(H)$ is $(H^*)^{\otimes g}$.
\smallskip\\
\indent We construct representations of the subalgebras of invariant elements $\mathcal{L}^{\mathrm{inv}}_{g,n}(H)$ in section \ref{RepInvariants} with a generalization of the method used in \cite{alekseev}. More precisely, for each representation $V$ of $\mathcal{L}_{g,n}(H)$ we associate a representation $\mathrm{Inv}(V) \subset V$ of $\mathcal{L}_{g,n}^{\mathrm{inv}}(H)$, defined by the requirement that the holonomy of a connection along the boundary of the unique face of the graph $\Gamma$ acts trivially on it.
\smallskip\\
\indent In section \ref{sectionRepMCG}, we construct a projective representation of the mapping class groups $\mathrm{MCG}(\Sigma_{g,0}^{\mathrm{o}})$ and $\mathrm{MCG}(\Sigma_{g,0})$ (we discuss the case $n > 0$ in subsection \ref{CasGeneral}). The idea of the construction is to associate an automorphism $\widetilde{f}$ of $\mathcal{L}_{g,0}(H)$ to each element $f$ of the mapping class group (Proposition \ref{liftHumphries}), called the lift of $f$. To define such a lift, we just replace generators of the fundamental group by matrices of generators of $\mathcal{L}_{g,0}(H)$ (up to some normalization), see \eqref{actionPi1} and \eqref{courbesDeviennentMatrices}. Since $\mathcal{L}_{g,0}(H)$ is isomorphic to a matrix algebra, this automorphism is inner and we get an element $\widehat{f} \in \mathcal{L}_{g,0}(H)$, unique up to scalar. Then to $f$ we associate the representation of $\widehat{f}$ on $(H^*)^{\otimes g}$ (in the case of $\Sigma_{g,0}^{\mathrm{o}}$) and on $\mathrm{Inv}\!\left((H^*)^{\otimes g}\right)$ (in the case of $\Sigma_{g,0}$). This construction was first introduced by Alekseev and Schomerus in \cite{AS} and \cite{AS2} in the semi-simple setting. Here we generalize and complete this approach with detailed proofs in the non-semi-simple setting. 
\smallskip\\
\indent  Finally, we give explicit formulas for the representation of the Dehn twists about the curves depicted in Figure \ref{figureCourbesCanoniques} (Theorem \ref{formulesExplicites}) and in particular this allows us to prove that the representation of the mapping class group described above is equivalent (Theorem \ref{thmEquivalenceReps}) to another one constructed by Lyubashenko using categorical techniques based on the coend of a ribbon category $\mathcal{C}$ satisfying some assumptions \cite{lyu95a, lyu95b, lyu96}. For this equivalence we take $\mathcal{C} = \mathrm{mod}_l(H)$, the category of finite-dimensional left modules. For works based on the Lyubashenko representation, see \textit{e.g.} \cite{FSS1, FSS2}.

\smallskip

\indent Although the two representations are equivalent, the combinatorial quantization provides additional structure and tools. Indeed, it also gives rise to a representation of the quantized version of the classical observables $\mathcal{L}_{g,n}^{\mathrm{inv}}(H)$; this is interesting because these quantum observables are related to skein theory \cite{BFKB, BFKB2}. Moreover, as a deformation of the algebra of functions on the character variety, combinatorial quantization is a natural and explicit setting to derive mapping class group representations. %, and the construction is less {\em ad hoc} than Lyubashenko theory. 
\smallskip\\
\indent To sum up, the main results of this paper are:
\begin{itemize}
%\item The construction of representations of the subalgebra of invariant elements $\mathcal{L}_{g,n}^{\mathrm{inv}}(H)$ (Theorem \ref{thmInv}),
\item The construction of a projective representation of $\mathrm{MCG}(\Sigma_{g,0}^{\mathrm{o}})$ and $\mathrm{MCG}(\Sigma_{g,0})$ (Theorem \ref{thmRepMCG}),
\item Explicit formulas for the representation of the Dehn twists about the curves of Figure \ref{figureCourbesCanoniques} (Theorem \ref{formulesExplicites}),
\item The equivalence with the Lyubashenko representation for $\mathrm{mod}_l(H)$ (Theorem \ref{thmEquivalenceReps}). 
\end{itemize}

\indent Let us conclude with a few remarks about our results and further work. First, as already said, all our constructions are explicit; this feature of the theory could be helpful to make computations when one studies the representation of the mapping class group for a given $H$ (see for instance the proof of \cite[Theorem 6.4]{Fai18} for computations in the case of the torus with $H = \overline{U}_q(\mathfrak{sl}(2))$). Second, for $H = \overline{U}_q(\mathfrak{sl}(2))$, our representations of the mapping class group should be associated to logarithmic conformal field theory in arbitrary genus. For the torus $\Sigma_{1,0}$, this is indeed the case: combining the results of \cite{FGST} and \cite{Fai18}, the projective representation of $\mathrm{SL}_2(\mathbb{Z})$ obtained via the combinatorial quantization is equivalent to that coming from logarithmic conformal field theory. Hence, a natural problem is to study in depth the representation of the mapping class group obtained for $H = \overline{U}_q(\mathfrak{sl}(2))$ (basis of the representation space, explicit formulas for the action on this basis and structure of the representation). Another question is to study the relation between $\mathcal{L}_{g,n}^{\mathrm{inv}}\!\left( \overline{U}_q(\mathfrak{sl}(2)) \right)$ and skein theory (work in progress).

%In \cite{Fai18}, we studied the case of the torus $\Sigma_{1,0}$; in particular we showed that the representation of $\mathrm{SL}_2(\mathbb{Z})$ is equivalent to that of Lyubashenko--Majid. and for the gauge algebra $H = \overline{U}_q(\mathfrak{sl}(2))$, we computed it explicitly and determined its structure. Thanks to the work of \cite{FGST}, it implies that in the case of the torus and for $H = \overline{U}_q(\mathfrak{sl}(2))$, the projective representation of $\mathrm{SL}_2(\mathbb{Z})$ obtained via the combinatorial quantization is equivalent to that obtained via logarithmic conformal field theory. A natural problem is to obtain such explicit formulas and description of the structure of the projective representation of the mapping class group in higher genus for $H = \overline{U}_q(\mathfrak{sl}(2))$. The first difficulty comes from the fact that we do not know a suitable basis of $\mathrm{Inv}\!\left(\overline{U}_q(\mathfrak{sl}(2))^{\otimes g}\right)$ (for $g=1$, we used the GTA basis and its multiplication rules \cite{GT, arike, F}). However, if a suitable basis is constructed, it would be easier to compute explicitly the representation in the setting of combinatorial quantization than in Lyubashenko setting (see the proof of \cite[Theorem 6.4]{Fai18} for the computation in the combinatorial quantization setting for the torus).

\bigskip

\noindent \textbf{Acknowledgments.}\hspace{4pt} I am grateful to my advisors, St\'ephane Baseilhac and Philippe Roche, for their regular support and their useful remarks. I thank A. Gainutdinov for inviting me to present this work at the algebra seminar of the University of Hamburg.
\medskip\\
\noindent \textbf{Notations.} \hspace{4pt} If $A$ is a $\mathbb{C}$-algebra, $V$ is a finite-dimensional $A$-module and $x \in A$, we denote by $\overset{V}{x} \in \End_{\mathbb{C}}(V)$ the representation of $x$ on the module $V$. Similarly, if $X \in A^{\otimes n}$ and if $V_1, \ldots , V_n$ are $A$-modules, we denote by $\overset{V_1 \ldots V_n}{X}$ the representation of $X$ on $V_1 \otimes \ldots \otimes V_n$. 
Here we consider only finite-dimensional representations, hence $H$-module implicitly means finite-dimensional $H$-module.
\smallskip\\
%\indent We use integral indices to describe embeddings $A^{\otimes m} \hookrightarrow A^{\otimes n}$ with $m \leq n$, see \cite[VIII.2]{kassel}.
%\smallskip\\ 
\indent Let $\mathrm{Mat}_m(A) = \mathrm{Mat}_m(\mathbb{C}) \otimes A$. Every $M \in \mathrm{Mat}_m(A)$ is written as $M = \sum_{i,j} E^i_j \otimes M^i_j$, where $E^i_j$ is the matrix with $1$ at the intersection of the $i$-th row and the $j$-th column and $0$ elsewhere. If $f : A \to A$ is a morphism, then we define $f(M)$ by $f(M) = \sum_{i,j} E^i_j \otimes f(M^i_j)$. Let moreover $N = \sum_{i,j}E^i_j \otimes N^i_j \in \mathrm{Mat}_n(\mathbb{C}) \otimes A$. We embed $M,N$ in $\mathrm{Mat}_m(\mathbb{C}) \otimes \mathrm{Mat}_n(\mathbb{C}) \otimes A = \mathrm{Mat}_{mn}(A)$ by
\[ M_1 = \sum_{i,j}E^i_j \otimes \mathbb{I}_n \otimes M^i_j \in \mathrm{Mat}_m(\mathbb{C}) \otimes \mathrm{Mat}_n(\mathbb{C}) \otimes A, \:\:\:\:
N_2 = \sum_{i,j} \mathbb{I}_m \otimes E^i_j \otimes N^i_j \in \mathrm{Mat}_m(\mathbb{C}) \otimes \mathrm{Mat}_n(\mathbb{C}) \otimes A. \]
where $\mathbb{I}_k$ is the identity matrix of size $k$. Then $M_1N_2$ (resp. $N_2M_1$) contains all the possible products of coefficients of $M$ (resp. of $N$) by coefficients of $N$ (resp. of $M$): $(M_1N_2)^{ik}_{j\ell} = M^i_jN^k_{\ell}$ (resp. $(N_2M_1)^{ik}_{j\ell} = N^k_{\ell}M^i_j$).
\smallskip\\
\indent In order to simplify notation we use Sweedler's notation (see \cite[Not. III.1.6]{kassel}) without summation sign for coproducts, that is we write
\[ \Delta(x) = x' \otimes x'', \:\:\: \Delta^{(2)}(x) = (\Delta \otimes \mathrm{id}) \circ \Delta(x) = x' \otimes x'' \otimes x''', \:\:\: \ldots, \:\:\: \Delta^{(n)}(x) = x^{(1)} \otimes \ldots \otimes x^{(n+1)}. \]
We write the universal $R$-matrix as $R = a_i \otimes b_i$ with implicit summation on $i$ and define $R' = b_i \otimes a_i$. We also denote $RR' = X_i \otimes Y_i$, $(RR')^{-1} = \overline{X}_i \otimes \overline{Y}_i$.
\smallskip\\ 
\indent The symbol ``?'' will mean a variable in functional constructions. For instance if $H$ is a finite-dimensional Hopf algebra and $\varphi, \psi \in H^*$, $a,b \in H$, then for all $x, y \in H$, $\varphi(?a) : x \mapsto \varphi(xa)$, $\varphi(?a) \otimes \psi(b?) : x\otimes y \mapsto \varphi(xa)\psi(by)$ and $\varphi(?a)\psi(b?) : x \mapsto \varphi(x'a)\psi(bx'')$ (thanks to the dual Hopf algebra structure on $H^*$, see below).

\section{Preliminaries}\label{preliminaries}
\indent In all this paper, $H$ is a finite-dimensional, factorizable, ribbon Hopf algebra.
\subsection{Factorizable ribbon Hopf algebras}
We recall basic facts about Hopf algebras. For more details, see \cite{kassel}.% (and note that we take back the notations of \cite[Section 2]{Fai18}). 
\smallskip\\
\indent If $I$ is a (finite-dimensional) $H$-module, we denote by $\overset{I}{T} \in \mathrm{Mat}_{\dim(I)}(H^*)$ the matrix defined by $\overset{I}{T}(x) = \overset{I}{x}$. Since $H$ is finite-dimensional, the coefficients of the matrices $\overset{I}{T}$ span $H^*$ when $I$ runs in the set of $H$-modules. We assume that $H$ is factorizable, which means that the coefficients of the matrices $(\overset{I}{T} \otimes \mathrm{id})(RR')$ span $H$ when $I$ runs in the set of $H$-modules. Let $R^{(+)} = R$, $R^{(-)} = R'^{-1}$, and consider the matrices $\overset{I}{L}{^{(\pm)}} = (\overset{I}{T} \otimes \mathrm{id})(R^{(\pm)})$. Since $H$ is factorizable, the coefficients of the matrices $\overset{I}{L}{^{(+)}}, \overset{I}{L}{^{(-)}}$ generate $H$ as an algebra when $I$ runs in the set of $H$-modules. As a consequence of the properties of the universal $R$-matrix (see \cite[VIII.2]{kassel}), we have the following relations: 
\begin{equation}\label{propertiesL}
\begin{aligned}
&\overset{I}{L} \,\!^{(\epsilon)}_1\overset{J}{L} \,\!^{(\epsilon)}_2 = \overset{\!\!\!\!\!I\otimes J}{L^{(\epsilon)}_{12}}, \:\:\:\:\: \Delta(\overset{I}{L} \,\!^{(\epsilon)}\,\!^a_b) = \sum_i\overset{I}{L} \,\!^{(\epsilon)}\,\!^i_b \otimes \overset{I}{L} \,\!^{(\epsilon)}\,\!^a_i,\\
&\overset{IJ}{R}\,\!^{(\epsilon)}_{12} \overset{I}{L} \,\!^{(\epsilon)}_1 \overset{J}{L} \,\!^{(\sigma)}_2 = \overset{J}{L} \,\!^{(\sigma)}_2 \overset{I}{L} \,\!^{(\epsilon)}_1 \overset{IJ}{R}\,\!^{(\epsilon)}_{12} \:\:\:\:\, \forall \, \epsilon, \sigma \in \{\pm\}\\
&\overset{IJ}{R}\,\!^{(\epsilon)}_{12} \overset{I}{L} \,\!^{(\sigma)}_1 \overset{J}{L} \,\!^{(\sigma)}_2 = \overset{J}{L} \,\!^{(\sigma)}_2 \overset{I}{L} \,\!^{(\sigma)}_1 \overset{IJ}{R}\,\!^{(\epsilon)}_{12} \:\:\:\: \forall \, \epsilon, \sigma \in \{\pm\}.
\end{aligned}
\end{equation}
Recall that the Drinfeld element $u$ (see \cite[VIII.4]{kassel}) and its inverse are:
\begin{equation}\label{u}
u = S(b_i)a_i = b_iS^{-1}(a_i), \:\:\:\:\: u^{-1} = S^{-2}(b_i)a_i = S^{-1}(b_i)S(a_i) = b_iS^2(a_i).
\end{equation}
We assume that $H$ contains a ribbon element $v$ (see \cite[XIV.6]{kassel}); it satisfies
\begin{equation}\label{ribbon}
v \text{ is central and invertible, } \:\:\:\:\: v^2 = uS(u), \:\:\:\:\: \Delta(v) = (R'R)^{-1} v \otimes v, \:\:\:\:\: S(v) = v.
\end{equation}
Then $H$ contains a canonical pivotal element $g = uv^{-1}$. It satisfies $\Delta(g) = g \otimes g$ and $S^2(x) = gxg^{-1}$ for all $x \in H$.
\smallskip\\
\indent We denote by $\mathcal{O}(H)$ the vector space $H^*$ endowed with the dual Hopf algebra structure, which in terms of matrix coefficients is:
\begin{equation}\label{structureOH}
\overset{I}{T_1}\overset{J}{T_2} = \overset{I \otimes J}{T}\!\!\!_{12}, \:\: 1_{H^*} = \overset{\mathbb{C}}{T}, \:\: \Delta(\overset{I}{T^{\, a}_{\, b}}) = \sum_i\overset{I}{T^{\, a}_{\, i}} \otimes \overset{I}{T^{\, i}_{\, b}}, \:\: \varepsilon(\overset{I}{T}) = I_{\dim(I)}, \:\: S(\overset{I}{T}) = \overset{I}{T}{^{-1}}
\end{equation}
where $\mathbb{C}$ is the trivial representation, so $\overset{\mathbb{C}}{T} = \varepsilon$, the counit of $H$. In particular, in $\mathcal{O}(H)$ holds the following exchange relation:
\begin{equation*}
\overset{IJ}{R}_{12} \overset{I}{T}_1 \overset{J}{T}_2 = \overset{J}{T}_2 \overset{I}{T}_1 \overset{IJ}{R}_{12}.
\end{equation*}
\indent Since $H$ is finite-dimensional, it exists right and left integrals $\mu^r, \mu^l \in \mathcal{O}(H)$ defined by
\[ \forall \, \varphi \in \mathcal{O}(H), \:\:\: \mu^r \varphi = \varepsilon(\varphi)\mu^r, \:\:\: \varphi \mu^l = \varepsilon(\varphi)\mu^l. \]
They are unique up to scalar and we fix $\mu^l = \mu^r \circ S$. Moreover, it holds
\begin{align}
&\forall \, h \in H, \:\:\: \mu^r(h?) \varphi = \mu^r\!\left(h'?\right) \varphi\!\left(S^{-1}(h'')\right), \label{integraleShifte}\\
&\mu^l = \mu^r(g^2 ?), \label{muLmuRgCarre}\\
&\forall\, x,y \in H, \:\:\: \mu^r(xy) = \mu^r(S^2(y)x), \:\:\: \mu^l(xy) = \mu^l(S^{-2}(y)x), \label{quasiCyclic}\\
&\forall\, x,y \in H, \:\:\: \mu^r(gxy) = \mu^r(gyx), \:\:\: \mu^l(g^{-1}xy) = \mu^l(g^{-1}yx) \label{integraleShifteSLF}.
\end{align}
\noindent These properties are well-known, for proofs see \textit{e.g.} \cite[Prop. 5.3, Lemma 5.9, Lemma 5.10]{Fai18} and the references therein; \eqref{integraleShifte} is easy, \eqref{integraleShifteSLF} is an obvious consequence of \eqref{quasiCyclic}. %Note that  \eqref{muLmuRgCarre}--\eqref{integraleShifteSLF} hold because of the above assumptions on $H$.

\subsection{Heinsenberg double of $\mathcal{O}(H)$}\label{heisenbergDouble}
\indent Let $H$ be a Hopf algebra. We recall the definition of the Heisenberg double $\mathcal{H}(\mathcal{O}(H))$ (see \textit{e.g.} \cite[4.1.10]{Mon}). As a vector space, $\mathcal{H}(\mathcal{O}(H)) = \mathcal{O}(H) \otimes H$. We identify $\psi \otimes 1$ with $\psi \in \mathcal{O}(H)$ and $1\otimes h$ with $h \in H$. Then the structure of algebra on $\mathcal{H}(\mathcal{O}(H))$ is defined by the following conditions:
\begin{itemize}
\item $\mathcal{O}(H) \otimes 1$ and $1 \otimes H$ are subalgebras of $\mathcal{H}(\mathcal{O}(H))$,
\item Under the previous identifications, we have the exchange rule 
\begin{equation}\label{relDefHeisenberg}
h\psi = \psi(?h')h''
\end{equation}
where $\psi(?z) \in \mathcal{O}(H)$ is defined by $\psi(?z)(x) = \psi(xz)$.
\end {itemize}
In terms of matrices, the exchange relation is
\begin{equation}\label{echangeHeisenberg}
\overset{I}{L}\,\!^{(\pm)}_1 \overset{J}{T}_2 = \overset{J}{T}_2 \overset{I}{L}\,\!^{(\pm)}_1 \overset{IJ}{R}\,\!^{(\pm)}_{12}.
\end{equation}
\indent There is a faithful representation $\triangleright$ of $\mathcal{H}(\mathcal{O}(H))$ on $\mathcal{O}(H)$ (see \cite[Lem. 9.4.2]{Mon}) defined by
\begin{equation}\label{repHO}
\psi \triangleright \varphi = \psi\varphi, \:\:\:\:\: h \triangleright \varphi = \varphi(?h).
\end{equation}
Hence we have an injective morphism $\rho : \mathcal{H}(\mathcal{O}(H)) \to \End_{\mathbb{C}}(H^*)$; by equality of the dimensions, it follows that 
\[ \mathcal{H}(\mathcal{O}(H)) \cong \End_{\mathbb{C}}(H^*). \]
In particular, the elements of $\mathcal{H}(\mathcal{O}(H))$ can be defined by their action on $\mathcal{O}(H)$ under $\triangleright$.  In terms of matrices, the representation $\triangleright$ is
\begin{equation}\label{repTriangleHeisenberg}
\overset{I}{T_1}\triangleright \overset{J}{T}_2 = \overset{I \otimes J}{T}\!\!\!_{12}, \:\:\:\:\:\:\: \overset{I}{L}\,\!^{(\pm)}_1 \triangleright \overset{J}{T}_2 = (\overset{I}{a_i^{(\pm)}})_1 \, b_i^{(\pm)} \triangleright \overset{J}{T}_2 = (\overset{I}{a_i^{(\pm)}})_1 \overset{J}{T}_2 (\overset{I}{b_i^{(\pm)}})_2 = \overset{J}{T}_2 \overset{IJ}{R}\,\!^{(\pm)}_{12}
\end{equation}
where $R^{(\pm)} = a_i^{(\pm)} \otimes b_i^{(\pm)}$.
\smallskip\\
\indent For $h \in H$, let $\widetilde{h}\in \mathcal{H}(\mathcal{O}(H))$ be defined by 
\begin{equation}\label{operateursTilde}
\widetilde{h} \triangleright \varphi = \varphi(S^{-1}(h)?).
\end{equation}
It is easy to see that
\begin{equation}\label{echangeTilde}
\forall\, g \in H, \: \forall \, \psi \in \mathcal{O}(H), \:\:\:\: \widetilde{g}\widetilde{h} = \widetilde{gh}, \:\:\: g\widetilde{h} = \widetilde{h}g, \:\:\: \widetilde{h}\psi = \psi\!\left(S^{-1}(h'')?\right)h'.
\end{equation}
Applying this to the matrices $\overset{I}{L}{^{(\pm)}}$ of generators of $H$, we define
\[ \overset{I}{\widetilde{L}}{^{(+)}} = \overset{I}{a_i} \widetilde{b_i}, \:\:\: \overset{I}{\widetilde{L}}{^{(-)}} = \overset{I}{S^{-1}(b_i)} \widetilde{a_i} \in \mathrm{Mat}_{\dim(I)}\!\left(\mathcal{H}(\mathcal{O}(H))\right) \]
or equivalently $\overset{I}{\widetilde{L}}{^{(\pm)}_1} \triangleright \overset{J}{T}_2 = \overset{IJ}{R}\,\!^{(\pm)-1}_{12} \overset{J}{T}_2$. Using the standard properties of the $R$-matrix, it is not difficult to show the following relations:
\begin{equation}\label{LTilde}
\begin{array}{l}
\overset{I}{\widetilde{L}}{^{(\epsilon)}_1} \overset{J}{\widetilde{L}}{^{(\epsilon)}_2} = \overset{I\otimes J}{\widetilde{L}}\!\!{^{(\epsilon)}_{12}}, \:\:\:\:\:\:\:\:\:\:\:\:\:\:\:\:\: \overset{I}{\widetilde{L}}{^{(\epsilon)}_1}\overset{J}{L}{^{(\sigma)}_2} =  \overset{J}{L}{^{(\sigma)}_2}\overset{I}{\widetilde{L}}{^{(\epsilon)}_1}, \:\:\:\:\:\:\:\:\:\:\:\:\:\:\:\:\:  \overset{IJ}{R}{^{(\epsilon)}_{12}} \overset{I}{\widetilde{L}}{^{(\epsilon)}_1}\overset{J}{T}_2 = \overset{J}{T}_2\overset{I}{\widetilde{L}}{^{(\epsilon)}_1},\\
\overset{IJ}{R}\,\!^{(\epsilon)}_{12} \overset{I}{\widetilde{L}} \,\!^{(\epsilon)}_1 \overset{J}{\widetilde{L}} \,\!^{(\sigma)}_2 = \overset{J}{\widetilde{L}} \,\!^{(\sigma)}_2 \overset{I}{\widetilde{L}} \,\!^{(\epsilon)}_1 \overset{IJ}{R}\,\!^{(\epsilon)}_{12} \:\:\:\: \forall\, \epsilon, \sigma \in \{\pm\}, \:\:\:\:\:\:\:\:\:  \overset{IJ}{R}\,\!^{(\epsilon)}_{12} \overset{I}{\widetilde{L}} \,\!^{(\sigma)}_1 \overset{J}{\widetilde{L}} \,\!^{(\sigma)}_2 = \overset{J}{\widetilde{L}} \,\!^{(\sigma)}_2 \overset{I}{\widetilde{L}} \,\!^{(\sigma)}_1 \overset{IJ}{R}\,\!^{(\epsilon)}_{12} \:\:\:\: \forall\, \epsilon, \sigma \in \{\pm\}.
\end{array}
\end{equation}

\section{Definition of $\mathcal{L}_{g,n}(H)$ and the Alekseev isomorphism}\label{sectionLgnAlekseev}
\indent Recall that $H$ is a finite-dimensional factorizable ribbon Hopf algebra. The algebras $\mathcal{L}_{g,n}(H)$ where introduced by Alekseev for $H = U_q(\mathfrak{g})$, which gave a presentation of them by generators and relations close to \eqref{PresentationLgn}. Here we will define $\mathcal{L}_{g,n}(H)$ using the braided tensor product, as in \cite{AS2}. This has the advantage to show immediately that $\mathcal{L}_{g,n}(H)$ is a $H$-module-algebra and to emphasize the role of the two building blocks of the theory, namely $\mathcal{L}_{0,1}(H)$ and $\mathcal{L}_{1,0}(H)$. We quickly recall the main properties of these building blocks, and we refer to \cite{Fai18} for more details about them under our assumptions on $H$.

\subsection{Definition and properties of $\mathcal{L}_{0,1}(H)$ and $\mathcal{L}_{1,0}(H)$}
\indent Let $\mathrm{T}(H^*)$ be the tensor algebra associated to $H^*$, which by definition is spanned by all the formal products $\psi_1 \cdots \psi_n$ of elements of $H^*$, modulo the obvious multilinear relations. There is a canonical injection $j : H^* \to \mathrm{T}(H^*)$ and we denote $\overset{I}{M} = j(\overset{I}{T})$.
\begin{definition}\label{defL01}
The loop algebra $\mathcal{L}_{0,1}(H)$ is the quotient of $\mathrm{T}(H^*)$ by the following fusion relations:
\[ \overset{I \otimes J}{M}\!_{12} = \overset{I}{M}_1\overset{IJ}{(R')}_{12}\overset{J}{M}_2\overset{IJ}{(R')}{^{-1}_{12}} \]
for all finite-dimensional $H$-modules $I,J$.
\end{definition}
See \eqref{produitL01Explicite} for an explicit description of the product in $\mathcal{L}_{0,1}(H)$. The matrix coefficients $\overset{I}{M}{^i_j}$ for all $I,i,j$ linearly span $\mathcal{L}_{0,1}(H)$. The following exchange relation, called reflection equation, holds in $\mathcal{L}_{0,1}(H)$:
\[ \overset{IJ}{R}_{12}\overset{I}{M}_1\overset{IJ}{(R')}_{12}\overset{J}{M}_2 = \overset{J}{M}_2\overset{IJ}{R}_{12}\overset{I}{M}_1\overset{IJ}{(R')}_{12}. \]
An important fact is that $\mathcal{L}_{0,1}(H)$ is endowed with a structure of left $\mathcal{O}(H)$-comodule-algebra $\Omega : \mathcal{L}_{0,1}(H) \to \mathcal{O}(H) \otimes \mathcal{L}_{0,1}(H)$ (\textit{i.e.} $\Omega$ is a morphism of algebras, see \cite[Def. III.7.1]{kassel}) defined by 
\[ \Omega(\overset{I}{M}{^a_b}) = \overset{I}{T^a_{\,i}}S(\overset{I}{T}\,\!^{j}_{b}) \otimes \overset{I}{M}\,\!^i_{j}. \]
If we view $\mathcal{O}(H)$ and $\mathcal{L}_{0,1}(H)$ as subalgebras of $\mathcal{O}(H) \otimes \mathcal{L}_{0,1}(H)$ in the canonical way, then $\Omega$ is simply written $\Omega(\overset{I}{M}) = \overset{I}{T}\overset{I}{M}S(\overset{I}{T})$. Equivalently, evaluating the coaction $\Omega$ on $H$, $\mathcal{L}_{0,1}(H)$ is endowed with a structure of right $H$-module-algebra (see \cite[Def. V.6.1]{kassel} for this notion) defined by 
\begin{equation}\label{actionL01}
\overset{I}{M} \cdot h = \overset{I}{h'}\overset{I}{M}\overset{I}{S(h'')}
\end{equation}
 for $h \in H$. Moreover, if we endow $H$ with the right adjoint action defined by $a \cdot h = S(h')ah''$ ($a, h \in H$), then 
\begin{equation}\label{isoPsi01}
\fonc{\Psi_{0,1}}{\mathcal{L}_{0,1}(H)}{H}{\overset{I}{M}}{(\overset{I}{T} \otimes \mathrm{id})(RR') = \overset{I}{L}{^{(+)}}\overset{I}{L}{^{(-)-1}}}
\end{equation}
is an isomorphism of $H$-module-algebras. In particular, $\mathcal{L}_{0,1}^{\mathrm{inv}}(H) \cong \mathcal{Z}(H)$, where $\mathcal{L}_{0,1}^{\mathrm{inv}}(H)$ is the space of coinvariants (that is, the elements such that $x \cdot h = \varepsilon(h)$ for all $h \in H$ or equivalently $\Omega(x) = 1 \otimes x$). Moreover, the matrices $\overset{I}{M}$ are invertible.
\medskip\\ 
\indent Now consider the free product $\mathcal{L}_{0,1}(H) \ast \mathcal{L}_{0,1}(H)$. Let $j_1$ (resp. $j_2$) be the canonical algebra embeddings in the first (resp. second) copy of $\mathcal{L}_{0,1}(H)$ in $\mathcal{L}_{0,1}(H) \ast \mathcal{L}_{0,1}(H)$, and define $\overset{I}{A} = j_1(\overset{I}{M})$, $\overset{I}{B} = j_2(\overset{I}{M})$.
\begin{definition}
The handle algebra $\mathcal{L}_{1,0}(H)$ is the quotient of $\mathcal{L}_{0,1}(H) \ast \mathcal{L}_{0,1}(H)$ by the following exchange relations:
\[ \overset{IJ}{R}_{12}\overset{I}{B}_1\overset{IJ}{(R')}_{12}\overset{J}{A}_2 = \overset{J}{A}_2\overset{IJ}{R}_{12}\overset{I}{B}_1\overset{IJ}{R}{^{-1}_{12}} \]
for all finite-dimensional $H$-modules $I, J$.
\end{definition}
Similarly to $\mathcal{L}_{0,1}(H)$, $\mathcal{L}_{1,0}(H)$ is endowed with a structure of left $\mathcal{O}(H)$-comodule-algebra structure $\Omega : \mathcal{L}_{1,0}(H) \to \mathcal{O}(H) \otimes \mathcal{L}_{1,0}(H)$ defined by
\[ \Omega(\overset{I}{A}) = \overset{I}{T}\overset{I}{A}S(\overset{I}{T}), \:\:\:\: \Omega(\overset{I}{B}) = \overset{I}{T}\overset{I}{B}S(\overset{I}{T}). \]
As previously, it is equivalent to deal with the right action defined by $\overset{I}{A} \cdot h = \overset{I}{h'}\overset{I}{A}\overset{I}{S(h'')}$, $\overset{I}{B} \cdot h = \overset{I}{h'}\overset{I}{B}\overset{I}{S(h'')}$. The map
\[\begin{array}{crll}
\Psi_{1,0} :& \mathcal{L}_{1,0}(H) & \rightarrow & \mathcal{H}(\mathcal{O}(H)) \\ 
                         & \overset{I}{A} &\mapsto & \overset{I}{L}\,\!^{(+)}\overset{I}{L}\,\!^{(-)-1}\\
                         & \overset{I}{B} &\mapsto &  \overset{I}{L}\,\!^{(+)}\overset{I}{T}\overset{I}{L}\,\!^{(-)-1}\\
\end{array}\]
is an isomorphism of algebras (see \cite{Fai18} for a proof). It follows that $\mathcal{L}_{1,0}(H)$ is isomorphic to a matrix algebra, and in particular has trivial center.

\subsection{Braided tensor product and definition of $\mathcal{L}_{g,n}(H)$}\label{defLgn}
Let $\mathrm{mod}_r(H)$ be the category of finite-dimensional right $H$-modules (or, equivalently, of finite-dimensional left $H$-comodules). The braiding in $\mathrm{mod}_r(H)$ is given by:
\[ \fonc{c_{I,J}}{I \otimes J}{J \otimes I}{v \otimes w}{w \cdot a_i \otimes v \cdot b_i} \]
with $R = a_i \otimes b_i$. Let $(A, m_A, 1_A)$ and $(B, m_B, 1_B)$ be two algebras in $\mathrm{mod}_r(H)$ (that is, $H$-module-algebras), and define:
\begin{align*}
&m_{A \widetilde{\otimes} B} = (m_A \otimes m_B) \circ (\mathrm{id}_A \otimes c_{B,A} \otimes \mathrm{id}_B) \: : \: A \otimes B \to A \otimes B,\\  
&1_{A \widetilde{\otimes} B} = 1_A \otimes 1_B \: : \: \mathbb{C} \to A \otimes B.
\end{align*}
This endows $A \otimes B$ with a structure of algebra in $\mathrm{mod}_r(H)$, denoted $A\,\widetilde{\otimes}\,B$ and called braided tensor product of $A$ and $B$ (see \cite[Lemma 9.2.12]{majid}).  Note that $\widetilde{\otimes}$ is associative.
\\\indent There are two canonical algebra embeddings $j_A, j_B : A, B \hookrightarrow A\,\widetilde{\otimes}\,B$ respectively defined by $j_A(x) = x \otimes 1_B$, $j_B(y) = 1_A \otimes y$. We identify $x \in A$ (resp. $y \in B$) with $j_A(x) \in A\,\widetilde{\otimes}\,B$ (resp. $j_B(y)$). Under these identifications, the multiplication rule in $A\,\widetilde{\otimes}\,B$ is entirely given by:
\begin{equation}\label{braidedProduct}
\forall \, x \in A, \forall \, y \in B, \:\: yx = (x \cdot a_i)(y \cdot b_i) .
\end{equation}
\indent Since $\mathcal{L}_{0,1}(H)$ and $\mathcal{L}_{1,0}(H)$ are algebras in $\mathrm{mod}_r(H)$, we can apply the braided tensor product to them.
\begin{definition}\label{definitionLgn}
$\mathcal{L}_{g,n}(H)$ is the $H$-module-algebra $\mathcal{L}_{1,0}(H)^{\widetilde{\otimes} g} \, \widetilde{\otimes} \, \mathcal{L}_{0,1}(H)^{\widetilde{\otimes} n}$.
\end{definition}
\noindent It is useful to keep in mind that the $H$-module-algebra $\mathcal{L}_{g,n}(H)$ is associated with the surface $\Sigma_{g,n}\setminus D$; in order to make this precise we now define the matrices introduced in Figure \ref{figureIntro}. There are canonical algebra embeddings $j_i : \mathcal{L}_{1,0}(H) \hookrightarrow \mathcal{L}_{g,n}(H)$ for $1 \leq i \leq g$ and $j_i : \mathcal{L}_{0,1}(H) \hookrightarrow \mathcal{L}_{g,n}(H)$ for $g+1 \leq i \leq g+n$, given by $j_i(x) = 1^{\otimes i-1} \otimes x \otimes 1^{\otimes g+n-i}$. Define 
\[ \overset{I}{A}(i) = j_i(\overset{I}{A}), \:\:\: \overset{I}{B}(i) = j_i(\overset{I}{B}) \: \text{ for } 1 \leq i \leq g \:\: \text{ and } \:\: \overset{I}{M}(i) = j_i(\overset{I}{M}) \: \text{ for } g+1 \leq i \leq g+n. \]
\indent Relation \eqref{braidedProduct} indicates that $\mathcal{L}_{g,n}(H)$ is an exchange algebra. Let us write the exchange relations in a matrix form. Let $\overset{I}{U}$ be $\overset{I}{A}$ or $\overset{I}{B}$ or $\overset{I}{M}$, let $\overset{J}{V}$ be $\overset{J}{A}$ or $\overset{J}{B}$ or $\overset{J}{M}$ and let $i < j$. Then, by definition of the right action and by \eqref{braidedProduct}:
\begin{align*}
\overset{J}{V}(j)_2 \overset{I}{U}(i)_1 = (\overset{I}{a'_k})_1 \overset{I}{U}(i)_1 \overset{I}{S(a''_k)}_1 (\overset{J}{b'_k})_2 \overset{J}{V}(j)_2 \overset{J}{S(b''_k)}_2 
= (\overset{I}{a_l})_1 \overset{IJ}{R}_{12} \overset{I}{U}(i)_1 \overset{IJ}{R}{^{-1}_{12}}\overset{J}{V}(j)_2 \overset{IJ}{R}_{12} \overset{J}{S(b_l)}_2
\end{align*}
where for the second equality we applied properties of the $R$-matrix and obvious commutation relations in $\mathrm{End}_{\mathbb{C}}(I) \otimes \mathrm{End}_{\mathbb{C}}(J) \otimes \mathcal{L}_{g,n}(H)$. Using that $a_m a_l \otimes S(b_l) b_m = 1 \otimes 1$ together with obvious commutation relations, we obtain the desired exchange relation:
\begin{equation*}
\overset{IJ}{R}_{12} \overset{I}{U}(i)_1 \overset{IJ}{R}{^{-1}_{12}}\overset{J}{V}(j)_2 = \overset{J}{V}(j)_2 \overset{IJ}{R}_{12} \overset{I}{U}(i)_1 \overset{IJ}{R}{^{-1}_{12}}.
\end{equation*}
To sum up, the presentation of $\mathcal{L}_{g,n}(H)$ by generators and relations is:
\begin{equation}\label{PresentationLgn}
  \left\{
      \begin{aligned}
& \overset{I \otimes J}{U}\!\!(i)_{12} = \overset{I}{U}(i)_1\,(\overset{IJ}{R'})_{12}\,\overset{J}{U}(i)_2\, (\overset{IJ}{R'}){_{12}^{-1}}& &\text{ for } 1 \leq i \leq g+n\\
& \overset{IJ}{R}_{12}\,\overset{I}{U}(i)_1\, \overset{IJ}{R}{^{-1}_{12}} \,\overset{J}{V}(j)_2 = \overset{J}{V}(j)_2\, \overset{IJ}{R}_{12}\,\overset{I}{U}(i)_1\, \overset{IJ}{R}{^{-1}_{12}}& &\text{ for } 1 \leq i < j \leq g+n\\
& \overset{IJ}{R}_{12}\,\overset{I}{B}(i)_1\, (\overset{IJ}{R'})_{12}\,\overset{J}{A}(i)_2 = \overset{J}{A}(i)_2\, \overset{IJ}{R}_{12}\,\overset{I}{B}(i)_1\, \overset{IJ}{R}{_{12}^{-1}} & &\text{ for } 1 \leq i \leq g
\end{aligned}
    \right.
\end{equation}
where $U(i)$ (resp. $V(i)$) is $A(i)$ or $B(i)$ if $1 \leq i \leq g$ and is $M(i)$ if $g+1 \leq i \leq g+n$. Such a presentation was first introduced in \cite{alekseev} and \cite{AGS}. Recall that the first line of relations is the $\mathcal{L}_{0,1}(H)$-fusion relation on each loop, the second line is the exchange relation of the braided tensor product and the third line is the $\mathcal{L}_{1,0}(H)$-exchange-relation.
\smallskip\\
\noindent \textbf{Notation.}~~ Let $\overset{I}{N} = \overset{I}{v}{^m} \overset{I}{N}{^{n_1}_1} \ldots \overset{I}{N}{^{n_l}_l} \in \mathrm{Mat}_{\dim(I)}\!\left(\mathcal{L}_{g,n}(H)\right)$, where $m, n_i \in \mathbb{Z}$  and each $N_i$ is one of the $A(j), B(j), M(k)$ for some $j$ or $k$. By definition of the right action on $\mathcal{L}_{g,n}(H)$, we have a morphism of $H$-modules $j_N : \mathcal{L}_{0,1}(H) \to \mathcal{L}_{g,n}(H)$ defined by $j_N(\overset{I}{M}) = \overset{I}{N}$. Let $x \in H \cong \mathcal{L}_{0,1}(H)$, then we denote $x_N = j_N(x)$. The following lemma is an obvious fact.
\begin{lemme}\label{injectionFusion}
If $N$ satisfies the fusion relation of $\mathcal{L}_{0,1}(H)$, $\overset{I \otimes J}{N}\!\!_{12} = \overset{I}{N}(i)_1\,\overset{IJ}{(R')}_{12}\,\overset{J}{N}(i)_2\, \overset{IJ}{(R')}{_{12}^{-1}}$, then $j_N$ is a morphism of $H$-module-algebras: $(xy)_N = x_N y_N$.
\end{lemme}
\noindent See e.g. \eqref{actionHsurFormes} for an application of this lemma.

\subsection{The Alekseev isomorphism}
\indent Consider the tensor product algebra $\mathcal{L}_{1,0}(H)^{\otimes g} \otimes \mathcal{L}_{0,1}(H)^{\otimes n}$. We have canonical algebra embeddings $j_i : \mathcal{L}_{1,0}(H) \hookrightarrow \mathcal{L}_{1,0}(H)^{\otimes g} \otimes \mathcal{L}_{0,1}(H)^{\otimes n}$ for $1 \leq i \leq g$ and $j_i : \mathcal{L}_{0,1}(H) \hookrightarrow \mathcal{L}_{1,0}(H)^{\otimes g} \otimes \mathcal{L}_{0,1}(H)^{\otimes n}$ for $g+1 \leq i \leq g+n$, defined by $j_i(x) = 1^{\otimes i-1} \otimes x \otimes 1^{\otimes g+n-i}$. Define $\overset{I}{\underline{A}}(i) = j_i(\overset{I}{A})$, $\overset{I}{\underline{B}}(i) = j_i(\overset{I}{B})$ for $1 \leq i \leq g$ and $\overset{I}{\underline{M}}(i) = j_i(\overset{I}{M})$ for $g+1 \leq i \leq g+n$. We underline these matrices to avoid confusion with prior matrices having coefficients in $\mathcal{L}_{g,n}(H)$. By definition, the exchange relation between copies in $\mathcal{L}_{1,0}(H)^{\otimes g} \otimes \mathcal{L}_{0,1}(H)^{\otimes n}$ is simply
\[ \overset{I}{\underline{U}}(i)_1\,\overset{J}{\underline{V}}(j)_2 = \overset{J}{\underline{V}}(j)_2\,\overset{I}{\underline{U}}(i)_1 \]
where $i \neq j$, $\underline{U}(i), \underline{V}(i)$ is $\underline{A}(i)$ or $\underline{B}(i)$ if $1 \leq i \leq g$ and is $\underline{M}(i)$ if $g+1 \leq i \leq g+n$.
\smallskip\\
\indent The next result is due to Alekseev (see \cite{alekseev}). Consider the matrices $\overset{I}{M}{^{(-)}} = \Psi_{0,1}^{-1}(\overset{I}{L}{^{(-)}})$ and $\overset{I}{C}{^{(-)}} = \Psi_{1,0}^{-1}(\overset{I}{L}{^{(-)}}\overset{I}{\widetilde{L}}{^{(-)}})$. Let
\begin{equation}\label{matricesAlekseev}
\begin{array}{ll}
\overset{I}{\Lambda}_1 = \mathbb{I}_{\dim(I)}, & \:\: \overset{I}{\Lambda}_i = \overset{I}{\underline{C}}{^{(-)}}(1) \ldots \overset{I}{\underline{C}}{^{(-)}}(i-1)  \:\: \text{ for } 2 \leq i \leq g+1,\\
\overset{I}{\Gamma}_{g+1} = \mathbb{I}_{\dim(I)},  & \:\: \overset{I}{\Gamma}_i = \overset{I}{\Lambda}_{g+1} \overset{I}{\underline{M}}{^{(-)}}(g+1) \ldots \overset{I}{\underline{M}}{^{(-)}}(i-1) \:\: \text{ for } g+2 \leq i \leq g+n.
\end{array}
\end{equation}
be matrices with coefficients in $\mathcal{L}_{1,0}(H)^{\otimes g} \otimes \mathcal{L}_{0,1}(H)^{\otimes n}$ (with $\mathbb{I}_s$ the identity matrix of size $s$).

\begin{proposition}\label{isoAlekseev}
The map
\[\begin{array}{crll}\alpha_{g,n} :& \mathcal{L}_{g,n}(H) = \mathcal{L}_{1,0}(H)^{\widetilde{\otimes} g} \, \widetilde{\otimes} \, \mathcal{L}_{0,1}(H)^{\widetilde{\otimes} n} & \rightarrow & \mathcal{L}_{1,0}(H)^{\otimes g} \otimes \mathcal{L}_{0,1}(H)^{\otimes n} \\
                         & \overset{I}{A}(i) &\mapsto &  \overset{I}{\Lambda}_i\,\overset{I}{\underline{A}}(i)\,\overset{I}{\Lambda}{_i^{-1}} \:\: \text{ for } 1 \leq i \leq g \\
                         & \overset{I}{B}(i) &\mapsto &  \overset{I}{\Lambda}_i\,\overset{I}{\underline{B}}(i)\,\overset{I}{\Lambda}{_i^{-1}} \:\: \text{ for } 1 \leq i \leq g \\
                         & \overset{I}{M}(i) &\mapsto &  \overset{I}{\Gamma}_i\,\overset{I}{\underline{M}}(i)\,\overset{I}{\Gamma}{_i^{-1}} \:\: \text{ for } g+1 \leq i \leq g+n \\
\end{array}\]
is an isomorphism of algebras, which we call the Alekseev isomorphism.
\end{proposition}
\debutDemo
In order to show that it is a morphism of algebras, one must check using various exchange relations that the defining relations \eqref{PresentationLgn} of $\mathcal{L}_{g,n}(H)$ are preserved under $\alpha_{g,n}$. This is a straightforward but tedious task and we will not give the details. Let us prove that $\alpha_{g,n}$ is bijective. We first show that $\alpha_{g,0}$ is surjective for all $g$ by induction. For $g=1$, $\alpha_{1,0}$ is the identity. For $g \geq 2$, we embed $\mathcal{L}_{g-1,0}(H)$ in $\mathcal{L}_{g,0}(H)$ in an obvious way by $\overset{I}{A}(i) \mapsto \overset{I}{A}(i)$ and $\overset{I}{B}(i) \mapsto \overset{I}{B}(i)$ for $1 \leq i \leq g-1$. Then the restriction of $\alpha_{g,0}$ to $\mathcal{L}_{g-1,0}(H)$ is $\alpha_{g-1,0}$, and by induction we assume that $\alpha_{g-1,0}(\mathcal{L}_{g-1,0}(H)) = \mathcal{L}_{1,0}(H)^{\otimes g-1}$. Since $\overset{I}{\Lambda}_i \in \mathrm{Mat}_{\dim(I)}\!\left( \mathcal{L}_{1,0}(H)^{\otimes i-1} \otimes \mathbb{C}^{\otimes g+1-i} \right)$, there exists matrices $\overset{I}{\mathcal{N}}_i$ $(1 \leq i \leq g)$ such that $\alpha_{g,0}(\overset{I}{\mathcal{N}}_i) = \overset{I}{\Lambda}_i$. Then $\alpha_{g,0}(\overset{I}{\mathcal{N}}{^{-1}_{i}} \overset{I}{U}(i) \overset{I}{\mathcal{N}}_{i}) = \overset{I}{\underline{U}}(i)$, with $U = A$ or $B$ and $\alpha_{g,0}$ is surjective. Similarly, for $g$ fixed and $n \geq 1$, we can embed $\mathcal{L}_{g,n-1}(H)$ into $\mathcal{L}_{g,n}(H)$ and reproduce the same reasoning. Hence $\alpha_{g,n}$ is surjective for all $g,n$. Since the domain and the range of $\alpha_{g,n}$ have the same dimension, it is an isomorphism.
\finDemo
\smallskip
\\\indent We can now generalize the isomorphisms $\Psi_{0,1}$ and $\Psi_{1,0}$ by
\begin{equation}\label{isoPsi}
\Psi_{g,n} = \left(\Psi_{1,0}^{\otimes g} \otimes \Psi_{0,1}^{\otimes n}\right) \circ \alpha_{g,n} \: : \: \mathcal{L}_{g,n}(H) \overset{\sim}{\rightarrow} \mathcal{H}(\mathcal{O}(H))^{\otimes g} \otimes H^{\otimes n}.
\end{equation}
In particular $\mathcal{L}_{g,0}(H)$ is a matrix algebra, since $\mathcal{H}(\mathcal{O}(H))$ is. 
\smallskip\\
\indent Thanks to $\Psi_{g,n}$, the representation theory of $\mathcal{L}_{g,n}(H)$ is entirely determined by the representation theory of $H$. Indeed, the only indecomposable (and simple) representation of $\mathcal{H}(\mathcal{O}(H)) \cong \End_{\mathbb{C}}(H^*)$ is $H^*$, thus it follows that the indecomposable representations of $\mathcal{L}_{g,n}(H)$ are of the form
\[ (H^*)^{\otimes g} \otimes I_1 \otimes \ldots \otimes I_n \]
where $I_1, \ldots, I_n$ are indecomposable representations of $H$. We will denote the action of $\mathcal{L}_{g,n}(H)$ on $(H^*)^{\otimes g} \otimes I_1 \otimes \ldots \otimes I_n$ by $\triangleright$, namely:
\begin{equation}\label{actionTriangle}
x \triangleright \varphi_1 \otimes \ldots \otimes \varphi_g \otimes v_1 \otimes \ldots \otimes v_n = \Psi_{g,n}(x) \cdot \varphi_1 \otimes \ldots \otimes \varphi_g \otimes v_1 \otimes \ldots \otimes v_n
\end{equation}
for $x \in \mathcal{L}_{g,n}(H)$, where $\cdot$ is the action component-by-component of $\Psi_{g,n}(x)$ on $(H^*)^{\otimes g} \otimes I_1 \otimes \ldots \otimes I_n$.

\section{Representation of $\mathcal{L}^{\mathrm{inv}}_{g,n}(H)$}\label{RepInvariants}
\indent Recall that an element $x \in \mathcal{L}_{g,n}(H)$ is invariant if $x \cdot h = \varepsilon(h)x$ for all $h \in H$, or equivalently, if $\Omega(x) = 1 \otimes x$. In this section we construct representations of the subalgebra of invariants $\mathcal{L}^{\mathrm{inv}}_{g,n}(H)$. For this, we use an idea introduced in \cite{alekseev} (the matrices $\overset{I}{C}$), but adaptated to our assumptions on $H$.% (recalled in \cite{BNR}). The difference is that in our general assumptions we can not invoke Gauss decompositions to define certain matrices, so we have to define them \textit{ad hoc}.

\subsection{The matrices $\overset{I}{C}_{g,n}$}
\indent We first consider the case of $\mathcal{L}_{1,0}(H)$. Let us define matrices
\begin{equation}\label{matriceCL10}
\overset{I}{C} = \overset{I}{v}{^{2}}\overset{I}{B}\overset{I}{A}{^{-1}}\overset{I}{B}{^{-1}}\overset{I}{A}, \:\:\:\:\: \overset{I}{C}{^{(\pm)}} = \Psi_{1,0}^{-1}(\overset{I}{L}{^{(\pm)}}\overset{I}{\widetilde{L}}{^{(\pm)}})
\end{equation}
\begin{lemme}\label{decGauss}
The following equality holds in $\mathcal{L}_{1,0}(H)$:
\[ \overset{I}{C} = \overset{I}{C}{^{(+)}}\overset{I}{C}{^{(-)-1}}. \]
Moreover, the matrices $\overset{I}{C}$ satisfy the fusion relation of $\mathcal{L}_{0,1}(H)$:
\[ \overset{I \otimes J}{C}\!\!_{12} = \overset{I}{C}_1 \, \overset{IJ}{(R')}_{12} \, \overset{J}{C}_2 \, \overset{IJ}{(R')}{^{-1}_{12}}. \]
\end{lemme}
\debutDemo
We have
\[ \Psi_{1,0}\!\left(\overset{I}{v}{^{2}}\overset{I}{B}\overset{I}{A}{^{-1}}\overset{I}{B}{^{-1}}\overset{I}{A}\right) = \overset{I}{L}{^{(+)}} \left(\overset{I}{v}{^{2}} \overset{I}{T}\overset{I}{L}{^{(+)-1}}\overset{I}{L}{^{(-)}} S(\overset{I}{T}) \right)\overset{I}{L}{^{(-)-1}}. \]
Let us simplify the middle term:
\begin{align*}
\overset{I}{v}{^{2}} \overset{I}{T} \overset{I}{S(a_i)} \overset{I}{S^{-1}(b_j)} b_i a_j S(\overset{I}{T}) &= \overset{I}{v}{^{2}} \overset{I}{T} \overset{I}{S(a_i)} \overset{I}{S^{-1}(b_j)} \overset{I}{S(a'_j)} \overset{I}{S(b'_i)} S(\overset{I}{T}) b''_i a''_j \\
&= \overset{I}{v}{^{2}}\overset{I}{T} \overset{I}{S(a_ia_k)}\overset{I}{S^{-1}(b_j b_{\ell})} \overset{I}{S(a_j)} \overset{I}{S(b_k)} S(\overset{I}{T}) b_ia_{\ell} = \overset{I}{T} \overset{I}{S(a_i)}\overset{I}{S(b_{\ell})} S(\overset{I}{T}) a_{\ell} b_i.
\end{align*}
The first equality is the exchange relation \eqref{relDefHeisenberg} in $\mathcal{H}(\mathcal{O}(H))$ and the second follows from the properties of the $R$-matrix. The third equality is obtained as follows: denoting $m : H \otimes H \to H$ the multiplication, we can write 
\begin{align*}
&v^2 S(a_ia_k) S^{-1}(b_{\ell}) S^{-1}(b_j) S(a_j) S(b_k) \otimes b_ia_{\ell} = v S(a_ia_k) S^{-1}(b_k b_{\ell}) g^{-1} \otimes b_ia_{\ell}\\
= \:&v \otimes 1 (m \circ (S \otimes S^{-1}) \otimes \mathrm{id})(R_{13} R_{12} R_{32}) g^{-1} \otimes 1 = v \otimes 1 (m \circ (S \otimes S^{-1}) \otimes \mathrm{id})(R_{32} R_{12} R_{13}) g^{-1} \otimes 1\\
= \:&v S(a_k a_i) S^{-1}(b_{\ell} b_k) g^{-1} \otimes a_{\ell} b_i = v S(a_i) S\!\left(S^{-2}(b_k)a_k\right) S^{-1}(b_{\ell})g^{-1} \otimes a_{\ell} b_i = S(a_i) S(b_{\ell}) \otimes a_{\ell}b_i.
\end{align*}
We used formula \eqref{u} for $u^{-1}$ twice, a Yang-Baxter relation and the standard properties for $g$ and $v$. Now, we have:
\[ \overset{I}{T}_1 \overset{I}{S(b_{\ell} a_i)}_1 S(\overset{I}{T})_1 a_{\ell} b_i \triangleright \overset{J}{T}_2 = \overset{I}{T}_1 \overset{I}{S(b_{\ell} a_i)}_1 S(\overset{I}{T})_1 \overset{J}{T}_2 \overset{J}{(a_{\ell} b_i)}_2
= \overset{I}{S(b_{\ell} a_i)}_1 \overset{J}{(a_{\ell} b_i)}_2 \overset{J}{T}_2 = \overset{I}{(a_i b_{\ell})}_1 \widetilde{b_i} \widetilde{a_{\ell}} \triangleright \overset{J}{T}_2. \]
For the second equality, we used that for any $h \in H$:
\[ \langle \overset{I}{S(b_{\ell} a_i)}_1 S(\overset{I}{T})_1 \overset{J}{T}_2 \overset{J}{(a_{\ell} b_i)}_2, h \rangle = \overset{I}{S(h'b_{\ell} a_i)}_1 \overset{J}{(h''a_{\ell} b_i)}_2 = \overset{I}{S(b_{\ell} a_ih')}_1 \overset{J}{(a_{\ell} b_ih'')}_2 = \langle S(\overset{I}{T})_1 \overset{I}{S(b_{\ell} a_i)}_1  \overset{J}{(a_{\ell} b_i)}_2 \overset{J}{T}_2, h \rangle. \]
Since $\triangleright$ is faithful, we finally get 
\[ \overset{I}{v}{^{2}} \overset{I}{T} \overset{I}{S(a_i)} \overset{I}{S^{-1}(b_j)} b_i a_j S(\overset{I}{T}) = \overset{I}{(a_{\ell}b_i)} \widetilde{b_{\ell}} \widetilde{a_i} = \overset{I}{\widetilde{L}}{^{(+)}}\overset{I}{\widetilde{L}}{^{(-)-1}}.\]
Hence 
\[ \Psi_{1,0}(\overset{I}{C}) = \overset{I}{L}{^{(+)}}\overset{I}{\widetilde{L}}{^{(+)}}(\overset{I}{L}{^{(-)}}\overset{I}{\widetilde{L}}{^{(-)}})^{-1} 
= \Psi_{1,0}(\overset{I}{C}{^{(+)}}\overset{I}{C}{^{(-)-1}}) \]
as desired. To prove the fusion relation, it suffices to consider $\Psi_{1,0}(\overset{I}{C}) = \overset{I}{L}{^{(+)}}\overset{I}{\widetilde{L}}{^{(+)}}\overset{I}{\widetilde{L}}{^{(-)-1}}\overset{I}{L}{^{(-)-1}}$ and to use the exchange relations in \eqref{LTilde}. This is a straightforward computation left to the reader.
\finDemo
\smallskip\\
\indent We now give the general definition. For $i \leq g$, let $\overset{I}{C}(i)$ be the embedding of $\overset{I}{C}$ previously defined on the $i$-th copy of $\mathcal{L}_{1,0}(H)$ in $\mathcal{L}_{g,n}(H)$.
\begin{definition} 
$\overset{I}{C}_{g,n} = \overset{I}{C}(1) \ldots \overset{I}{C}(g) \overset{I}{M}(g+1) \ldots \overset{I}{M}(g+n)$.
\end{definition}
\noindent Geometrically (see Figure \ref{figureIntro}), for each $I$ the matrix $\overset{I}{C}_{g,n}$ corresponds to the holonomy along the boundary of the unique face of the graph $\Gamma$ defined in the Introduction.
\smallskip\\
\indent There is a decomposition analogous to Lemma \ref{decGauss}, which was the case $g=1, n=0$. %Consider the matrices $\overset{I}{C}{^{(\pm)}}(1)$,  $\ldots$,  $\overset{I}{C}{^{(\pm)}}(g)$,  $\overset{I}{M}{^{(\pm)}}(g+1)$,  $\ldots$,  $\overset{I}{M}{^{(\pm)}}(g+n)$ $\in \mathrm{Mat}\!\left(\mathcal{L}_{1,0}(H)^{\otimes g} \otimes \mathcal{L}_{0,1}(H)^{\otimes n}\right)$
Indeed, let
\[ \overset{I}{C}{_{g,n}^{(\pm)}} = \alpha_{g,n}^{-1}\!\left(\overset{I}{\underline{C}}{^{(\pm)}}(1) \ldots \overset{I}{\underline{C}}{^{(\pm)}}(g) \overset{I}{\underline{M}}{^{(\pm)}}(g+1) \ldots \overset{I}{\underline{M}}{^{(\pm)}}(g+n)\right) \in \mathrm{Mat}_{\dim(I)}\!\left(\mathcal{L}_{g,n}(H)\right). \]
\begin{proposition}
The following equality holds in $\mathcal{L}_{g,n}(H)$:
\[ \overset{I}{C}_{g,n} = \overset{I}{C}{^{(+)}_{g,n}}\,\overset{I}{C}{^{(-)-1}_{g,n}}. \]
Moreover, the matrices $\overset{I}{C}_{g,n}$ satisfy the fusion relation of $\mathcal{L}_{0,1}(H)$:
\[ (\overset{I \otimes J}{C}\!\!_{g,n})_{12} = (\overset{I}{C}_{g,n})_1 \, \overset{IJ}{(R')}_{12} \, (\overset{J}{C}_{g,n})_2 \, \overset{IJ}{(R'^{-1})}_{12}. \]
\end{proposition}
\debutDemo
The first claim is a simple consequence of the definition of $\alpha_{g,n}$ and of Lemma \ref{decGauss}. The fusion relation is a consequence of a more general fact which is easy to show, namely: if $i_1 < \ldots < i_k$ and if $\overset{I}{X^1}(i_1), \ldots, \overset{I}{X^k}(i_k)$ are matrices satisfying the fusion relation of $\mathcal{L}_{0,1}(H)$, then their product $\overset{I}{X^1}(i_1) \ldots \overset{I}{X^k}(i_k)$ also satisfies the fusion relation of $\mathcal{L}_{0,1}(H)$.
\finDemo

\indent The image of these matrices have simple expressions in $\mathcal{H}(\mathcal{O}(H))^{\otimes g} \otimes H^{\otimes n}$:
\begin{lemme}\label{expressionM}
It holds
\begin{align*}
\Psi_{g,n}( \overset{I}{C}{_{g,n}^{(+)}} ) &= \overset{I}{a_i}\:\widetilde{b_i^{(2g -1 + n)}} b_i^{(2g + n)} \otimes \ldots \otimes \widetilde{b_i^{(1+n)}}b_i^{(2+n)} \otimes b_i^{(n)} \otimes \ldots \otimes b^{(1)}_i\\
\Psi_{g,n}( \overset{I}{C}{_{g,n}^{(-)}} ) &= \overset{I}{S^{-1}(b_i)} \: \widetilde{a_i^{(2g -1 + n)}} a_i^{(2g + n)} \otimes \ldots \otimes \widetilde{a_i^{(1+n)}}a_i^{(2+n)} \otimes a_i^{(n)} \otimes \ldots \otimes a^{(1)}_i\\
\Psi_{g,n}( \overset{I}{C}_{g,n} ) &= \overset{I}{X_i}\: \widetilde{Y_i^{(2g -1 + n)}} Y_i^{(2g + n)} \otimes \ldots \otimes \widetilde{Y_i^{(1+n)}}Y_i^{(2+n)} \otimes Y_i^{(n)} \otimes \ldots \otimes Y^{(1)}_i
\end{align*}
where $X_i \otimes Y_i = RR'$ and the superscripts mean iterated coproduct.
\end{lemme}
\debutDemo
As an immediate consequence of quasitriangularity, we have for all $n \geq 2$ 
\[ (\mathrm{id} \otimes \Delta^{(n-1)})(R) = a_i \otimes b^{(1)}_i \otimes \ldots \otimes b^{(n)}_i = a_{i_1} \ldots a_{i_n} \otimes b_{i_n} \otimes \ldots \otimes b_{i_1}. \]
with implicit summation on $i_1, \ldots, i_n $. It follows that
\begin{align*}
\Psi_{g,n}( \overset{I}{C}{_{g,n}^{(+)}} ) &= \overset{I}{L}{^{(+)}}(1)\overset{I}{\widetilde{L}}{^{(+)}}(1) \ldots \overset{I}{L}{^{(+)}}(g)\overset{I}{\widetilde{L}}{^{(+)}}(g) \overset{I}{L}{^{(+)}}(g+1) \ldots \overset{I}{L}{^{(+)}}(g+n)\\
&= \overset{I}{a_{i_1}} \ldots \overset{I}{a_{i_{2g+n}}} \: \widetilde{b_{i_2}} b_{i_1} \otimes \ldots \otimes \widetilde{b_{i_{2g}}} b_{i_{2g-1}} \otimes b_{i_{2g+1}} \otimes \ldots \otimes b_{i_{2g+n}}\\
&= \overset{I}{a_i}\:\widetilde{b_i^{(2g -1 + n)}} b_i^{(2g + n)} \otimes \ldots \otimes \widetilde{b_i^{(1+n)}}b_i^{(2+n)} \otimes b_i^{(n)} \otimes \ldots \otimes b^{(1)}_i
\end{align*}
as desired. The second is shown similarly since $R'^{-1}$ is also an universal $R$-matrix. The third is an immediate consequence.
\finDemo
\smallskip\\
\indent The matrices $\overset{I}{C}_{g,n}$ satisfying the fusion relation of $\mathcal{L}_{0,1}(H)$, we can apply Lemma \ref{injectionFusion} and define a representation of $H$ on $V = (H^*)^{\otimes g} \otimes I_1 \otimes \ldots \otimes I_n$ by
\begin{equation}\label{actionHsurFormes}
h \cdot v = h_{C_{g,n}} \triangleright v.
\end{equation}
Since $H$ is factorizable, each $h \in H$ is a linear combination of coefficients of the matrices $\overset{I}{X}_i Y_i$. Hence, $h_{C_{g,n}}$ is a linear combination of coefficients of the matrices $\overset{I}{C}_{g,n}$. It follows from Lemma \ref{expressionM}  that this representation is explicitly given by
\begin{equation}\label{actionH}
\begin{split}
&h \cdot \varphi_1 \, \otimes \, \ldots \, \otimes \, \varphi_g \, \otimes v_1 \otimes \, \ldots \, \otimes v_n\\
&= \varphi_1\!\left(S^{-1}\!\left(h^{(2g -1 + n)}\right) ? h^{(2g + n)}\right)\,  \otimes \ldots \otimes \, \varphi_g\!\left(S^{-1}\!\left(h^{(1+n)}\right) ? h^{(2+n)}\right)
\otimes \,h^{(n)}v_1 \otimes \ldots \otimes \, h^{(1)}v_n. 
\end{split}
\end{equation}

\subsection{Determination and representation of $\mathcal{L}_{g,n}^{\mathrm{inv}}(H)$}
\indent The matrices $\overset{I}{C}_{g,n}$ introduced above allow one to give a simple characterization of the invariant elements of $\mathcal{L}_{g,n}^{\mathrm{inv}}(H)$ and to construct representations of them. We begin with a technical lemma.
\begin{lemme}\label{conjugaisonM}
It holds
\[ (\overset{I}{C}{^{(\pm)}_{g,n}})_1 \, \overset{J}{U}(i)_2 \, (\overset{I}{C}{^{(\pm)}_{g,n}})^{-1}_1 \:\:= \overset{IJ}{R}{^{(\pm)-1}_{12}} \, \overset{J}{U}(i)_2 \, \overset{IJ}{R}{^{(\pm)}_{12}} \]
where $U$ is $A$ or $B$.
\end{lemme}
\debutDemo
Applying the isomorphisms $\Psi_{0,1}$ and $\Psi_{1,0}$ and using relations \eqref{propertiesL}, \eqref{echangeHeisenberg} and \eqref{LTilde}, it is easy to show the result for $\mathcal{L}_{0,1}(H)$ and $\mathcal{L}_{1,0}(H)$. We get similarly:
\[ \overset{IJ}{R}{^{(\pm)}_{12}} \, \overset{I}{C}{^{(\pm)}_{1}} \, \overset{J}{C}{^{(-)}_{2}} = \overset{J}{C}{^{(-)}_{2}} \, \overset{I}{C}{^{(\pm)}_{1}} \, \overset{IJ}{R}{^{(\pm)}_{12}}, 
\:\:\:\: \overset{IJ}{R}{^{(\pm)}_{12}} \, \overset{I}{M}{^{(\pm)}_{1}} \, \overset{J}{M}{^{(-)}_{2}} = \overset{J}{M}{^{(-)}_{2}} \, \overset{I}{M}{^{(\pm)}_{1}} \, \overset{IJ}{R}{^{(\pm)}_{12}}. \]
Using these preliminary facts, we can carry out the general computation. For instance, for $i \leq g$ 
{\footnotesize\begin{align*}
&\alpha_{g,n}\!\left( (\overset{I}{C}{^{(\pm)}_{g,n}})_1 \, \overset{J}{U}(i)_2 \, (\overset{I}{C}{^{(\pm)}_{g,n}})^{-1}_1\:\, \right) \\
& = \overset{I}{\underline{C}}{^{(\pm)}}(1)_1 \, \ldots \, \overset{I}{\underline{C}}{^{(\pm)}}(i)_1 \, \overset{J}{\underline{C}}{^{(-)}}(1)_2 \, \ldots \, \overset{J}{\underline{C}}{^{(-)}}(i-1)_2 \,
\overset{J}{\underline{U}}(i)_2 \, \overset{J}{\underline{C}}{^{(-)}}(i-1)^{-1}_2 \, 
 \ldots \, \overset{J}{\underline{C}}{^{(-)}}(1)^{-1}_2 \, \overset{I}{\underline{C}}{^{(\pm)}}(i)^{-1}_1 \, \ldots \, \overset{I}{\underline{C}}{^{(\pm)}}(1)^{-1}_1\\
& = \overset{I}{\underline{C}}{^{(\pm)}}(1)_1 \, \overset{J}{\underline{C}}{^{(-)}}(1)_2 \, \ldots \, \overset{I}{\underline{C}}{^{(\pm)}}(i-1)_1 \, \overset{J}{\underline{C}}{^{(-)}}(i-1)_2 \, 
\overset{I}{\underline{C}}{^{(\pm)}}(i)_1 \, \overset{J}{\underline{U}}(i)_2 \, 
 \overset{I}{\underline{C}}{^{(\pm)}}(i)^{-1}_1 \, \overset{J}{\underline{C}}{^{(-)}}(i-1)^{-1}_2 \, \overset{I}{\underline{C}}{^{(\pm)}}(i-1)^{-1}_1 \, \ldots\\
& \:\:\:\:\:\: \overset{J}{\underline{C}}{^{(-)}}(1)^{-1}_2 \, \overset{I}{\underline{C}}{^{(\pm)}}(1)^{-1}_1\\
& = \overset{I}{\underline{C}}{^{(\pm)}}(1)_1 \, \overset{J}{\underline{C}}{^{(-)}}(1)_2 \, \ldots \overset{I}{\underline{C}}{^{(\pm)}}(i-1)_1 \, \overset{J}{\underline{C}}{^{(-)}}(i-1)_2 \, 
\overset{IJ}{R}{^{(\pm)-1}_{12}} \, \overset{J}{\underline{U}}(i)_2 \, 
\overset{IJ}{R}{^{(\pm)}_{12}} \, \overset{J}{\underline{C}}{^{(-)}}(i-1)^{-1}_2 \, \overset{I}{\underline{C}}{^{(\pm)}}(i-1)^{-1}_1 \, \ldots \\
& \:\:\:\:\:\: \overset{J}{\underline{C}}{^{(-)}}(1)^{-1}_2 \, \overset{I}{\underline{C}}{^{(\pm)}}(1)^{-1}_1\\ 
& =\overset{IJ}{R}{^{(\pm)-1}_{12}}  \, \overset{J}{\underline{C}}{^{(-)}}(1)_2 \, \overset{I}{\underline{C}}{^{(\pm)}}(1)_1 \, \ldots \overset{J}{\underline{C}}{^{(-)}}(i-1)_2 \, 
\overset{I}{\underline{C}}{^{(\pm)}}(i-1)_1 \, \overset{J}{\underline{U}}(i)_2 \, 
 \overset{I}{\underline{C}}{^{(\pm)}}(i-1)^{-1}_1 \, \overset{J}{\underline{C}}{^{(-)}}(i-1)^{-1}_2 \, \ldots \\
& \:\:\:\:\:\: \overset{I}{\underline{C}}{^{(\pm)}}(1)^{-1}_1 \, \overset{J}{\underline{C}}{^{(-)}}(1)^{-1}_2 \, \overset{IJ}{R}{^{(\pm)}_{12}}\\ 
& =\overset{IJ}{R}{^{(\pm)-1}_{12}} \, \overset{J}{\underline{C}}{^{(-)}}(1)_2 \, \ldots \overset{J}{\underline{C}}{^{(-)}}(i-1)_2 \,
 \overset{J}{\underline{U}}(i)_2 \,  \overset{J}{\underline{C}}{^{(-)}}(i-1)^{-1}_2 \, \ldots \overset{J}{\underline{C}}{^{(-)}}(1)^{-1}_2 \, \overset{IJ}{R}{^{(\pm)}_{12}}
= \alpha_{g,n}\!\left(  \overset{IJ}{R}{^{(\pm)-1}_{12}} \, \overset{J}{U}(i)_2 \, \overset{IJ}{R}{^{(\pm)}_{12}} \right).
\end{align*}}
The case $i > g$ is treated in a similar way.
\finDemo
\smallskip\\
\indent For $(V, \triangleright)$ a representation of $\mathcal{L}_{g,n}(H)$, let
\[ \mathrm{Inv}(V) = \left\{ v \in V \, \left| \, \forall \, I, \: \overset{I}{C}_{g,n} \triangleright v = \mathbb{I}_{\dim(I)}v \right.\right\} = \left\{ v \in V \, \left| \, \forall \, h \in H, \:\: h \cdot v = \varepsilon(h)v \right.\right\} \]
where $\mathbb{I}_k$ is the identity matrix of size $k$, and the action $\cdot$ of $H$ on $V$ is defined in \eqref{actionHsurFormes} and \eqref{actionH}.

\begin{theoreme}\label{thmInv}
1) An element $x \in \mathcal{L}_{g,n}(H)$ is invariant under the action of $H$ (or equivalently under the coaction $\Omega$ of $\mathcal{O}(H)$) if, and only if, for every $H$-module $I$, $\overset{I}{C}_{g,n}x = x\overset{I}{C}_{g,n}$. \\
2) Let $V$ be a representation of $\mathcal{L}_{g,n}(H)$. Then $\mathrm{Inv}(V)$ is stable under the action of invariant elements and thus provides a representation of $\mathcal{L}^{\mathrm{inv}}_{g,n}(H)$.
\end{theoreme}
\debutDemo
1) The right action of $H$ on $\mathcal{L}_{g,n}(H)$ is by definition 
\[ \overset{I}{U}(i) \cdot h = \overset{I}{h'} \overset{I}{U}(i) \overset{I}{S(h'')}. \]
where $U(i)$ is $A(i)$ or $B(i)$ if $i \leq g$ and is $M(i)$ if $i > g$. Then, denoting $a^{(\pm)}_j \otimes b^{(\pm)}_j = R^{(\pm)}$, 
\begin{equation*}
\begin{split}
\overset{J}{U}(i)_2 \, \cdot \, S^{-1}(\overset{I}{L}{^{(\pm)}})_1 &= \overset{J}{U}(i)_2 \, \cdot \, S^{-1}(b^{(\pm)}_j) \, \overset{I}{(a_j^{(\pm)})}_1 = \overset{J}{S^{-1}(b_j^{(\pm)})}_2 \overset{J}{U}(i)_2 \overset{J}{(b^{(\pm)}_k)}_2 (\overset{I}{a_j^{(\pm)}} \overset{I}{a_k^{(\pm)}})_1 \\ &= \overset{IJ}{R}{^{(\pm)-1}_{12}} \overset{J}{U}(i)_2 \overset{IJ}{R}{^{(\pm)}_{12}} = (\overset{I}{C}{^{(\pm)}_{g,n}})_1 \, \overset{J}{U}(i)_2 \, (\overset{I}{C}{^{(\pm)}_{g,n}})_1^{-1}.
\end{split}
\end{equation*}
where the last equality is Lemma \ref{conjugaisonM}. Observe that the matrix $\overset{J}{U}(i)_2 \, \cdot \, S^{-1}(\overset{I}{L}{^{(\pm)}})_1$ contains all the elements obtained by acting by the coefficients of $S^{-1}(\overset{I}{L}{^{(\pm)}})$ on the coefficients of $\overset{J}{U}(i)$. The coefficients of the $S^{-1}(\overset{I}{L}{^{(\pm)}}) = \overset{I}{L}{^{(\pm)-1}}$ generate $H$ as an algebra, hence we deduce that an $x \in \mathcal{L}_{g,n}(H)$ is invariant if, and only if, $\overset{I}{C}{^{(\pm)}_{g,n}}x = x\overset{I}{C}{^{(\pm)}_{g,n}}$. Since $H$ is factorizable, it is generated by the coefficients of the matrices $\overset{I}{X}_i \otimes Y_i$ and we see by Lemma \ref{expressionM} that the subalgebra of $\mathcal{L}_{g,n}(H)$ generated by the coefficients of the matrices $\overset{I}{C}{^{(\pm)}_{g,n}}$ equals the subalgebra generated by the coefficients of the matrices $\overset{I}{C}_{g,n}$. Thus $x$ commutes with the $\overset{I}{C}{^{(\pm)}_{g,n}}$ if, and only if, $x$ commutes with the $\overset{I}{C}_{g,n}$.
\smallskip\\
2) If $v \in \mathrm{Inv}(V)$ and $x \in \mathcal{L}_{g,n}^{\mathrm{inv}}(H)$, we have
\[ \overset{I}{C}_{g,n} \triangleright (x \triangleright v) = \overset{I}{C}_{g,n} x \triangleright v = x \overset{I}{C}_{g,n} \triangleright v = \mathbb{I}_{\dim(I)} (x \triangleright v) \]
which shows that $x \triangleright v \in \mathrm{Inv}(V)$.
\finDemo
%\smallskip\\
%\indent Let $V = (H^*)^{\otimes g} \otimes I_1 \otimes \ldots \otimes I_n$ be a representation of $\mathcal{L}_{g,n}(H)$. Recall that by \eqref{actionH}, an element $\sum_i \varphi_{1,i} \otimes \ldots \otimes \varphi_{g,i} \otimes v_{1,i} \otimes \ldots \otimes v_{n,i} \in V$ belongs to $\mathrm{Inv}(V)$ if, and only if, for all $h \in H$, 
%\begin{equation*}
%\begin{split}
%&\sum_i \varphi_{1,i}\!\left(S^{-1}\!\left(h^{(2g -1 + n)}\right) ? h^{(2g + n)}\right)\,  \otimes \ldots \otimes \, \varphi_{g,i}\!\left(S^{-1}\!\left(h^{(1+n)}\right) ? h^{(2+n)}\right)
%\otimes \,h^{(n)}v_{1,i} \otimes \ldots \otimes \, h^{(1)}v_{n,i}\\
%&= \varepsilon(h) \sum_i \varphi_{1,i} \otimes \ldots \otimes \varphi_{g,i} \otimes v_{1,i} \otimes \ldots \otimes v_{n,i}. 
%\end{split}
%\end{equation*}

By definition, a $H$-connection $\nabla = (h_e)_{e \in E}$ (with $E = \{ b_1, a_1, \ldots, b_g, a_g, m_{g+1}, \ldots, m_{g+n} \}$, the set of edges) is flat if its holonomy along the boundary $c = b_1 a_1^{-1} b_1^{-1} a_1 \ldots b_g a_g^{-1} b_g^{-1} a_g m_{g+1} \ldots m_{g+n}$ of the unique face of the graph $\Gamma$ is trivial:
\[ \mathrm{Hol}^{\nabla}(C) = h_{b_1} h_{a_1}^{-1} h_{b_1}^{-1} h_{a_1} \ldots h_{b_g} h_{a_g}^{-1} h_{b_g}^{-1} h_{a_g} h_{m_{g+1}} \ldots h_{m_{g+n}} = 1.\]
Hence, the subrepresentation $\mathrm{Inv}(V)$ implements this flatness constraint. This constraint was directly implemented in $\mathcal{L}_{g,n}(H)$  (and not just on representations) in \cite{AGS2, AS} by means of characteristic projectors, giving rise to the moduli algebra, a quantum analogue of $\mathbb{C}[\mathcal{A}_f/\mathcal{G}]$ (see Introduction). However, the definition of these projectors requires the $S$-matrix, which has nice properties in the semi-simple case only. In \cite{AS}, the representation space of the mapping class group is the moduli algebra. Here we do not consider the moduli algebra; in particular we will not need to generalize these projectors to construct the projective representation of the mapping class group.

\section{Projective representations of mapping class groups}
\indent Let $\Sigma_{g,n}$ be the compact orientable surface of genus $g$ with $n$ open disks removed. For simplicity we consider the case of $\Sigma_g$ ($n = 0$). The particular features in this case are that the presentation of the mapping class group is easier and that the associated algebra $\mathcal{L}_{g,0}(H) \cong \mathcal{H}(\mathcal{O}(H))^{\otimes g}$ is isomorphic to a matrix algebra. 
\smallskip\\
\indent We will discuss the case of $n>0$ in subsection \ref{CasGeneral}.

\subsection{Mapping class group of $\Sigma_g$}\label{sectionRepMCG}
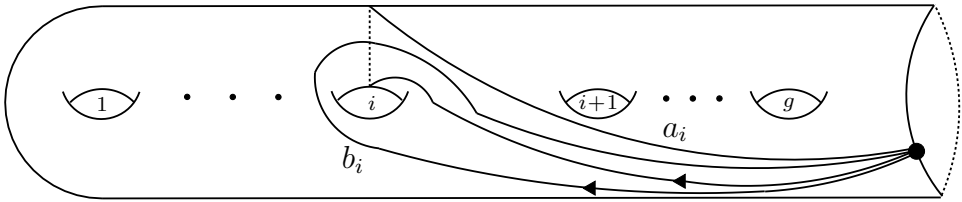
\begin{figure}[h]
\centering
\begin{tikzpicture}
\draw [shift={(6.469564580924096,13.159138980991877)},line width=0.7pt]  plot[domain=3.4622754799227886:5.956508902908017,variable=\t]({1*0.5320950869782717*cos(\t r)+0*0.5320950869782717*sin(\t r)},{0*0.5320950869782717*cos(\t r)+1*0.5320950869782717*sin(\t r)});
\draw [shift={(10.00269783209018,13.165039608452211)},line width=0.7pt]  plot[domain=3.4622754799227886:5.956508902908016,variable=\t]({1*0.5320950869782717*cos(\t r)+0*0.5320950869782717*sin(\t r)},{0*0.5320950869782717*cos(\t r)+1*0.5320950869782717*sin(\t r)});
\draw [shift={(10.00269783208794,12.464688987415235)},line width=0.7pt]  plot[domain=0.7739312631124724:2.3657977977209477,variable=\t]({1*0.5910130783015503*cos(\t r)+0*0.5910130783015503*sin(\t r)},{0*0.5910130783015503*cos(\t r)+1*0.5910130783015503*sin(\t r)});
\draw [shift={(13.002240809249017,13.168532653273042)},line width=0.7pt]  plot[domain=3.4622754799227886:5.956508902908016,variable=\t]({1*0.5320950869782717*cos(\t r)+0*0.5320950869782717*sin(\t r)},{0*0.5320950869782717*cos(\t r)+1*0.5320950869782717*sin(\t r)});
\draw [shift={(15.513845165580465,13.163431268478787)},line width=0.7pt]  plot[domain=3.4622754799227886:5.956508902908013,variable=\t]({1*0.5320950869782717*cos(\t r)+0*0.5320950869782717*sin(\t r)},{0*0.5320950869782717*cos(\t r)+1*0.5320950869782717*sin(\t r)});
\draw [shift={(6.482806603614594,12.43395118907535)},line width=0.7pt]  plot[domain=0.7739312631124735:2.3657977977209477,variable=\t]({1*0.5910130783015496*cos(\t r)+0*0.5910130783015496*sin(\t r)},{0*0.5910130783015496*cos(\t r)+1*0.5910130783015496*sin(\t r)});
\draw [shift={(12.999224683787697,12.431668913407975)},line width=0.7pt]  plot[domain=0.7739312631124746:2.3657977977209477,variable=\t]({1*0.5910130783015489*cos(\t r)+0*0.5910130783015489*sin(\t r)},{0*0.5910130783015489*cos(\t r)+1*0.5910130783015489*sin(\t r)});
\draw [shift={(15.51512446620466,12.429566421432781)},line width=0.7pt]  plot[domain=0.773931263112466:2.3657977977209517,variable=\t]({1*0.5910130783015439*cos(\t r)+0*0.5910130783015439*sin(\t r)},{0*0.5910130783015439*cos(\t r)+1*0.5910130783015439*sin(\t r)});
\draw [shift={(6.480309438845568,12.848959304762982)},line width=0.7pt]  plot[domain=1.5707963267948966:4.71238898038469,variable=\t]({1*1.2742709713406324*cos(\t r)+0*1.2742709713406324*sin(\t r)},{0*1.2742709713406324*cos(\t r)+1*1.2742709713406324*sin(\t r)});
\draw [line width=0.7pt] (6.480309438845568,14.123230276103614)-- (17.422556685283112,14.135909589251282);
\draw [line width=0.7pt] (6.480309438845568,11.574688333422351)-- (17.51131187731679,11.58736764657002);
\draw [shift={(19.184981212808975,12.944054153370493)},line width=0.7pt]  plot[domain=2.546978576128653:3.822767330570879,variable=\t]({1*2.127594791134662*cos(\t r)+0*2.127594791134662*sin(\t r)},{0*2.127594791134662*cos(\t r)+1*2.127594791134662*sin(\t r)});
\draw [shift={(14.721862984829814,12.779223082450809)},line width=0.7pt,dash pattern=on 1pt off 1pt]  plot[domain=-0.39351986934882:0.4655237794072971,variable=\t]({1*3.03225953065302*cos(\t r)+0*3.03225953065302*sin(\t r)},{0*3.03225953065302*cos(\t r)+1*3.03225953065302*sin(\t r)});
\draw [shift={(14.871897980950925,17.074816802592615)},line width=0.7pt]  plot[domain=4.773165963111252:5.156552164893947,variable=\t]({1*5.418821348503365*cos(\t r)+0*5.418821348503365*sin(\t r)},{0*5.418821348503365*cos(\t r)+1*5.418821348503365*sin(\t r)});
\draw [shift={(14.341129164738266,27.013490072683705)},line width=0.7pt]  plot[domain=4.434380004142646:4.768359539998192,variable=\t]({1*15.36871747587277*cos(\t r)+0*15.36871747587277*sin(\t r)},{0*15.36871747587277*cos(\t r)+1*15.36871747587277*sin(\t r)});
\draw [shift={(10.225482873570378,13.179940429435394)},line width=0.7pt]  plot[domain=3.078612991156848:4.604699109551869,variable=\t]({1*0.9469443022116039*cos(\t r)+0*0.9469443022116039*sin(\t r)},{0*0.9469443022116039*cos(\t r)+1*0.9469443022116039*sin(\t r)});
\draw [shift={(9.910460565236708,12.958573401957683)},line width=0.7pt]  plot[domain=1.3347726481689977:2.722115191834436,variable=\t]({1*0.6917948831092703*cos(\t r)+0*0.6917948831092703*sin(\t r)},{0*0.6917948831092703*cos(\t r)+1*0.6917948831092703*sin(\t r)});
\draw [shift={(9.757206469290598,11.749568867271712)},line width=0.7pt]  plot[domain=0.5235757641220097:1.4049139737965757,variable=\t]({1*1.9072381002667376*cos(\t r)+0*1.9072381002667376*sin(\t r)},{0*1.9072381002667376*cos(\t r)+1*1.9072381002667376*sin(\t r)});
\draw [shift={(15.08193507989958,21.619561834284326)},line width=0.7pt]  plot[domain=4.321635583891609:4.932736147087521,variable=\t]({1*9.643301163013478*cos(\t r)+0*9.643301163013478*sin(\t r)},{0*9.643301163013478*cos(\t r)+1*9.643301163013478*sin(\t r)});
\draw [line width=0.7pt,dash pattern=on 1pt off 1pt] (10.000166143032205,13.055696643265817)-- (10.000128180926128,14.124443941377749);
\draw [shift={(14.955133874245728,17.24602134616958)},line width=0.7pt]  plot[domain=4.555917771113938:5.1296043544394205,variable=\t]({1*5.505266920503724*cos(\t r)+0*5.505266920503724*sin(\t r)},{0*5.505266920503724*cos(\t r)+1*5.505266920503724*sin(\t r)});
\draw [shift={(14.775250561197215,19.55683005994665)},line width=0.7pt]  plot[domain=4.181115566542873:4.625111267435227,variable=\t]({1*7.783950677704752*cos(\t r)+0*7.783950677704752*sin(\t r)},{0*7.783950677704752*cos(\t r)+1*7.783950677704752*sin(\t r)});
\draw [shift={(10.333516446691569,12.638241096542256)},line width=0.7pt]  plot[domain=0.3947911196997582:2.244639591154181,variable=\t]({1*0.539649939145546*cos(\t r)+0*0.539649939145546*sin(\t r)},{0*0.539649939145546*cos(\t r)+1*0.539649939145546*sin(\t r)});
\draw [shift={(15.56330898184788,20.711407038664966)},line width=0.7pt]  plot[domain=4.010397061863714:4.901431663365731,variable=\t]({1*8.622932082711552*cos(\t r)+0*8.622932082711552*sin(\t r)},{0*8.622932082711552*cos(\t r)+1*8.622932082711552*sin(\t r)});
\draw[color=black] (9.780922284021113,12.0896509555724) node {$b_i$};
\draw[color=black] (14.016964707628347,12.453485939830028) node {$a_i$};
\begin{scriptsize}
\draw[color=black] (6.481430737229788,12.826529098498881) node {$1$};
\draw[color=black] (10.0304404379994802,12.837681256432324) node {$i$};
\draw[color=black] (13.007378022759817,12.837681256432324) node {$i\!+\!1$};
\draw[color=black] (15.513570399115622,12.826529098498881) node {$g$};
\draw [fill=black] (17.191817204379586,12.199771390729145) circle (3pt);
\draw [fill=black,shift={(14.096747319484662,11.802507297517565)},rotate=90] (0,0) ++(0 pt,3pt) -- ++(2.598076211353316pt,-4.5pt)--++(-5.196152422706632pt,0 pt) -- ++(2.598076211353316pt,4.5pt);
\draw [fill=black,shift={(12.912235568100115,11.711341859354846)},rotate=90] (0,0) ++(0 pt,3pt) -- ++(2.598076211353316pt,-4.5pt)--++(-5.196152422706632pt,0 pt) -- ++(2.598076211353316pt,4.5pt);
\draw [fill=black] (13.9,12.9) circle (1pt);
\draw [fill=black] (14.6,12.9) circle (1pt);
\draw [fill=black] (14.25,12.9) circle (1pt);
\draw [fill=black] (7.599886378457062,12.92127233955011) circle (1pt);
\draw [fill=black] (8.79930680637218,12.91904706973951) circle (1pt);
\draw [fill=black] (8.19959659241462,12.92015970464481) circle (1pt);
\end{scriptsize}
\end{tikzpicture}
\caption{Surface $\Sigma_g^{\mathrm{o}}$ and generators of $\pi_1(\Sigma_g^{\mathrm{o}})$.}
\label{figureSurface}
\end{figure}

\indent Let $D \subset \Sigma_g$ be an embedded open disk and define $\Sigma_g^{\mathrm{o}} = \Sigma_g \setminus D$. We put a basepoint on the boundary circle $c = \partial(\Sigma_g^{\mathrm{o}})$ and we take the curves $a_i, b_i$ ($1 \leq i \leq g$) represented in Figure \ref{figureSurface} as generators for the free group $\pi_1(\Sigma_g^{\mathrm{o}})$. 

%%%%%%%%%%%%%%%%%%%%%%%%%%%%%%%%%%%%%%%%%%%%%%%%%%%%%%
%\begin{comment}

%\end{comment}
%%%%%%%%%%%%%%%%%%%%%%%%%%%%%%%%%%%%%%%%%%%%%%%%%%%%%%

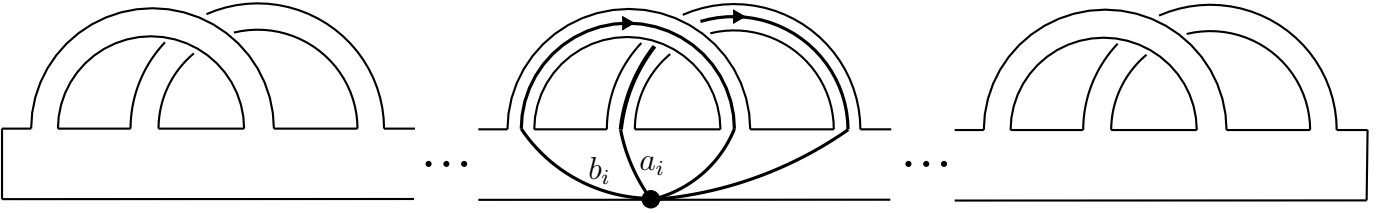
\begin{figure}[h]
\begin{tikzpicture}[scale=0.935]
\draw [shift={(7.116832178429767,6)},line width=0.8pt]  plot[domain=0:3.141592653589793,variable=\t]({1*1.7075037526311627*cos(\t r)+0*1.7075037526311627*sin(\t r)},{0*1.7075037526311627*cos(\t r)+1*1.7075037526311627*sin(\t r)});
\draw [shift={(7.095787711044769,6)},line width=0.8pt]  plot[domain=0:3.141592653589793,variable=\t]({1*1.309575195744216*cos(\t r)+0*1.309575195744216*sin(\t r)},{0*1.309575195744216*cos(\t r)+1*1.309575195744216*sin(\t r)});
\draw [shift={(8.596936063235535,6)},line width=0.8pt]  plot[domain=0:1.8118338027760237,variable=\t]({1*1.4030639367644646*cos(\t r)+0*1.4030639367644646*sin(\t r)},{0*1.4030639367644646*cos(\t r)+1*1.4030639367644646*sin(\t r)});
\draw [shift={(8.596936063235535,6)},line width=0.8pt]  plot[domain=0:1.9894438844649243,variable=\t]({1*1.776642756206547*cos(\t r)+0*1.776642756206547*sin(\t r)},{0*1.776642756206547*cos(\t r)+1*1.776642756206547*sin(\t r)});
\draw [shift={(8.596936063235535,6)},line width=0.8pt]  plot[domain=2.2704846801922587:3.141592653589793,variable=\t]({1*1.3953711574291767*cos(\t r)+0*1.3953711574291767*sin(\t r)},{0*1.3953711574291767*cos(\t r)+1*1.3953711574291767*sin(\t r)});
\draw [shift={(8.596936063235535,6)},line width=0.8pt]  plot[domain=2.3877511379579315:3.141592653589793,variable=\t]({1*1.789237366642512*cos(\t r)+0*1.789237366642512*sin(\t r)},{0*1.789237366642512*cos(\t r)+1*1.789237366642512*sin(\t r)});
\draw [line width=0.8pt] (5,6)-- (5,5);
\draw [line width=0.8pt] (5,5)-- (10.794708954246598,5.00203777401749);
\draw [line width=0.8pt] (5,6)-- (5.409328425798604,6);
\draw [line width=0.8pt] (5.786212515300553,6)-- (6.8080282517522,6);
\draw [line width=0.8pt] (7.200815462007147,6)-- (8.405362906788985,6);
\draw [line width=0.8pt] (8.82433593106093,6)-- (10,6);
\draw [line width=0.8pt] (10.373578819442082,6)-- (10.80786927095924,6);
\draw [shift={(13.816832178429763,5.99)},line width=0.8pt]  plot[domain=0:3.141592653589793,variable=\t]({1*1.7075037526311627*cos(\t r)+0*1.7075037526311627*sin(\t r)},{0*1.7075037526311627*cos(\t r)+1*1.7075037526311627*sin(\t r)});
\draw [shift={(13.795787711044765,5.99)},line width=0.8pt]  plot[domain=0:3.141592653589793,variable=\t]({1*1.309575195744216*cos(\t r)+0*1.309575195744216*sin(\t r)},{0*1.309575195744216*cos(\t r)+1*1.309575195744216*sin(\t r)});
\draw [shift={(15.296936063235531,5.99)},line width=0.8pt]  plot[domain=0:1.8118338027760237,variable=\t]({1*1.4030639367644646*cos(\t r)+0*1.4030639367644646*sin(\t r)},{0*1.4030639367644646*cos(\t r)+1*1.4030639367644646*sin(\t r)});
\draw [shift={(15.296936063235531,5.99)},line width=0.8pt]  plot[domain=0:1.9894438844649243,variable=\t]({1*1.776642756206547*cos(\t r)+0*1.776642756206547*sin(\t r)},{0*1.776642756206547*cos(\t r)+1*1.776642756206547*sin(\t r)});
\draw [shift={(15.296936063235531,5.99)},line width=0.8pt]  plot[domain=2.2704846801922587:3.141592653589793,variable=\t]({1*1.3953711574291767*cos(\t r)+0*1.3953711574291767*sin(\t r)},{0*1.3953711574291767*cos(\t r)+1*1.3953711574291767*sin(\t r)});
\draw [shift={(15.296936063235531,5.99)},line width=0.8pt]  plot[domain=2.3877511379579315:3.141592653589793,variable=\t]({1*1.789237366642512*cos(\t r)+0*1.789237366642512*sin(\t r)},{0*1.789237366642512*cos(\t r)+1*1.789237366642512*sin(\t r)});
\draw [shift={(14.14958069516861,7.227726613631109)},line width=1.2pt]  plot[domain=3.729595257137362:4.701699287460754,variable=\t]({1*2.2187359355090455*cos(\t r)+0*2.2187359355090455*sin(\t r)},{0*2.2187359355090455*cos(\t r)+1*2.2187359355090455*sin(\t r)});
\draw [shift={(13.802618659443407,5.996993354832272)},line width=1.2pt]  plot[domain=0:3.141592653589793,variable=\t]({1*1.4991378524730479*cos(\t r)+0*1.4991378524730479*sin(\t r)},{0*1.4991378524730479*cos(\t r)+1*1.4991378524730479*sin(\t r)});
\draw [shift={(16.112966893148847,6.443978813579067)},line width=1.2pt]  plot[domain=3.327560085751379:3.7729570768128284,variable=\t]({1*2.455301608888626*cos(\t r)+0*2.455301608888626*sin(\t r)},{0*2.455301608888626*cos(\t r)+1*2.455301608888626*sin(\t r)});
\draw [shift={(16.032112045973285,5.7372208481955465)},line width=1.6pt]  plot[domain=2.484254690839589:3.0336235085501446,variable=\t]({1*2.3423271124786536*cos(\t r)+0*2.3423271124786536*sin(\t r)},{0*2.3423271124786536*cos(\t r)+1*2.3423271124786536*sin(\t r)});
\draw [shift={(13.755663711022823,6.614413362953582)},line width=1.2pt]  plot[domain=4.936459140732812:5.903246179167029,variable=\t]({1*1.6651881573061873*cos(\t r)+0*1.6651881573061873*sin(\t r)},{0*1.6651881573061873*cos(\t r)+1*1.6651881573061873*sin(\t r)});
\draw [shift={(15.290022790798636,5.978921534198891)},line width=1.2pt]  plot[domain=0.006881023462768673:1.8643557340986399,variable=\t]({1*1.6100153249433073*cos(\t r)+0*1.6100153249433073*sin(\t r)},{0*1.6100153249433073*cos(\t r)+1*1.6100153249433073*sin(\t r)});
\draw [shift={(13.724070995031564,10.489990747896217)},line width=1.2pt]  plot[domain=4.785288706191459:5.326972569147958,variable=\t]({1*5.513782387193336*cos(\t r)+0*5.513782387193336*sin(\t r)},{0*5.513782387193336*cos(\t r)+1*5.513782387193336*sin(\t r)});
\draw [line width=0.8pt] (11.7,4.99)-- (17.494708954246594,4.99203777401749);
\draw [line width=0.8pt] (11.7,5.99)-- (12.1093284257986,5.99);
\draw [line width=0.8pt] (12.486212515300549,5.99)-- (13.508028251752195,5.99);
\draw [line width=0.8pt] (13.900815462007142,5.99)-- (15.105362906788981,5.99);
\draw [line width=0.8pt] (15.524335931060925,5.99)-- (16.7,5.99);
\draw [line width=0.8pt] (17.073578819442076,5.99)-- (17.507869270959233,5.99);
\draw [shift={(20.51683217842976,5.98)},line width=0.8pt]  plot[domain=0:3.141592653589793,variable=\t]({1*1.7075037526311627*cos(\t r)+0*1.7075037526311627*sin(\t r)},{0*1.7075037526311627*cos(\t r)+1*1.7075037526311627*sin(\t r)});
\draw [shift={(20.495787711044763,5.98)},line width=0.8pt]  plot[domain=0:3.141592653589793,variable=\t]({1*1.309575195744216*cos(\t r)+0*1.309575195744216*sin(\t r)},{0*1.309575195744216*cos(\t r)+1*1.309575195744216*sin(\t r)});
\draw [shift={(21.996936063235527,5.98)},line width=0.8pt]  plot[domain=0:1.8118338027760237,variable=\t]({1*1.4030639367644646*cos(\t r)+0*1.4030639367644646*sin(\t r)},{0*1.4030639367644646*cos(\t r)+1*1.4030639367644646*sin(\t r)});
\draw [shift={(21.996936063235527,5.98)},line width=0.8pt]  plot[domain=0:1.9894438844649243,variable=\t]({1*1.776642756206547*cos(\t r)+0*1.776642756206547*sin(\t r)},{0*1.776642756206547*cos(\t r)+1*1.776642756206547*sin(\t r)});
\draw [shift={(21.996936063235527,5.98)},line width=0.8pt]  plot[domain=2.2704846801922587:3.141592653589793,variable=\t]({1*1.3953711574291767*cos(\t r)+0*1.3953711574291767*sin(\t r)},{0*1.3953711574291767*cos(\t r)+1*1.3953711574291767*sin(\t r)});
\draw [shift={(21.996936063235527,5.98)},line width=0.8pt]  plot[domain=2.3877511379579315:3.141592653589793,variable=\t]({1*1.789237366642512*cos(\t r)+0*1.789237366642512*sin(\t r)},{0*1.789237366642512*cos(\t r)+1*1.789237366642512*sin(\t r)});
\draw [line width=0.8pt] (18.4,4.98)-- (24.19470895424659,4.98203777401749);
\draw [line width=0.8pt] (18.4,5.98)-- (18.809328425798597,5.98);
\draw [line width=0.8pt] (19.186212515300546,5.98)-- (20.20802825175219,5.98);
\draw [line width=0.8pt] (20.600815462007137,5.98)-- (21.80536290678898,5.98);
\draw [line width=0.8pt] (22.22433593106092,5.98)-- (23.4,5.98);
\draw [line width=0.8pt] (23.773578819442072,5.98)-- (24.20786927095923,5.98);
\draw [line width=0.8pt] (24.20786927095923,5.98)-- (24.19470895424659,4.98203777401749);
%\begin{scriptsize}
\draw[color=black] (13.406499243040994,5.405300111476876) node {$b_i$};
\draw[color=black] (14.146452087701703,5.480060431237907) node {$a_i$};
\draw [fill=black,shift={(13.783362211564688,7.496007527524644)},rotate=270] (0,0) ++(0 pt,3pt) -- ++(2.598076211353316pt,-4.5pt)--++(-5.196152422706632pt,0 pt) -- ++(2.598076211353316pt,4.5pt);
\draw [fill=black,shift={(15.343159148919554,7.588059773753824)},rotate=270] (0,0) ++(0 pt,3pt) -- ++(2.598076211353316pt,-4.5pt)--++(-5.196152422706632pt,0 pt) -- ++(2.598076211353316pt,4.5pt);
\draw [fill=black] (14.126199467643557,4.998184518464418) circle (3.5pt);
\draw [fill=black] (17.75,5.5) circle (1pt);
\draw [fill=black] (18.25,5.5) circle (1pt);
\draw [fill=black] (18,5.5) circle (1pt);
\draw [fill=black] (11,5.5) circle (1pt);
\draw [fill=black] (11.5,5.5) circle (1pt);
\draw [fill=black] (11.25,5.5) circle (1pt);
%\end{scriptsize}
\end{tikzpicture}
\caption{Surface $\Sigma_g^{\mathrm{o}}$ viewed as a ribbon graph.}
\label{figureSurfaceRuban}
\end{figure}

With these generators, 
\[ c = b_1 a_1^{-1} b_1^{-1} a_1 \ldots b_g a_g^{-1} b_g^{-1} a_g. \]
Recall from the Introduction the embedded oriented graph $\Gamma = \left(\{\bullet\}, \{b_1, a_1, \ldots, b_g, a_g\}\right)$ with vertex $\bullet$ and edges the generators of the fundamental group. It is readily seen that $\Sigma_g^{\mathrm{o}}$ is homeomorphic to the thickening of $\Gamma$ represented in Figure \ref{figureSurfaceRuban}.

%%%%%%%%%%%%%%%%%%%%%%%%%%%%%%%%%%%%%%%%%%%%%%%%%%%%%%
%\begin{comment}

%\end{comment}
%%%%%%%%%%%%%%%%%%%%%%%%%%%%%%%%%%%%%%%%%%%%%%%%%%%%%%

Simple closed curves on a surface will simply be called circles. Elements of $\pi_1(\Sigma_g^{\mathrm{o}})$ will be called loops. We consider circles up to free homotopy. In particular, if $\gamma \in \pi_1(\Sigma_g^{\mathrm{o}})$, we denote by $[\gamma]$ the free homotopy class of $\gamma$. For $\alpha$ a circle, we denote by $\tau_{\alpha}$ the Dehn twist about it. If $\gamma \in \pi_1(\Sigma_g^{\mathrm{o}})$, then $\tau_{\gamma}$ is a shortand for $\tau_{[\gamma]}$, thus defined as follows: consider a circle $\gamma'$ freely homotopic to $\gamma$ and which does not intersect the boundary circle $c$; then $\tau_{\gamma} = \tau_{\gamma'}$. %Each continuous map $f : \Sigma_g^{\mathrm{o}} \to \Sigma_g^{\mathrm{o}}$ which fixes the boundary acts on $\pi_1(\Sigma_g^{\mathrm{o}})$ by $f(c) = f \circ c$.
\smallskip\\
\indent If $S$ is a compact oriented surface, we denote by $\mathrm{MCG}(S)$ its mapping class group, that is the group of isotopy classes of orientation preserving homeomorphisms of $S$ which fix the boundary pointwise. 
\smallskip\\
\indent There exists presentations of $\mathrm{MCG}(\Sigma_g)$ and $\mathrm{MCG}(\Sigma_g^{\mathrm{o}})$ due to Wajnryb \cite{wajnryb} (also see \cite[Sect. 5.2.1]{FM}). Let
\begin{align*}
& d_1 = a_1, \:\:\: d_i = a_{i-1} b_{i} a_{i}^{-1} b_{i}^{-1} \:\:\: \text{ for } 2 \leq i \leq g, \\
& e_1 = a_1, \:\:\: e_i = b_1 a_1^{-1} b_1^{-1} a_1 \ldots b_{i-1} a_{i-1}^{-1} b_{i-1}^{-1} a_{i-1} b_i a_i^{-1} b_i^{-1} \:\:\: \text{ for } 2 \leq i \leq g.
\end{align*}
The correspondence of notations with \cite[Figure 5.7]{FM} is $c_0 = [e_2]$, $c_{2j} = [b_j]$, $c_{2j-1} = [d_j]$. The free homotopy class of these loops are depicted in Figure \ref{figureCourbesCanoniques} below.

%%%%%%%%%%%%%%%%%%%%%%%%%%%%%%%%%%%%%%%%%%%%%%%%%%%
%\begin{comment}

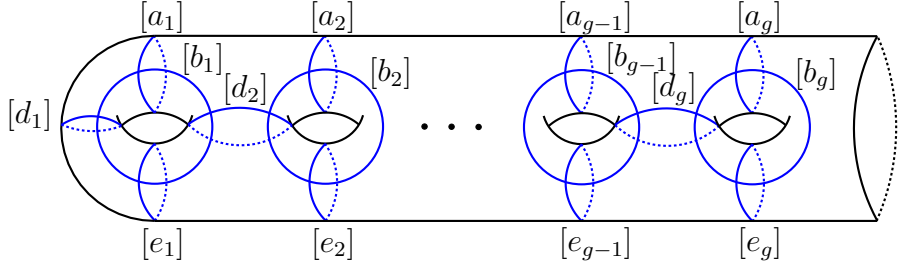
\begin{figure}[h]
\centering
\begin{tikzpicture}
\draw [shift={(0.09650151379600214,1.5429378179389268)},line width=0.8pt]  plot[domain=3.421902658602023:6.006573738596067,variable=\t]({1*0.51666717817069*cos(\t r)+0*0.51666717817069*sin(\t r)},{0*0.51666717817069*cos(\t r)+1*0.51666717817069*sin(\t r)});
\draw [shift={(0.10911145931580171,0.7755953316876043)},line width=0.8pt]  plot[domain=0.8615978203133784:2.312375095503757,variable=\t]({1*0.6611530960524286*cos(\t r)+0*0.6611530960524286*sin(\t r)},{0*0.6611530960524286*cos(\t r)+1*0.6611530960524286*sin(\t r)});
\draw [shift={(0.607947495552267,0.5150456636016608)},line width=0.8pt,color=blue]  plot[domain=2.3506836341863977:3.9339293327304574,variable=\t]({1*0.719075436642927*cos(\t r)+0*0.719075436642927*sin(\t r)},{0*0.719075436642927*cos(\t r)+1*0.719075436642927*sin(\t r)});
\draw [shift={(-0.6689237270985572,0.5071637424741923)},line width=0.8pt,dash pattern=on 1pt off 1pt,color=blue]  plot[domain=-0.5830776988158339:0.5924723398848708,variable=\t]({1*0.9211182116186638*cos(\t r)+0*0.9211182116186638*sin(\t r)},{0*0.9211182116186638*cos(\t r)+1*0.9211182116186638*sin(\t r)});
\draw [shift={(0.6020596132847671,1.9470980286279413)},line width=0.8pt,color=blue]  plot[domain=2.3506836341863977:3.9339293327304574,variable=\t]({1*0.719075436642927*cos(\t r)+0*0.719075436642927*sin(\t r)},{0*0.719075436642927*cos(\t r)+1*0.719075436642927*sin(\t r)});
\draw [shift={(-0.6748116093660571,1.939216107500473)},line width=0.8pt,dash pattern=on 1pt off 1pt,color=blue]  plot[domain=-0.5830776988158339:0.5924723398848708,variable=\t]({1*0.9211182116186638*cos(\t r)+0*0.9211182116186638*sin(\t r)},{0*0.9211182116186638*cos(\t r)+1*0.9211182116186638*sin(\t r)});
\draw [shift={(2.8472037611940717,1.9469618891429243)},line width=0.8pt,color=blue]  plot[domain=2.3506836341863977:3.9339293327304574,variable=\t]({1*0.719075436642927*cos(\t r)+0*0.719075436642927*sin(\t r)},{0*0.719075436642927*cos(\t r)+1*0.719075436642927*sin(\t r)});
\draw [shift={(1.5703325385432474,1.939079968015456)},line width=0.8pt,dash pattern=on 1pt off 1pt,color=blue]  plot[domain=-0.5830776988158339:0.5924723398848708,variable=\t]({1*0.9211182116186638*cos(\t r)+0*0.9211182116186638*sin(\t r)},{0*0.9211182116186638*cos(\t r)+1*0.9211182116186638*sin(\t r)});
\draw [shift={(2.3542556072251064,0.7754591922025873)},line width=0.8pt]  plot[domain=0.8615978203133784:2.312375095503757,variable=\t]({1*0.6611530960524286*cos(\t r)+0*0.6611530960524286*sin(\t r)},{0*0.6611530960524286*cos(\t r)+1*0.6611530960524286*sin(\t r)});
\draw [shift={(2.3416456617053067,1.5428016784539098)},line width=0.8pt]  plot[domain=3.421902658602023:6.006573738596067,variable=\t]({1*0.51666717817069*cos(\t r)+0*0.51666717817069*sin(\t r)},{0*0.51666717817069*cos(\t r)+1*0.51666717817069*sin(\t r)});
\draw [shift={(2.8530916434615716,0.5149095241166438)},line width=0.8pt,color=blue]  plot[domain=2.3506836341863977:3.9339293327304574,variable=\t]({1*0.719075436642927*cos(\t r)+0*0.719075436642927*sin(\t r)},{0*0.719075436642927*cos(\t r)+1*0.719075436642927*sin(\t r)});
\draw [shift={(1.5762204208107473,0.5070276029891753)},line width=0.8pt,dash pattern=on 1pt off 1pt,color=blue]  plot[domain=-0.5830776988158339:0.5924723398848708,variable=\t]({1*0.9211182116186638*cos(\t r)+0*0.9211182116186638*sin(\t r)},{0*0.9211182116186638*cos(\t r)+1*0.9211182116186638*sin(\t r)});
\draw [shift={(1.21306815089848,0.35154293843021955)},line width=0.8pt,color=blue]  plot[domain=0.9265211537579431:2.19965921018045,variable=\t]({1*1.1485995723359301*cos(\t r)+0*1.1485995723359301*sin(\t r)},{0*1.1485995723359301*cos(\t r)+1*1.1485995723359301*sin(\t r)});
\draw [shift={(1.2281202223333165,2.00394503576041)},line width=0.8pt,dash pattern=on 1pt off 1pt,color=blue]  plot[domain=3.953952961918977:5.455765981919788,variable=\t]({1*1.0009623273059072*cos(\t r)+0*1.0009623273059072*sin(\t r)},{0*1.0009623273059072*cos(\t r)+1*1.0009623273059072*sin(\t r)});
\draw [shift={(6.2158455220957,1.9393099892368282)},line width=0.8pt,color=blue]  plot[domain=2.3506836341863977:3.9339293327304574,variable=\t]({1*0.719075436642927*cos(\t r)+0*0.719075436642927*sin(\t r)},{0*0.719075436642927*cos(\t r)+1*0.719075436642927*sin(\t r)});
\draw [shift={(4.938974299444876,1.9314280681093599)},line width=0.8pt,dash pattern=on 1pt off 1pt,color=blue]  plot[domain=-0.5830776988158339:0.5924723398848708,variable=\t]({1*0.9211182116186638*cos(\t r)+0*0.9211182116186638*sin(\t r)},{0*0.9211182116186638*cos(\t r)+1*0.9211182116186638*sin(\t r)});
\draw [shift={(8.460989670005006,1.9391738497518112)},line width=0.8pt,color=blue]  plot[domain=2.3506836341863977:3.9339293327304574,variable=\t]({1*0.719075436642927*cos(\t r)+0*0.719075436642927*sin(\t r)},{0*0.719075436642927*cos(\t r)+1*0.719075436642927*sin(\t r)});
\draw [shift={(7.18411844735418,1.9312919286243428)},line width=0.8pt,dash pattern=on 1pt off 1pt,color=blue]  plot[domain=-0.5830776988158339:0.5924723398848708,variable=\t]({1*0.9211182116186638*cos(\t r)+0*0.9211182116186638*sin(\t r)},{0*0.9211182116186638*cos(\t r)+1*0.9211182116186638*sin(\t r)});
\draw [shift={(5.722897368126735,0.7678072922964911)},line width=0.8pt]  plot[domain=0.8615978203133784:2.312375095503757,variable=\t]({1*0.6611530960524286*cos(\t r)+0*0.6611530960524286*sin(\t r)},{0*0.6611530960524286*cos(\t r)+1*0.6611530960524286*sin(\t r)});
\draw [shift={(5.710287422606935,1.5351497785478136)},line width=0.8pt]  plot[domain=3.421902658602023:6.006573738596067,variable=\t]({1*0.51666717817069*cos(\t r)+0*0.51666717817069*sin(\t r)},{0*0.51666717817069*cos(\t r)+1*0.51666717817069*sin(\t r)});
\draw [shift={(6.826854059709413,0.3437548990391064)},line width=0.8pt,color=blue]  plot[domain=0.9265211537579431:2.19965921018045,variable=\t]({1*1.1485995723359301*cos(\t r)+0*1.1485995723359301*sin(\t r)},{0*1.1485995723359301*cos(\t r)+1*1.1485995723359301*sin(\t r)});
\draw [shift={(6.84190613114425,1.9961569963692969)},line width=0.8pt,dash pattern=on 1pt off 1pt,color=blue]  plot[domain=3.953952961918977:5.455765981919788,variable=\t]({1*1.0009623273059072*cos(\t r)+0*1.0009623273059072*sin(\t r)},{0*1.0009623273059072*cos(\t r)+1*1.0009623273059072*sin(\t r)});
\draw [shift={(7.96804151603604,0.7676711528114741)},line width=0.8pt]  plot[domain=0.8615978203133784:2.312375095503757,variable=\t]({1*0.6611530960524286*cos(\t r)+0*0.6611530960524286*sin(\t r)},{0*0.6611530960524286*cos(\t r)+1*0.6611530960524286*sin(\t r)});
\draw [shift={(7.9554315705162395,1.5350136390627966)},line width=0.8pt]  plot[domain=3.421902658602023:6.006573738596067,variable=\t]({1*0.51666717817069*cos(\t r)+0*0.51666717817069*sin(\t r)},{0*0.51666717817069*cos(\t r)+1*0.51666717817069*sin(\t r)});
\draw [shift={(8.466877552272505,0.5071214847255306)},line width=0.8pt,color=blue]  plot[domain=2.3506836341863977:3.9339293327304574,variable=\t]({1*0.719075436642927*cos(\t r)+0*0.719075436642927*sin(\t r)},{0*0.719075436642927*cos(\t r)+1*0.719075436642927*sin(\t r)});
\draw [shift={(7.1900063296216805,0.4992395635980621)},line width=0.8pt,dash pattern=on 1pt off 1pt,color=blue]  plot[domain=-0.5830776988158339:0.5924723398848708,variable=\t]({1*0.9211182116186638*cos(\t r)+0*0.9211182116186638*sin(\t r)},{0*0.9211182116186638*cos(\t r)+1*0.9211182116186638*sin(\t r)});
\draw [shift={(6.2217334043632,0.5072576242105477)},line width=0.8pt,color=blue]  plot[domain=2.3506836341863977:3.9339293327304574,variable=\t]({1*0.719075436642927*cos(\t r)+0*0.719075436642927*sin(\t r)},{0*0.719075436642927*cos(\t r)+1*0.719075436642927*sin(\t r)});
\draw [shift={(4.944862181712376,0.4993757030830791)},line width=0.8pt,dash pattern=on 1pt off 1pt,color=blue]  plot[domain=-0.5830776988158339:0.5924723398848708,variable=\t]({1*0.9211182116186638*cos(\t r)+0*0.9211182116186638*sin(\t r)},{0*0.9211182116186638*cos(\t r)+1*0.9211182116186638*sin(\t r)});
\draw [line width=0.8pt,color=blue] (0.09852713674962052,1.250536628718387) circle (0.7574954102991686cm);
\draw [line width=0.8pt,color=blue] (2.343671284658925,1.25040048923337) circle (0.7574954102991686cm);
\draw [line width=0.8pt,color=blue] (5.7123130455605535,1.2427485893272738) circle (0.7574954102991686cm);
\draw [line width=0.8pt,color=blue] (7.957457193469859,1.2426124498422568) circle (0.7574954102991686cm);
\draw [shift={(0.09631377499285078,1.227295175105653)},line width=0.8pt]  plot[domain=1.5737998535724984:4.715392507162291,variable=\t]({1*1.2273007109475733*cos(\t r)+0*1.2273007109475733*sin(\t r)},{0*1.2273007109475733*cos(\t r)+1*1.2273007109475733*sin(\t r)});
\draw [line width=0.8pt] (0.09262754998570155,2.454590350211306)-- (9.6,2.449787628);
\draw [line width=0.8pt] (0.1,0)-- (9.6,0);
\draw [shift={(11.892365852439406,1.2268715761423163)},line width=0.8pt]  plot[domain=2.651526337357502:3.6330012492875543,variable=\t]({1*2.5981656743405397*cos(\t r)+0*2.5981656743405397*sin(\t r)},{0*2.5981656743405397*cos(\t r)+1*2.5981656743405397*sin(\t r)});
\draw [shift={(6.7413489409840714,1.23934029514781)},line width=0.8pt,dash pattern=on 1pt off 1pt]  plot[domain=-0.41005552442894366:0.4005426772305154,variable=\t]({1*3.1087631516065994*cos(\t r)+0*3.1087631516065994*sin(\t r)},{0*3.1087631516065994*cos(\t r)+1*3.1087631516065994*sin(\t r)});
\draw [shift={(-0.7121121976730084,2.2393759283163495)},line width=0.8pt, dash pattern=on 1pt off 1pt, color=blue]  plot[domain=4.302386010838742:5.078168584050612,variable=\t]({1*1.0481989230313646*cos(\t r)+0*1.0481989230313646*sin(\t r)},{0*1.0481989230313646*cos(\t r)+1*1.0481989230313646*sin(\t r)});
\draw [shift={(-0.7406297858915413,0.6780379733528217)},line width=0.8pt, color=blue]  plot[domain=0.969994260950464:2.146348567799668,variable=\t]({1*0.7109635065368765*cos(\t r)+0*0.7109635065368765*sin(\t r)},{0*0.7109635065368765*cos(\t r)+1*0.7109635065368765*sin(\t r)});

\draw [fill=black] (3.647040179783018,1.251684065629662) circle (1pt);
\draw [fill=black] (4.002494967436406,1.248870822572264) circle (1pt);
\draw [fill=black] (4.357949755089795,1.246057579514866) circle (1pt);

\draw[color=black] (0.20600112035952642,-0.3) node {$[e_1]$};
\draw[color=black] (0.19669607174480683,2.7) node {$[a_1]$};
\draw[color=black] (-1.55,1.4) node {$[d_1]$};
\draw[color=black] (0.73, 2.15) node {$[b_1]$};
\draw[color=black] (2.4392127878922274,2.7) node {$[a_2]$};
\draw[color=black] (2.448517836506947,-0.3) node {$[e_2]$};
\draw[color=black] (3.2, 1.9547606120430459) node {$[b_2]$};
\draw[color=black] (5.891385823953194,2.7) node {$[a_{g-1}]$};
\draw[color=black] (5.900690872567914,-0.3) node {$[e_{g-1}]$};
\draw[color=black] (6.5,2.15) node {$[b_{g-1}]$};
\draw[color=black] (8.05015710256814,2.7) node {$[a_g]$};
\draw[color=black] (8.059462151182858,-0.3) node {$[e_g]$};
\draw[color=black] (8.8,1.9454555634283333) node {$[b_g]$};
\draw[color=black] (1.294691808281718,1.75) node {$[d_2]$};
\draw[color=black] (6.924246220187068,1.75) node {$[d_g]$};
\end{tikzpicture}
\caption{A canonical set of curves on the surface $\Sigma_g^{\mathrm{o}}$.}\label{figureCourbesCanoniques}
\end{figure}

%\end{comment}
%%%%%%%%%%%%%%%%%%%%%%%%%%%%%%%%%%%%%%%%%%%%%%%%%%%%%%

\begin{figure}[h]
\centering
\begin{tikzpicture}
\draw [line width=0.8pt] (1.9956673267530805,1.4577673242750369) to [bend left=30] (1.7414945936295918,2.234406231041254);
\draw [line width=0.8pt] (1,1)-- (3,1);
\draw [line width=0.8pt, <-] (1.9956673267530805,1.4577673242750369) to [bend right=20] (3,1);
\draw [line width=0.8pt] (5,1)-- (3,1);
\draw [line width=0.8pt] (4.004332673246919,1.4577673242750366) to [bend left=20] (3,1);
\draw [line width=0.8pt, <-] (4.004332673246919,1.4577673242750366) to [bend right=30] (4.258505406370409,2.234406231041254);
\draw [line width=0.8pt,dash pattern=on 1pt off 1pt] (1.7414945936295918,2.234406231041254) to [bend left=50] (4.258505406370409,2.234406231041254);
\draw [fill=black] (3,1) circle (2.5pt);
\draw (3,1.8) node{$\circlearrowright$}; 
\end{tikzpicture}
\caption{A positively oriented loop in the neighborhood of the basepoint fixed in Figure \ref{figureSurfaceRuban}.}
\label{figureCourbeOriente}
\end{figure}
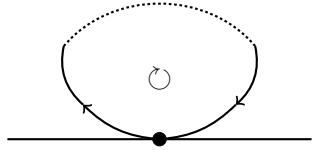

The Dehn twists $\tau_{e_2}, \tau_{b_i}, \tau_{d_i}$ are called the Humphries generators.
Then $\mathrm{MCG}(\Sigma_g^{\mathrm{o}})$ is generated by the Humphries generators together with four families of relations called disjointness relations, braid relations, 3-chain relation and lantern relation, see \cite[Theorem 5.3]{FM}. 
%\begin{itemize}
%\item (Disjointness relations)  $st = ts$ if $s$ and $t$ are disjoint.
%\item (Braid relations) $sts= tst$ if $s$ and $t$ intersect exactly once.
%\item (3-chain relation) 
%\item (Lantern relation)
%\end{itemize}
The presentation of $\mathrm{MCG}(\Sigma_g)$ is obtained as the quotient of $\mathrm{MCG}(\Sigma_g^{\mathrm{o}})$ by the hyperelliptic relation:
\begin{equation}\label{hyperelliptic}
\left( \tau_{b_g} \tau_{d_g} \ldots \tau_{b_1} \tau_{d_1} \tau_{d_1} \tau_{b_1} \ldots \tau_{d_g} \tau_{b_g} \right) \omega = \omega  \left( \tau_{b_g} \tau_{d_g} \ldots \tau_{b_1} \tau_{d_1} \tau_{d_1} \tau_{b_1} \ldots \tau_{d_g} \tau_{b_g} \right)
\end{equation}
where $\omega$ is any word in the Humphries generators which equals $\tau_{a_g}$.
\\\indent The action of the Humphries generators on the fundamental group is easily computed. We just indicate the non-trivial actions:
\begin{equation}\label{actionPi1}
\begin{split}
& \tau_{e_2}(a_1) = e_2^{-1} a_1 e_2, \:\: \tau_{e_2}(b_1) = e_2^{-1} b_1 e_2, \:\: \tau_{e_2}(b_2) = e_2^{-1} b_2, \\
& \tau_{b_i}(a_i) = b^{-1}_i a_i, \\
& \tau_{a_1}(b_1) = b_1 a_1 \:\:\:\: (\text{recall that } a_1 = d_1),\\
& \tau_{d_i}(a_{i-1}) = d_i^{-1} a_{i-1} d_i, \:\: \tau_{d_i}(b_{i-1}) = b_{i-1} d_i, \:\: \tau_{d_i}(b_i) = d_i^{-1} b_i \:\:\: \text{(with } i \geq 2\text{)}.
\end{split}
\end{equation}
\indent In the sequel, we will be concerned with positively oriented, non-separating simple loops in $\pi_1(\Sigma_g^{\mathrm{o}})$. %We say that a loop is simple if it is homotopic to a loop which does not have self-crossings. 
We say that a simple loop is positively oriented if its orientation is clockwise, as indicated in Figure \ref{figureCourbeOriente}. Recall that a loop is non-separating if it does not cut the surface into two connected components and that it is simple if it does not contains self-crossings (up to homotopy). %which shows the loop in a neighbourhood of the basepoint in the thickening of $\Gamma$: %Finally, we say that a loop $\gamma$ is non-separating if $(\Sigma_g^{\mathrm{o}})\setminus \gamma$ is connected.
 It is clear that these properties are preserved by Dehn twists, hence the set of such loops is stable under the action of $\mathrm{MCG}(\Sigma_g^{\mathrm{o}})$. Note that the loops $a_i, b_i, d_i, e_i$ satisfy these properties.

\subsection{Dehn twists as automorphisms of $\mathcal{L}_{g,0}(H)$}
\indent In $\pi_1(\Sigma_g^{\mathrm{o}})$ we have the curves $a_i, b_i$ while in $\mathcal{L}_{g,0}(H)$ we have the matrices $\overset{I}{A}(i), \overset{I}{B}(i)$. Using this, we can lift the action of  the Humphries generators on $\pi_1(\Sigma_g^{\mathrm{o}})$ to $\mathcal{L}_{g,0}(H)$ by replacing loops by matrices, up to some normalization, as we shall see now.
\smallskip\\
\indent First, we lift some of the loops introduced above. We replace the generators of $\pi_1(\Sigma_g^{\mathrm{o}})$ by matrices of generators of $\mathcal{L}_{g,0}(H)$ (see Figure \ref{figureIntro}), up to some normalization. More precisely, we define 
\[ \overset{I}{D}_j = \overset{I}{v}{^2}\overset{I}{A}(j-1) \overset{I}{B}(j) \overset{I}{A}(j){^{-1}} \overset{I}{B}(j){^{-1}}, \:\:\:\:\: \overset{I}{E}_2 = \overset{I}{v}{^4}\overset{I}{B}(1) \overset{I}{A}(1){^{-1}} \overset{I}{B}(1){^{-1}} \overset{I}{A}(1) \overset{I}{B}(2) \overset{I}{A}(2){^{-1}} \overset{I}{B}(2){^{-1}} \]
with $2 \leq j \leq g$. Then $\overset{I}{D}_j$ and $\overset{I}{E}_2$ satisfy the fusion relation of $\mathcal{L}_{0,1}(H)$:
\[ (\overset{I \otimes J}{D_j})_{12} = (\overset{I}{D}_j)_1 \overset{IJ}{(R')}_{12} (\overset{J}{D}_j)_2 (\overset{IJ}{R'})^{-1}_{12}, \:\:\:\:\: (\overset{I \otimes J}{E_2})_{12} = (\overset{I}{E}_2)_1 \overset{IJ}{(R')}_{12} (\overset{J}{E}_2)_2 (\overset{IJ}{R'})^{-1}_{12}. \]
This is easy to check: we observe that $\overset{I}{D}_j = \overset{I}{A}(1) \overset{I}{C}(2) \overset{I}{A}(2){^{-1}}$, $\overset{I}{E}_2 = \overset{I}{C}(1) \overset{I}{C}(2) \overset{I}{A}(2){^{-1}}$ and we use Lemma \ref{decGauss} and relations \eqref{PresentationLgn} to write the fusion and reorder the matrices. Note that the normalizations by powers of $v$ are necessary to have the $\mathcal{L}_{0,1}(H)$-fusion relation on these elements (see the proof of Proposition \ref{liftHumphries} below for an example of computation).
\smallskip\\
Now, we lift the action of the Humphries generators on the fundamental group \eqref{actionPi1}. More precisely, let us define maps $\widetilde{\tau_{e_2}}, \widetilde{\tau_{b_i}}, \widetilde{\tau_{d_j}} : \mathcal{L}_{g,0}(H) \to \mathcal{L}_{g,0}(H)$ by:
\begin{equation}\label{courbesDeviennentMatrices}
\begin{split}
& \widetilde{\tau_{e_2}}(\overset{I}{A}(1)) = \overset{I}{E}{_2^{-1}} \overset{I}{A}(1) \overset{I}{E}_2, \:\: \widetilde{\tau_{e_2}}(\overset{I}{B}(1)) = \overset{I}{E}{_2^{-1}} \overset{I}{B}(1) \overset{I}{E}_2, \:\: \widetilde{\tau_{e_2}}(\overset{I}{B}(2)) = \overset{I}{v}\overset{I}{E}{_2^{-1}} \overset{I}{B}(2), \\
& \widetilde{\tau_{b_i}}(\overset{I}{A}(i)) = \overset{I}{v} \overset{I}{B}(i){^{-1}} \overset{I}{A}(i), \\
& \widetilde{\tau_{a_1}}(\overset{I}{B}(1)) = \overset{I}{v}{^{-1}} \overset{I}{B}(1) \overset{I}{A}(1) \:\:\:\: (\text{recall that } a_1 = d_1),\\
& \widetilde{\tau_{d_j}}(\overset{I}{A}(j-1)) = \overset{I}{D}{_j^{-1}} \overset{I}{A}(j-1) \overset{I}{D}_j, \:\: \widetilde{\tau_{d_j}}(\overset{I}{B}(j-1)) = \overset{I}{v}{^{-1}}\overset{I}{B}(j-1) \overset{I}{D}_j, \:\: \widetilde{\tau_{d_j}}(\overset{I}{B}(j)) = \overset{I}{v}\overset{I}{D}{_j^{-1}}\overset{I}{B}(j),
\end{split}
\end{equation}
for $j \geq 2$, and the other matrices are fixed.

\begin{proposition}\label{liftHumphries} 1) The maps $\widetilde{\tau_{e_2}}, \widetilde{\tau_{b_i}}, \widetilde{\tau_{d_i}}$ are automorphisms of $\mathcal{L}_{g,0}(H)$.
\\2) The assignment
\[ \tau_{e_2} \mapsto \widetilde{\tau_{e_2}}, \:\:\: \tau_{b_i} \mapsto \widetilde{\tau_{b_i}}, \:\:\: \tau_{d_i} \mapsto \widetilde{\tau_{d_i}} \]
extends to a morphism of groups $\mathrm{MCG}(\Sigma_g^{\mathrm{o}}) \to \mathrm{Aut}(\mathcal{L}_{g,0}(H))$.
\end{proposition}
\debutDemo
1) We have to check that these maps are compatible with the defining relations \eqref{PresentationLgn}. This relies on straightforward but tedious computations. For instance, let us show that $\widetilde{\tau_{d_j}}(\overset{I}{B}(j-1))$ satisfies the fusion relation. First, it is easy to establish the following exchange relation:
\[ (\overset{IJ}{R'})_{12} \overset{J}{B}(j-1)_2 \overset{IJ}{R}_{12} (\overset{I}{D}_j)_1 (\overset{IJ}{R'})_{12} = (\overset{I}{D}_j)_1 (\overset{IJ}{R'})_{12} \overset{J}{B}(j-1)_2. \]
Hence 
\begin{align*}
\overset{I \otimes J}{B}\!(j-1)_{12} (\overset{I \otimes J}{D_j})_{12} &= \overset{I}{B}(j-1)_1 (\overset{IJ}{R'})_{12} \overset{J}{B}(j-1)_2 (\overset{IJ}{R'})^{-1}_{12} (\overset{I}{D}_j)_1 (\overset{IJ}{R'})_{12} (\overset{J}{D}_j)_2 (\overset{IJ}{R'})^{-1}_{12}\\
&= \overset{I}{B}(j-1)_1 (\overset{IJ}{R'})_{12} \overset{J}{B}(j-1)_2 (\overset{IJ}{R'})^{-1}_{12} \overset{I \otimes J}{v}_{\!\!\! 12} (\overset{IJ}{R'})_{12} \overset{IJ}{R}_{12} \overset{I}{v}{_1^{-1}} \overset{J}{v}{_2^{-1}} (\overset{I}{D}_j)_1 (\overset{IJ}{R'})_{12} (\overset{J}{D}_j)_2 (\overset{IJ}{R'})^{-1}_{12}\\
&= \overset{I \otimes J}{v}_{\!\!\! 12} \overset{I}{v}{_1^{-1}} \overset{J}{v}{_2^{-1}} \overset{I}{B}(j-1)_1 (\overset{IJ}{R'})_{12} \overset{J}{B}(j-1)_2 \overset{IJ}{R}_{12} (\overset{I}{D}_j)_1 (\overset{IJ}{R'})_{12} (\overset{J}{D}_j)_2 (\overset{IJ}{R'})^{-1}_{12}\\
&= \overset{I \otimes J}{v}_{\!\!\! 12} \overset{I}{v}{_1^{-1}} \overset{I}{B}(j-1)_1 (\overset{I}{D}_j)_1 (\overset{IJ}{R'})_{12} \overset{J}{v}{_2^{-1}} \overset{J}{B}(j-1)_2 (\overset{J}{D}_j)_2 (\overset{IJ}{R'})^{-1}_{12}.
\end{align*}
%For the first equality, we used that $B(j-1)$ and $D_j$ satisfy the fusion relation
For the second equality we applied a trick based on \eqref{ribbon}. The aim is to replace $R'^{-1}$ by $R$ in order to apply the previously established exchange relation. It follows that $\overset{I}{v}{^{-1}} \overset{I}{B}(j-1) \overset{I}{D}_j$ satisfies the fusion relation, as desired. This computation shows how the normalizations by powers of $v$ arises in order to satisfy the fusion relation. These normalizations have no importance when one checks the compatibility with the other defining relations of $\mathcal{L}_{g,0}(H)$.
\\\noindent 2) Straightforward verification using Wajnryb's relations. %Note that the only fact to check is that the normalizations by powers of $v$ agree with the relations, since formally this is the only difference with the action of $\mathrm{MCG}(\Sigma_g^{\mathrm{o}})$ on $\pi_1(\Sigma_g^{\mathrm{o}})$.
\finDemo

\begin{definition}
The lift of an element $f \in \mathrm{MCG}(\Sigma_g^{\mathrm{o}})$, denoted by $\widetilde{f}$, is its image by the morphism of Proposition \ref{liftHumphries}.
\end{definition}
\noindent Let $u_i$ be one of the generators of $\pi_1(\Sigma_g^{\mathrm{o}})$, let $f \in \mathrm{MCG}(\Sigma_g^{\mathrm{o}})$ and let $f(u_i) = a_{i_1}^{m_1} b_{j_1}^{n_1} \ldots a_{i_k}^{m_k} b_{j_k}^{n_k}$ with $m_{\ell}, n_{\ell} \in \mathbb{Z}$. Then it follows from the definition of $\widetilde{\tau_{e_2}}, \widetilde{\tau_{b_i}}, \widetilde{\tau_{d_i}}$ that 
\begin{equation}\label{expressionLift}
f(\overset{I}{U}(i)) = \overset{I}{v}{^N}\overset{I}{A}(i_1)^{m_1} \overset{I}{B}(j_1)^{n_1} \ldots \overset{I}{A}(i_k)^{m_k} \overset{I}{B}(j_k)^{n_k}
\end{equation}
where $U(i)=A(i)$ (resp. $U(i)=B(i)$) if $u_i = a_i$ (resp. $u_i = b_i$) and for some $N \in \mathbb{Z}$. In other words, $f$ and $\widetilde{f}$ are formally identical except for the normalizations by some power of $v$.
\smallskip\\
\indent Recall that $\mathcal{L}_{g,0}(H) \cong \mathrm{End}_{\mathbb{C}}\!\left((H^*)^{\otimes g}\right)$ is a matrix algebra. By the Skolem-Noether theorem, every automorphism of $\mathcal{L}_{g,0}(H)$ is inner. Hence to each $f \in \mathrm{MCG}(\Sigma_g^{\mathrm{o}})$ is associated an element $\widehat{f} \in \mathcal{L}_{g,0}(H)$, unique up to scalar, such that
\begin{equation}\label{conjugaison}
\forall\, x \in \mathcal{L}_{g,0}(H), \:\:\: \widetilde{f}(x) = \widehat{f}x\widehat{f}^{-1}.
\end{equation}
\indent We now determine the elements $\widehat{\tau_{\gamma}}$ associated to Dehn twists.

\begin{lemme}\label{lemmeA1}
We have $\widehat{\tau_{a_1}} = v_{A(1)}^{-1}$. In other words:
\[ \forall\, x \in \mathcal{L}_{g,0}(H), \:\: \widetilde{\tau_{a_1}}(x) = v_{A(1)}^{-1} \, x \, v_{A(1)}. \]
\end{lemme}
\debutDemo
We have $v_{A(1)}^{-1} \overset{I}{A}(1) = \overset{I}{A}(1) v_{A(1)}^{-1} = \widetilde{\tau_{a_1}}(\overset{I}{A}(1))v_{A(1)}^{-1}$. Indeed, since $v^{-1}$ is central in $H$, $v_{A(1)}^{-1}$ is central in the subalgebra generated by the coefficients of the matrices $\overset{I}{A}(1)$. Next, let $j_1 : \mathcal{H}(\mathcal{O}(H)) \to \mathcal{H}(\mathcal{O}(H))^{\otimes g}$ be the canonical embedding on the first copy. Observe that for all $x \in H$, $\Psi_{g,0}(x_{A(1)}) = j_1(x)$. Then:
%\[ \Psi_{g,0}(v_{A(1)}^{-1} \overset{I}{A}(1)) = j_1(v^{-1}\overset{I}{L}{^{(+)}}\overset{I}{L}{^{(-)-1}}) = \overset{I}{D}_j \,j_1(v^{-1}Y_j) = \overset{I}{D}_j \,j_1(Y_j v^{-1}) = j_1(\overset{I}{L}{^{(+)}}\overset{I}{L}{^{(-)-1}}v^{-1}) = \Psi_{g,0}(\overset{I}{A}(1) v_{A(1)}^{-1}). \]
\begin{align*}
\Psi_{g,0}(v_{A(1)}^{-1} \overset{I}{B}(1)) &= j_1(v^{-1} \overset{I}{L}{^{(+)}} \, \overset{I}{T} \, \overset{I}{L}{^{(-)-1}}) = j_1(\overset{I}{L}{^{(+)}} \, \overset{I}{T} \, \overset{I}{(v'^{-1})} \, v''^{-1} \,\overset{I}{L}{^{(-)-1}})\\
& = j_1(\overset{I}{L}{^{(+)}} \, \overset{I}{T} \, \overset{I}{v}{^{-1}}\overset{I}{b_i} \overset{I}{a_j} a_i b_j \, v^{-1} \,\overset{I}{L}{^{(-)-1}}) = j_1(\overset{I}{v}{^{-1}} \overset{I}{L}{^{(+)}} \, \overset{I}{T} \, \overset{I}{L}{^{(-)-1}} \, \overset{I}{L}{^{(+)}} \,\overset{I}{L}{^{(-)-1}} v^{-1})\\
& = \Psi_{g,0}(\overset{I}{v}{^{-1}} \overset{I}{B}(1) \overset{I}{A}(1) v_{A(1)}^{-1}) =  \Psi_{g,0}( \widetilde{\tau_{a_1}}(\overset{I}{B}(1)) v_{A(1)}^{-1}).
\end{align*}
We used the exchange relation \eqref{relDefHeisenberg} of $\mathcal{H}(\mathcal{O}(H))$ together with \eqref{ribbon} and the definition of the matrices $\overset{I}{L}{^{(\pm)}}$. Finally, recall the matrices \eqref{matricesAlekseev} which occur in the definition of the Alekseev isomorphism. The same argument as in the proof of Lemma \ref{expressionM} shows that
\[ \Psi_{1,0}^{\otimes g}(\overset{I}{\Lambda}_i) = \overset{I}{S^{-1}(b_{\ell})} \: \widetilde{a_{\ell}^{(2i -1)}} a_{\ell}^{(2i)} \otimes \ldots \otimes \widetilde{a_{\ell}^{(1)}}b_{\ell}^{(2)}. \]
From this we see that $j_1(v^{-1})$ commutes with $\Psi_{1,0}^{\otimes g}(\overset{I}{\Lambda}_i)$. Eventually it follows that $\Psi_{g,0}(v^{-1}_{A(1)})$ commutes with $\Psi_{g,0}(\overset{I}{U}(i)) = \Psi_{g,0}(\widetilde{\tau_{a_1}}(\overset{I}{U}(i)))$, where $U$ is $A$ or $B$.
%\begin{align*}
%\Psi_{g,0}(v_{A(1)}^{-1} \overset{I}{U}(i)) &= j_1(v^{-1}) \Psi_{1,0}^{\otimes g}(\overset{I}{\Lambda}_i) j_i(\Psi_{1,0}(\overset{I}{U})) \Psi_{1,0}^{\otimes g}(\overset{I}{\Lambda}_i)^{-1} = \Psi_{1,0}^{\otimes g}(\overset{I}{\Lambda}_i) j_i(\Psi_{1,0}(\overset{I}{U})) \Psi_{1,0}^{\otimes g}(\overset{I}{\Lambda}_i)^{-1}  j_1(v^{-1}) 
%\\&= \Psi_{g,0}(\overset{I}{U}(i) v_{A(1)}^{-1}) = \Psi_{g,0}(\widetilde{\tau_{a_1}}(\overset{I}{U}(i)) v_{A(1)}^{-1})
%\end{align*}
\finDemo

\indent If $\gamma_1, \gamma_2$ are non-separating circles on a surface, it is well known that there exists a homeomorphism $f$ such that $f(\gamma_1)$ is freely homotopic to $\gamma_2$ (see e.g. \cite[Sect. 1.3.1]{FM}). Here we need to consider fixed-point homotopies.

\begin{lemme}\label{transfoLoop}
Let $\gamma_1, \gamma_2$ be positively oriented, non-separating simple loops in $\pi_1(\Sigma_g^{\mathrm{o}})$, then there exists $f \in \mathrm{MCG}(\Sigma_g^{\mathrm{o}})$ such that $f(\gamma_1) = \gamma_2$ in $\pi_1(\Sigma_g^{\mathrm{o}})$.
\end{lemme}
\debutDemo
We know that there exists $\eta \in \mathrm{MCG}(\Sigma_g^{\mathrm{o}})$ such that $\eta(\gamma_1) = \gamma_2' = \alpha^{\varepsilon} \gamma_2^{\pm 1} \alpha^{-\varepsilon}$ in $\pi_1(\Sigma_g^{\mathrm{o}})$ for some loop $\alpha$ and some $\varepsilon \in \{\pm 1\}$. $\gamma'_2$ is positively oriented, non-separating and simple since $\gamma_1$ is, and thus we can assume that $\alpha$ is simple and does not intersect $\gamma_2$ (except at the basepoint). There are six possible configurations for the loops $\alpha$ and $\gamma_2$ in a neighbourhood of the basepoint:

%%%%%%%%%%%%%%%%%%%%%%%%%%%%%%%%%%%%%%%%%%%%%%%%%%%%%%
%\begin{comment}

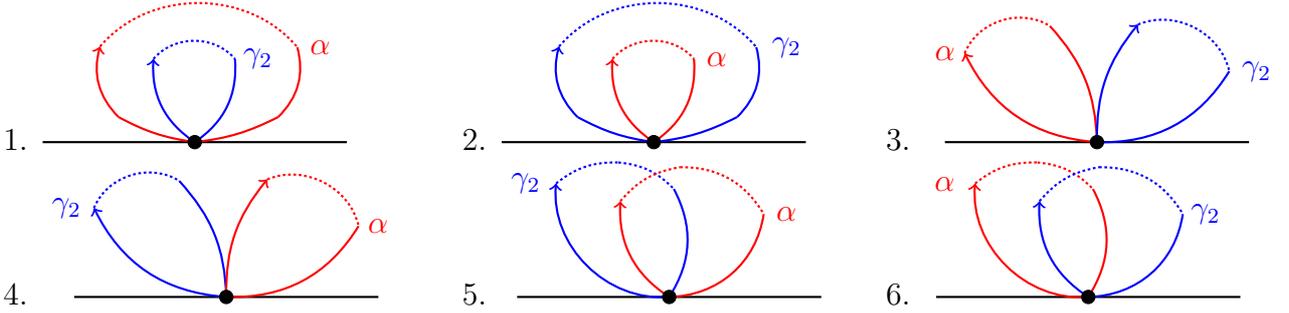
\begin{center}
\begin{tabular}{l l l}
1. \begin{tikzpicture}
\draw [line width=0.8pt] (1,1)-- (5,1);
\draw [line width=0.8pt, color=red] (3,1) to [bend left=10] (1.9983888861853318,1.332201709197648);
\draw [line width=0.8pt, color=red, ->] (1.9983888861853318,1.332201709197648) to [bend left=30] (1.7552819764080068,2.26861350982142);
\draw [line width=0.8pt,dash pattern=on 1pt off 1pt, color=red] (1.7552819764080068,2.26861350982142) to [bend left=50] (4.3484223473661405,2.2596095502000373);
\draw [line width=0.8pt, color=red] (3,1) to [bend right=10] (4.096311477967434,1.332201709197648);
\draw [line width=0.8pt, color=red] (4.096311477967434,1.332201709197648) to [bend right=30] (4.3484223473661405,2.2596095502000373);
\draw [line width=0.8pt, color=blue, ->] (3,1) to [bend left=30] (2.4575908268758346,2.103680132957095);
\draw [line width=0.8pt, color=blue] (3,1) to [bend right=30] (3.5290620218203412,2.103680132957095);
\draw [line width=0.8pt,dash pattern=on 1pt off 1pt, color=blue] (2.4575908268758346,2.103680132957095) to [bend left=50] (3.5290620218203412,2.103680132957095);
\draw [fill=black] (3,1) circle (2.5pt);
\draw[color=blue] (3.5290620218203412+0.3,2.103680132957095) node{$\gamma_2$};
\draw[color=red] (4.3484223473661405+0.3,2.2596095502000373) node{$\alpha$};
\end{tikzpicture}
&
~~~2. \begin{tikzpicture}
\draw [line width=0.8pt] (1,1)-- (5,1);
\draw [line width=0.8pt, color=blue] (3,1) to [bend left=10] (1.9983888861853318,1.332201709197648);
\draw [line width=0.8pt, color=blue, ->] (1.9983888861853318,1.332201709197648) to [bend left=30] (1.7552819764080068,2.26861350982142);
\draw [line width=0.8pt,dash pattern=on 1pt off 1pt, color=blue] (1.7552819764080068,2.26861350982142) to [bend left=50] (4.3484223473661405,2.2596095502000373);
\draw [line width=0.8pt, color=blue] (3,1) to [bend right=10] (4.096311477967434,1.332201709197648);
\draw [line width=0.8pt, color=blue] (4.096311477967434,1.332201709197648) to [bend right=30] (4.3484223473661405,2.2596095502000373);
\draw [line width=0.8pt, color=red, ->] (3,1) to [bend left=30] (2.4575908268758346,2.103680132957095);
\draw [line width=0.8pt, color=red] (3,1) to [bend right=30] (3.5290620218203412,2.103680132957095);
\draw [line width=0.8pt,dash pattern=on 1pt off 1pt, color=red] (2.4575908268758346,2.103680132957095) to [bend left=50] (3.5290620218203412,2.103680132957095);
\draw [fill=black] (3,1) circle (2.5pt);
\draw[color=red] (3.5290620218203412+0.3,2.103680132957095) node{$\alpha$};
\draw[color=blue] (4.3484223473661405+0.4,2.2596095502000373) node{$\gamma_2$};
\end{tikzpicture}
&
~~~3. \begin{tikzpicture}
\draw [line width=0.8pt] (1,1)-- (5,1);
\draw [line width=0.8pt, color=red, ->] (3,1) to [bend left=30] (1.26,2.18);
\draw [line width=0.8pt, color=red] (3,1) to [bend right=20] (2.38,2.54);
\draw [line width=0.8pt, color=red, dash pattern=on 1pt off 1pt] (1.26,2.18) to [bend left=50] (2.38,2.54);
\draw [line width=0.8pt, color=blue, ->] (3,1) to [bend left=20] (3.54,2.56);
\draw [line width=0.8pt,dash pattern=on 1pt off 1pt, color=blue] (3.54,2.56) to [bend left=50] (4.74,1.94);
\draw [line width=0.8pt, color=blue] (3,1) to [bend right=30] (4.74,1.94);
\draw [fill=black] (3,1) circle (2.5pt);
\draw [color=red] (1,2.18) node{$\alpha$};
\draw[color=blue] (5.1,1.94) node{$\gamma_2$};
\end{tikzpicture}
\\
4. \begin{tikzpicture}
\draw [line width=0.8pt] (1,1)-- (5,1);
\draw [line width=0.8pt, color=blue, ->] (3,1) to [bend left=30] (1.26,2.18);
\draw [line width=0.8pt, color=blue] (3,1) to [bend right=20] (2.38,2.54);
\draw [line width=0.8pt, color=blue, dash pattern=on 1pt off 1pt] (1.26,2.18) to [bend left=50] (2.38,2.54);
\draw [line width=0.8pt, color=red, ->] (3,1) to [bend left=20] (3.54,2.56);
\draw [line width=0.8pt,dash pattern=on 1pt off 1pt, color=red] (3.54,2.56) to [bend left=50] (4.74,1.94);
\draw [line width=0.8pt, color=red] (3,1) to [bend right=30] (4.74,1.94);
\draw [fill=black] (3,1) circle (2.5pt);
\draw [color=blue] (1-0.1,2.18) node{$\gamma_2$};
\draw[color=red] (5,1.94) node{$\alpha$};
\end{tikzpicture}
&
~~~5. \begin{tikzpicture}
\draw [line width=0.8pt] (1,1)-- (5,1);
\draw [line width=0.8pt, color=blue, ->] (3,1) to [bend left=50] (1.5121750666306817,2.3046293483069493+0.2); 
\draw [line width=0.8pt,dash pattern=on 1pt off 1pt, color=blue] (1.5121750666306817,2.3046293483069493+0.2) to [bend left=20] (2.3315353921764808,2.5837520965698046+0.2);
\draw [line width=0.8pt,dash pattern=on 1pt off 1pt, color=blue] (2.3315353921764808,2.5837520965698046+0.2) to [bend left=20] (3.060856121508456,2.23259767133589+0.2);
\draw [line width=0.8pt, color=blue] (3.060856121508456,2.23259767133589+0.2) to [bend left=30] (3,1);
\draw [line width=0.8pt, color=red, ->] (3,1) to [bend left=30] (2.3585472710406283,2.0705263981510065+0.2);
\draw [line width=0.8pt,dash pattern=on 1pt off 1pt, color=red] (2.3585472710406283,2.0705263981510065+0.2) to [bend left=20] (3.141891758100898,2.5207243792201277+0.2);
\draw [line width=0.8pt,dash pattern=on 1pt off 1pt, color=red] (3.141891758100898,2.5207243792201277+0.2) to [bend left=30] (3.943244164403932+0.3,1.8994511653447406+0.2);
\draw [line width=0.8pt, color=red] (3.943244164403932+0.3,1.8994511653447406+0.2) to [bend left=40] (3,1);
\draw [fill=black] (3,1) circle (2.5pt);
\draw[color=blue] (1.5121750666306817-0.4,2.3046293483069493+0.2) node{$\gamma_2$};
\draw[color=red] (3.943244164403932+0.6,1.8994511653447406+0.2) node{$\alpha$};
\end{tikzpicture}
&
~~~6. \begin{tikzpicture}
\draw [line width=0.8pt] (1,1)-- (5,1);
\draw [line width=0.8pt, color=red, ->] (3,1) to [bend left=50] (1.5121750666306817,2.3046293483069493+0.2); 
\draw [line width=0.8pt,dash pattern=on 1pt off 1pt, color=red] (1.5121750666306817,2.3046293483069493+0.2) to [bend left=20] (2.3315353921764808,2.5837520965698046+0.2);
\draw [line width=0.8pt,dash pattern=on 1pt off 1pt, color=red] (2.3315353921764808,2.5837520965698046+0.2) to [bend left=20] (3.060856121508456,2.23259767133589+0.2);
\draw [line width=0.8pt, color=red] (3.060856121508456,2.23259767133589+0.2) to [bend left=30] (3,1);
\draw [line width=0.8pt, color=blue, ->] (3,1) to [bend left=30] (2.3585472710406283,2.0705263981510065+0.2);
\draw [line width=0.8pt,dash pattern=on 1pt off 1pt, color=blue] (2.3585472710406283,2.0705263981510065+0.2) to [bend left=20] (3.141891758100898,2.5207243792201277+0.2);
\draw [line width=0.8pt,dash pattern=on 1pt off 1pt, color=blue] (3.141891758100898,2.5207243792201277+0.2) to [bend left=30] (3.943244164403932+0.3,1.8994511653447406+0.2);
\draw [line width=0.8pt, color=blue] (3.943244164403932+0.3,1.8994511653447406+0.2) to [bend left=40] (3,1);
\draw [fill=black] (3,1) circle (2.5pt);
\draw[color=red] (1.5121750666306817-0.4,2.3046293483069493+0.2) node{$\alpha$};
\draw[color=blue] (3.943244164403932+0.6,1.8994511653447406+0.2) node{$\gamma_2$};
\end{tikzpicture}
\end{tabular}
\end{center}

%\end{comment}
%%%%%%%%%%%%%%%%%%%%%%%%%%%%%%%%%%%%%%%%%%%%%%%%%%%%%%
In case 1, $\gamma_2' = \alpha \gamma_2 \alpha^{-1}$, and then $\tau_{\alpha}(\gamma'_2) = \alpha^{-1} \gamma'_2 \alpha = \gamma_2$. Case 2 is impossible because none of the four possible loops $\alpha^{\varepsilon} \gamma_2^{\pm 1} \alpha^{-\varepsilon}$ is simple. In case 3, $\gamma'_2 = \alpha \gamma_2 \alpha^{-1}$. For $\beta = \alpha \gamma_2$, we have $\tau_{\beta}(\alpha) = \beta^{-1}\alpha \beta$, $\tau_{\beta}(\gamma_2) = \beta^{-1} \gamma_2 \beta$, and thus $\tau_{\beta}(\gamma'_2) = \gamma_2$. In case 4, $\gamma'_2 = \alpha^{-1} \gamma_2 \alpha$. For $\delta = \gamma_2 \alpha$, we get similarly to case 3 that $\tau_{\delta}^{-1}(\gamma'_2) = \delta \gamma'_2 \delta^{-1} = \gamma_2$. In case 5, $\gamma'_2 = \alpha^{-1} \gamma_2^{-1} \alpha$. Observe that $\tau_{\alpha}(\gamma_2) = \gamma_2 \alpha$, $\tau_{\gamma_2}(\alpha) = \gamma_2^{-1}\alpha$, and then
\[ \tau_{\alpha}^{-1} \tau_{\gamma_2}^{-2} \tau_{\alpha}^{-1}(\alpha^{-1} \gamma^{-1}_2 \alpha) = \tau_{\alpha}^{-1} \tau_{\gamma_2}^{-2}(\gamma^{-1}_2 \alpha) = \tau_{\alpha}^{-1}(\gamma_2 \alpha) =  \gamma_2.\]
In case 6, $\gamma'_2 = \alpha \gamma_2^{-1} \alpha^{-1}$, and we get similarly to case 5 that $\tau_{\alpha} \tau_{\gamma_2}^{2} \tau_{\alpha}(\gamma'_2) = \gamma_2$.
\finDemo

\begin{exemple}
We have
\[
\tau_{b_i} \tau_{a_i}(b_i) = a_i, \:\:\:\:\:
\tau_{d_i}^{-1} \tau_{b_{i-1}}^{-1}(d_i) = b_{i-1}, \:\:\:\:\:
\tau_{y_i}^{-1} \tau_{a_i}^{-1} \tau_{b_{i-1}}^{-1} \tau_{y_i}^{-1}(a_i) = b_{i-1}, \:\:\:\:\:
\tau_{y_2}^{-1}\tau_{b_1}^{-1}\tau_{e_2}\tau_{y_2}(e_2) = b_1
\]
where $y_i=a_{i-1}b_i$. This allows to transform any of the loops $a_i, b_i, d_i, e_2$ into $a_1$.
\end{exemple}

\begin{lemme}\label{liftUnique}
%Let $\gamma = s_1^{\pm 1} \ldots s_n^{\pm 1} \in \pi_1(\Sigma_g^{\mathrm{o}})$, where the $s_i$ are generators of the fundamental group. A lift of $\gamma$ is $\overset{I}{\widetilde{\gamma}} = \overset{I}{v}{^k} \overset{I}{S}_1^{\pm 1} \ldots \overset{I}{S}_n^{\pm 1}$ for all $I$, which satisfies the fusion relation, and where $S_i = A(j)$ if $s_i = a_j$ and $S_i = B(j)$ if $s_i = b_j$.
Let $f, g \in \mathrm{MCG}(\Sigma_g^{\mathrm{o}})$ such that $f(a_1) = g(a_1)$. Then $\widetilde{f}(\overset{I}{A}(1)) = \widetilde{g}(\overset{I}{A}(1))$.
\end{lemme}
\debutDemo
Let $\eta \in \mathrm{MCG}(\Sigma_g^{\mathrm{o}})$ be such that $\eta(a_1) = a_1$. A priori, $\widetilde{\eta}(\overset{I}{A}(1)) = \overset{I}{v}{^N}\overset{I}{A}(1)$ (see \eqref{expressionLift}) and we must show that $N=0$. Let $\mathrm{MCG}(\Sigma_g^{\mathrm{o}})_{[a_1]}$ be the stabilizer of the free homotopy class $[a_1]$. There is a surjection $p : \mathrm{MCG}(\Sigma_g^{\mathrm{o}} \setminus [a_1]) \to \mathrm{MCG}(\Sigma_g^{\mathrm{o}})_{[a_1]}$. $\mathrm{MCG}(\Sigma_g^{\mathrm{o}} \setminus [a_1])$ is generated by $\tau_{d_i}, \tau_{b_i}, \tau_{e_2}, \tau_{a_g}, \tau_{e_g}, \tau_{\beta}$ with $i \geq 2$, where $\tau_{\beta}$ is such that $p(\tau_{\beta}) = \tau_{e_g}$ (see \cite[Figure 4.10]{FM}). It follows that $\mathrm{MCG}(\Sigma_g^{\mathrm{o}})_{[a_1]} = \langle \tau_{d_i}, \tau_{b_i}, \tau_{e_2}, \tau_{a_g}, \tau_{e_g} \rangle_{i \geq 2}$. For each of these generators $h$, it is possible to verify directly that 
\begin{align*}
&h(a_1) = \gamma_h a_1 \gamma_h^{-1} \: \text{ in } \: \pi_1(\Sigma_g^{\mathrm{o}}), \: \text{ with } \: \gamma_h = a_{i_1}^{m_1} b_{j_1}^{n_1} \ldots a_{i_k}^{m_k} b_{j_k}^{n_k} \:\: (m_{\ell}, n_{\ell} \in \mathbb{Z}),\\
&\widetilde{h}(\overset{I}{A}(1)) = \overset{I}{\Gamma}_h \overset{I}{A}(1) \overset{I}{\Gamma}{^{-1}_h} \: \text{ in } \: \mathcal{L}_{g,0}(H), \: \text{ with } \: \Gamma_h = \overset{I}{v}{^r}\overset{I}{A}(i_1)^{m_1} \overset{I}{B}(j_1)^{n_1} \ldots \overset{I}{A}(i_k)^{m_k} \overset{I}{B}(j_k)^{n_k} \:\: (r, m_{\ell}, n_{\ell} \in \mathbb{Z}).
\end{align*}
In other words, $h(a_1)$ and $\widetilde{h}(\overset{I}{A}(1))$ are formally identical, without any power of $v$ in the expression of $\widetilde{h}(\overset{I}{A}(1))$ (note that the normalization of $\Gamma_h$ by a power of $v$ vanishes in the conjugation). We deduce that this property is true for any $h \in \mathrm{MCG}(\Sigma_g^{\mathrm{o}})_{[a_1]}$. Since $\eta(a_1) = a_1$ in $\pi_1(\Sigma_g^{\mathrm{o}})$, then in particular $\eta \in \mathrm{MCG}(\Sigma_g^{\mathrm{o}})_{[a_1]}$, and thus $\widetilde{\eta}(\overset{I}{A}(1)) = \overset{I}{A}(1)$, as desired.
\finDemo
\smallskip\\
This lemma justifies the following definition.
\begin{definition}
Let $\gamma \in \pi_1(\Sigma_g^{\mathrm{o}})$ be a positively-oriented, non-separating simple loop, and let $f \in \mathrm{MCG}(\Sigma_g^{\mathrm{o}})$ be such that $f(a_1) = \gamma$. The lift of $\gamma$ is $\overset{I}{\widetilde{\gamma}} = \widetilde{f}(\overset{I}{A}(1))$.
\end{definition}
\noindent Some comments are in order. First, if $\gamma = a_{i_1}^{m_1} b_{j_1}^{n_1} \ldots a_{i_k}^{m_k} b_{j_k}^{n_k}$ with $m_{\ell}, n_{\ell} \in \mathbb{Z}$ is a positively oriented, non-separating simple loop, then $\overset{I}{\widetilde{\gamma}} = \overset{I}{v}{^N}\overset{I}{A}(i_1)^{m_1} \overset{I}{B}(j_1)^{n_1} \ldots \overset{I}{A}(i_k)^{m_k} \overset{I}{B}(j_k)^{n_k}$ with $N \in \mathbb{Z}$. In other words, $\gamma$ and $\widetilde{\gamma}$ are formally identical except for the normalization by a power of $v$. Note that by definition every lift satisfies the fusion relation of $\mathcal{L}_{0,1}(H)$. We mention again that the normalization by a power of $v$ is required to satisfy the fusion relation (see e.g. the proof of Proposition \ref{liftHumphries}). %The power can be determined explicitly (see \cite{AGS2, AS}), but we do not need it.
\smallskip\\
\indent We can now answer the question of what are the elements implementing lifting of Dehn twists by conjugation. We use the notation introduced at the end of subsection 3.2.
\begin{proposition}\label{propDehnTwist}
For any non-separating circle $\gamma$ on $\Sigma_g^{\mathrm{o}}$, we have $\widehat{\tau_{\gamma}} = v_{\widetilde{\gamma}}^{-1}$. In other words:
\[ \forall\, x \in \mathcal{L}_{g,0}(H), \:\: \widetilde{\tau_{\gamma}}(x) = v_{\widetilde{\gamma}}^{-1} \, x \, v_{\widetilde{\gamma}}. \]
If $\gamma, \delta \in \pi_1(\Sigma_g^{\mathrm{o}})$ are positively oriented non-separating simple loops such that $[\gamma] = [\delta]$, then $v_{\widetilde{\gamma}}$ is proportional to $v_{\widetilde{\delta}}$
\end{proposition}
\debutDemo
We represent the circle $[\gamma]$ by a positively-oriented, non-separating simple loop $\gamma \in \pi_1(\Sigma_g^{\mathrm{o}})$. Let $f \in \mathrm{MCG}(\Sigma_g^{\mathrm{o}})$ be such that $f(a_1) = \gamma$, then
\[ \widetilde{\tau_{\gamma}} = \widetilde{\tau_{f(a_1)}} = \widetilde{f \tau_{a_1}f^{-1}} = \widetilde{f} \, \widetilde{\tau_{a_1}} \, \widetilde{f}^{-1}. \]
Hence, by Lemma \ref{lemmeA1},
\[ \widetilde{\tau_{\gamma}}\!\left(\widetilde{f}(x)\right) = \widetilde{f}\!\left(\widetilde{\tau_{a_1}}(x)\right) = \widetilde{f}\!\left(v_{A(1)}^{-1} x v_{A(1)}\right) = v_{\widetilde{\gamma}}^{-1} \widetilde{f}(x) v_{\widetilde{\gamma}}. \]
Replacing $x$ by $\widetilde{f}^{-1}(x)$, we get the result. The second claim follows from a similar reasoning together with the fact that $\tau_{\gamma}$ depends only of the free homotopy class of $\gamma$.
\finDemo

Note that an analogous result in the semi-simple setting has been stated without proof in \cite{AS}. The notation $v_{\widetilde{\gamma}}^{-1}$ does not appear in their work; instead, they express this element in a basis of characters, which is possible in the semi-simple case only.

\begin{corollaire}
For all $f \in \mathrm{MCG}(\Sigma_g^{\mathrm{o}})$, it holds $\widehat{f} \in \mathcal{L}_{g,0}^{\mathrm{inv}}(H)$.
\end{corollaire}
\debutDemo
Let $\gamma$ be a positively-oriented, non-separating simple loop. Then $\overset{I}{\widetilde{\gamma}}$ satisfies the fusion relation of $\mathcal{L}_{0,1}(H)$, and thus $j_{\widetilde{\gamma}}$ is a morphism of $H$-module-algebras (Lemma \ref{injectionFusion}). Hence, since $v^{-1} \in \mathcal{Z}(H) = \mathcal{L}_{0,1}^{\mathrm{inv}}(H)$, we have $v_{\widetilde{\gamma}}^{-1} \in \mathcal{L}_{g,0}^{\mathrm{inv}}(H)$. In particular, the statement is true for the Humphries generators thanks to Proposition \ref{propDehnTwist} and thus for any $f$.
\finDemo

\subsection{Representation of the mapping class group}\label{sectionRepMCG}
\indent The only additional fact needed is the following lemma.
\begin{lemme}\label{lemmevA1Moins1}
It holds: $v_{A(g)}^{-1} = v_{A(g)^{-1}}^{-1}$.
\end{lemme}
\debutDemo
Denote as usual $X_i \otimes Y_i = RR'$, $\overline{X}_i \otimes \overline{Y}_i = (RR')^{-1}$ and let $\mu^l$ be the left integral on $H$ (unique up to scalar). Using basic facts about integrals and \eqref{ribbon}, we have  (see \cite[Prop. 5.3]{Fai18}):
\[ \mu^l(vX_i)Y_i = \mu^l(v\overline{X}_i)\overline{Y}_i = \mu^l(v)v^{-1}. \]
Let us write $\mu^l(v)^{-1}\mu^l(v?) = \sum_{i,j,I} c_{I,i}^j \overset{I}{T}{^i_j}$ with $c_{I,i}^j \in \mathbb{C}$. Then, using the identification $\overset{I}{M} = (\overset{I}{X_i})Y_i$ between $\mathcal{L}_{0,1}(H)$ and $H$, the fact that $\overset{I}{M}{^{-1}} = (\overset{I}{\overline{X}_i})\overline{Y}_i$ and the equalities above, we get
\[ v_{A(g)}^{-1} = j_{A(g)}\!\left( \sum_{i,j,I} c_{I,i}^j \overset{I}{M}{^i_j} \right) = j_{A(g)}\!\left( \sum_{i,j,I} c_{I,i}^j (\overset{I}{M}{^{-1}})^i_j \right) = j_{A(g)^{-1}}\!\left( \sum_{i,j,I} c_{I,i}^j \overset{I}{M}{^i_j} \right) = v_{A(g)^{-1}}^{-1} \]
where the morphisms $j_{\bullet}$ are defined at the end of subsection \ref{defLgn}. We used that $j_{A(g)}$ is a morphism of algebras (see Lemma \ref{injectionFusion}).
\finDemo

It is clear that the lemma holds for the lift of any positively oriented, non-separating simple loop, but we do not need this.
\smallskip\\  
\indent Recall that we have a representation of $\mathcal{L}_{g,0}(H)$ on $(H^*)^{\otimes g}$, let us denote it $\rho$. We also have the induced representation of $\mathcal{L}^{\mathrm{inv}}_{g,0}(H)$ on $\mathrm{Inv}\!\left((H^*)^{\otimes g}\right)$, let us denote it $\rho_{\mathrm{inv}}$. Also recall that the elements $\widehat{f}$ are defined in \eqref{conjugaison}. We can now state the representation of the mapping class groups $\mathrm{MCG}(\Sigma_g^{\mathrm{o}})$ and $\mathrm{MCG}(\Sigma_g)$. An analogous result in the semi-simple setting has been given without proof in \cite{AS}.

\begin{theoreme}\label{thmRepMCG}
1) The map
\[ \fleche{\mathrm{MCG}(\Sigma_g^{\mathrm{o}})}{\mathrm{GL}\!\left((H^*)^{\otimes g}\right)}{f}{\rho(\widehat{f})} \]
is a projective representation.
\\2) The map
\[ \fleche{\mathrm{MCG}(\Sigma_g)}{\mathrm{GL}\!\left(\mathrm{Inv}\!\left((H^*)^{\otimes g}\right)\right)}{f}{\rho_{\mathrm{inv}}(\widehat{f})}\]
is a projective representation.
\end{theoreme}
\debutDemo
1) This is an immediate consequence of Proposition \ref{liftHumphries}.
\\2) We must show that the hyperelliptic relation \eqref{hyperelliptic} is projectively satisfied. The word $\omega$ can be constructed as follows: take $f \in \mathrm{MCG}(\Sigma_g^{\mathrm{o}})$ such that $f(a_1) = a_g$ and express it as a word in the Humphries generators $f = \tau_{\gamma_1} \ldots \tau_{\gamma_n}$. Then $\tau_{a_g} = f \tau_{a_1} f^{-1}$. The automorphism $\widetilde{\tau_{a_g}}$ is implemented by conjugation by $\widehat{f} v_{A(1)}^{-1} \widehat{f}^{-1}$ and also by conjugation by $v_{A(g)}^{-1}$ (Proposition \ref{propDehnTwist}). Hence, $\widehat{f} v_{A(1)}^{-1} \widehat{f}^{-1} \sim v_{A(g)}^{-1}$, where $\sim$ means proportional. Now, let $H = \tau_{b_g} \tau_{d_g} \ldots \tau_{b_1} \tau_{d_1} \tau_{d_1} \tau_{b_1} \ldots \tau_{d_g} \tau_{b_g}$. A computation gives $\widetilde{H}(\overset{I}{A}(g)) = \overset{I}{A}(g){^{-1}} \overset{I}{C}_{g,0}$. Thus
\[ \widehat{H} \widehat{f} v_{A(1)}^{-1} \widehat{f}^{-1} \widehat{H}^{-1} \sim \widehat{H} v_{A(g)}^{-1} \widehat{H}^{-1} = \widetilde{H}(v_{A(g)}^{-1}) = v_{A(g)^{-1}C_{g,0}}^{-1}. \]
By definition of $\mathrm{Inv}\!\left((H^*)^{\otimes g}\right)$ and Lemma \ref{lemmevA1Moins1}, we have 
\[ \rho_{\mathrm{inv}}(v_{A(g)^{-1}C_{g,0}}^{-1}) = \rho_{\mathrm{inv}}(v_{A(g)^{-1}}^{-1}) = \rho_{\mathrm{inv}}(v_{A(g)}^{-1}). \]
It follows that $\rho_{\mathrm{inv}}\!\left(\widehat{H} \left(\widehat{f} v_{A(1)}^{-1} \widehat{f}^{-1}\right) \widehat{H}^{-1}\right) \sim \rho_{\mathrm{inv}}\!\left(\widehat{f} v_{A(1)}^{-1} \widehat{f}^{-1}\right)$. This shows that the map is well-defined since $\mathrm{MCG}(\Sigma_g)$ is the quotient of $\mathrm{MCG}(\Sigma_g^{\mathrm{o}})$ by the hyperelliptic relation and that it is a projective representation.
\finDemo
%\smallskip\\
%\indent In \cite{lyu95b}, Lyubashenko has constructed a projective representation of $\mathrm{MCG}(\Sigma_{g,n}^{\mathrm{o}})$ and $\mathrm{MCG}(\Sigma_{g,n})$ by completely different techniques based on the coend of a ribbon category satisfying some assumptions. Our assumptions on $H$ allow to apply his construction to $\mathrm{mod}_l(H)$, the category of finite-dimensional left $H$-modules. A natural conjecture is that the representation constructed here is equivalent to Lyubashenko's representation for $\mathrm{mod}_l(H)$. By \cite{Fai18}, we know that the restrictions of these representations agree on each of the $g$ subtori of $\Sigma_g$. Showing the equivalence on the whole surface is technically more difficult and is the object of a work in progress.

\subsection{Discussion for the case $n > 0$}\label{CasGeneral}
\indent Let us consider the general case $n>0$, see Figure \ref{figureIntro}. Denote $\Sigma_{g,n}^{\mathrm{o}} = \Sigma_{g,n} \setminus D$, where $D$ is an embedded open disk. Recall that by definition the mapping class group fixes pointwise the
boundary.
\\
\indent In general, $\mathcal{L}_{g,n}(H)$ is not a matrix algebra and we cannot claim directly the existence and unicity up to scalar of the elements $\widehat{f}$. Nevertheless, we now describe an extension of the previous construction which should not be difficult to apply.
\smallskip\\
\indent \textbullet ~ Consider a generating set $\tau_{c_1}, \ldots, \tau_{c_k}$ of $\mathrm{MCG}(\Sigma^{\mathrm{o}}_{g,n})$ (this consists only of Dehn twists when $g > 1$, see \cite[Figure 4.10]{FM}) and compute the action on $\pi_1(\Sigma_{g,n}^{\mathrm{o}})$. Here the $c_i$ are loops in $\pi_1(\Sigma_{g,n}^{\mathrm{o}})$ written in terms of the generators depicted in Figure \ref{figureIntro}.
\smallskip\\
\indent  \textbullet ~ Determine the lifts $\widetilde{c_i}$ of the loops $c_i$ (\textit{i.e.} replace generators of $\pi_1(\Sigma_{g,n}^{\mathrm{o}})$ by matrices of generators of $\mathcal{L}_{g,n}(H)$ and then determine the normalisations by powers of $v$ needed to satisfy the fusion relation).
\smallskip\\
\indent \textbullet ~ Determine the lifts $\widetilde{\tau_{c_i}}$ of the generators as automorphisms of $\mathcal{L}_{g,n}(H)$ (\textit{i.e.} replace generators of $\pi_1(\Sigma_{g,n}^{\mathrm{o}})$ by matrices of generators of $\mathcal{L}_{g,n}(H)$ in the action of $\tau_{c_1}, \ldots, \tau_{c_k}$ on $\pi_1(\Sigma_{g,n}^{\mathrm{o}})$, and then determine the normalisations by powers of $v$ needed to satisfy the fusion relation and check that the other relations of \eqref{PresentationLgn} are satisfied).
\smallskip\\
\indent \textbullet ~ Show that the assignment $\tau_{c_k} \mapsto \widetilde{\tau_{c_k}}$ extends to a morphism of groups $\mathrm{MCG}(\Sigma_{g,n}^{\mathrm{o}}) \to \mathrm{Aut}(\mathcal{L}_{g,n}(H))$ (this is a just a tedious verification using a presentation of $\mathrm{MCG}(\Sigma_{g,n}^{\mathrm{o}})$). Thus the lift $\widetilde{f}$ of $f \in \mathrm{MCG}(\Sigma_{g,n}^{\mathrm{o}})$ is still defined.
\smallskip\\
\indent \textbullet ~ It is clear that Lemma \ref{transfoLoop} still holds, so in particular for each $i$ there exists $f_i \in \mathrm{MCG}(\Sigma^{\mathrm{o}}_{g,n})$ such that $f_i(a_1) = c_i$.
\smallskip\\
\indent \textbullet ~ It is clear that Lemma \ref{lemmeA1} still holds, so that $\widetilde{\tau}_{c_i}(x) = v_{\widetilde{c_i}}^{-1} x v_{\widetilde{c_i}}$ (by reproducing the proof of Proposition \ref{propDehnTwist} with the $f_i$). Since the $\tau_{c_i}$ are a generating set, it follows that for each $f \in \mathrm{MCG}(\Sigma_{g,n}^{\mathrm{o}})$, there exists an element $\widehat{f} \in \mathcal{L}^{\mathrm{inv}}_{g,n}(H)$, unique up to an invertible central element such that $\widetilde{f}(x) = \widehat{f} x \widehat{f}^{-1}$.
\smallskip\\
\indent \textbullet ~ Since $\mathcal{Z}\!\left(\mathcal{H}({\mathcal{O}}(H))^{\otimes g}\right) \cong \mathbb{C}$, we have for all $c \in \mathcal{Z}\!\left(\mathcal{L}_{g,n}(H)\right)$: 
\[\Psi_{g,n}(c) = 1 \otimes \ldots \otimes 1 \otimes c_1 \otimes \ldots \otimes c_n\]
with $c_i \in \mathcal{Z}(H)$. Let $V = (H^*)^{\otimes g} \otimes S_1 \otimes \ldots \otimes S_n$, where $S_1, \ldots, S_n$ are simple representations of $H$. Then $\Psi_{g,n}(c)$ acts by scalar on $V$ thanks to Schur lemma. Let $\rho$ (resp. $\rho_{\mathrm{inv}}$) be the representation of $\mathcal{L}_{g,n}(H)$ (resp. $\mathcal{L}_{g,n}^{\mathrm{inv}}(H)$) on $V$ (resp. $\mathrm{Inv}(V)$). Then the elements $\rho(\widehat{f})$ and thus $\rho_{\mathrm{inv}}(\widehat{f})$ are unique up to scalar. It should not be difficult to check that the corresponding generalisation of Theorem \ref{thmRepMCG} is true.

\subsection{Explicit formulas for the representation of some Dehn twists}\label{sectionFormulesExplicites}
\indent We will compute explicitly the representation on $(H^*)^{\otimes g}$ of the Dehn twists $\tau_{\gamma}$, where the curves $\gamma$ are represented in Figure \ref{figureCourbesCanoniques}. Thanks to Proposition \ref{propDehnTwist}, this amounts to compute the action of $v_{\widetilde{\gamma}}^{-1}$ on $(H^*)^{\otimes g}$.
\smallskip\\
\indent We recall that the action $\triangleright$ of $\mathcal{L}_{g,0}(H)$ on $(H^*)^{\otimes g}$ is defined using $\Psi_{g,0}$ in \eqref{actionTriangle} and that we denote the associated representation by $\rho$. Also recall the definition of the elements $\widetilde{h}$ in \eqref{operateursTilde} and the notation  $RR' = X_i \otimes Y_i$. Note that
\begin{equation}\label{coproduitRR}
X_i \otimes Y_i' \otimes Y_i'' = a_j X_i b_k \otimes Y_i \otimes b_j a_k.
\end{equation}

\indent Recall from \cite{Fai18} the elements $v_A^{-1}, v_B^{-1} \in \mathcal{L}_{1,0}(H)$ and their action on $H^*$:
\begin{equation}\label{actionvAvB}
\begin{split}
 v_A^{-1} \triangleright \varphi &= \varphi^{v^{-1}} = \varphi(v^{-1} ?),\\
 v_B^{-1} \triangleright \varphi &= \mu^l(v)^{-1}\left(\mu^l\!\left(g^{-1}v\,?\right) \varphi^v\right)^{v^{-1}} %= \mu^l(v)^{-1}\mu^l\!\left(g^{-1}v^{-1}v'\,?\right)\varphi(g^{-1}vS(v'')) 
\end{split}
\end{equation}
where $\varphi^h = \varphi(h?)$ for $h \in H$ and $\mu^l$ is the left integral on $H$. 
\smallskip\\
We will need the following generalization of \cite[Lemma 5.7]{Fai18} (in which we restricted to $\varphi \in \mathrm{SLF}(H)$).

\begin{lemme}\label{lemmeActionTresseAB}
For all $\varphi \in H^*$:
\begin{align*}
\left( v_{A}^{-1} v_{B}^{-1} v_{A}^{-1} \right)^2 \triangleright \varphi &= \frac{\mu^l(v^{-1})}{\mu^l(v)} \varphi\!\left( S^{-1}(a_i) g^{-1}v^{-1} S(?) b_i \right)\\
\left( v_{A}^{-1} v_{B}^{-1} v_{A}^{-1} \right)^{-2} \triangleright \varphi &= \frac{\mu^l(v)}{\mu^l(v^{-1})} \varphi\!\left( b_j S^{-1}(?) a_j g^{-1}v \right)
\end{align*}
\end{lemme}
\debutDemo
Write $\varphi = \sum_{I,i,j} \Phi_{I,i}^j \overset{I}{T}{_j^i} = \sum_I \mathrm{tr}\!\left( \Phi_I \overset{I}{T} \right)$ with $\Phi_{I,i}^j \in \mathbb{C}$ and let $z(\varphi) = \sum_I \mathrm{tr}\!\left( \overset{I}{b_i} \Phi_I \overset{I}{S^{-1}(a_i)} \overset{I}{M} \right)$ $\in \mathcal{L}_{0,1}(H)$. Then $z(\varphi)_B \triangleright \varepsilon = \varphi$ (where $z(\varphi)_B = j_B(z(\varphi))$, see notation at the end of section \ref{defLgn}, and $\varepsilon$ is the counit of $H$). Indeed
\[ z(\varphi)_B \triangleright \varepsilon 
%= \sum_I \mathrm{tr}\!\left( \overset{I}{b_i} \Phi_I \overset{I}{S^{-1}(a_i)} \Psi_{1,0}(\overset{I}{M}) \right) \triangleright \varepsilon 
= \sum_I \mathrm{tr}\!\left( \overset{I}{b_i} \Phi_I \overset{I}{S^{-1}(a_i)} \overset{I}{L}{^{(+)}} \overset{I}{T} \overset{I}{L}{^{(-)-1}} \triangleright \varepsilon \right)
= \sum_I \mathrm{tr}\!\left( \overset{I}{b_i} \Phi_I \overset{I}{S^{-1}(a_i)} \overset{I}{a_j} \overset{I}{T} \overset{I}{b_j}  \right) = \sum_I \mathrm{tr}\!\left(\Phi_I  \overset{I}{T} \right) = \varphi.\]
We simply used \eqref{actionTriangle}, \eqref{repTriangleHeisenberg}, the cyclicity of the trace and the equality $S^{-1}(a_i)a_j \otimes b_j b_i = 1 \otimes 1$. Observe that 
\[ \left(\widetilde{\tau}_a \widetilde{\tau}_b \widetilde{\tau}_a\right)^2 (\overset{I}{B}) = \overset{I}{v}{^2}\overset{I}{A}{^{-1}}\overset{I}{B}{^{-1}}\overset{I}{A} = \overset{I}{B}{^{-1}} \overset{I}{C} \]
where $\overset{I}{C} = \overset{I}{C}_{1,0}$ is defined in \eqref{matriceCL10}. Hence:
\begin{align*}
\left( v_{A}^{-1} v_{B}^{-1} v_{A}^{-1} \right)^2 \triangleright \varphi &= \left( v_{A}^{-1} v_{B}^{-1} v_{A}^{-1} \right)^2 z(\varphi)_{B} \triangleright \varepsilon = z(\varphi)_{B^{-1}C} \left( v_{A}^{-1} v_{B}^{-1} v_{A}^{-1} \right)^2 \triangleright \varepsilon\\
&= \frac{\mu^l(v^{-1})}{\mu^l(v)} z(\varphi)_{B^{-1}C} \triangleright \varepsilon = \frac{\mu^l(v^{-1})}{\mu^l(v)} z(\varphi)_{B^{-1}} \triangleright \varepsilon.
\end{align*}
We used Proposition \ref{propDehnTwist}, the formula of \cite[Lemma 5.7]{Fai18} applied to $\varepsilon$, and the fact that $\overset{I}{C} \triangleright \varepsilon = \mathbb{I}_{\dim(I)}\varepsilon$ (which follows from \ref{actionH}). Now we compute
\begin{align*}
z(\varphi)_{B^{-1}} \triangleright \varepsilon &= \sum_I \mathrm{tr}\!\left( \overset{I}{b_i} \Phi_I \overset{I}{S^{-1}(a_i)} \overset{I}{L}{^{(-)}} S(\overset{I}{T}) \overset{I}{L}{^{(+)-1}} \triangleright \varepsilon \right)
= \sum_I \mathrm{tr}\!\left( \overset{I}{b_i} \Phi_I \overset{I}{S^{-1}(a_i)} \overset{I}{S^{-1}(b_j)}a_j \triangleright S(\overset{I}{T}) \right)\\
&= \sum_I \mathrm{tr}\!\left( \overset{I}{b_i} \Phi_I \overset{I}{S^{-1}(a_i)} \overset{I}{S^{-1}(b_j)} \overset{I}{S(a_j)} S(\overset{I}{T}) \right) = \sum_I \mathrm{tr}\!\left( \Phi_I \overset{I}{S^{-1}(a_i)} \overset{I}{g}{^{-1}} \overset{I}{v}{^{-1}} S(\overset{I}{T}) \overset{I}{b_i} \right)\\
&= \varphi\!\left( S^{-1}(a_i) g^{-1}v^{-1} S(?) b_i \right).
\end{align*}
We used \eqref{repHO} and \eqref{u}. The second formula is easily checked.
\finDemo

\begin{theoreme}\label{formulesExplicites}
The following formulas hold:
\begin{align*}
v_{A(i)}^{-1} \triangleright \left(\varphi_1 \otimes \ldots \otimes \varphi_g\right) &= \varphi_1 \otimes \ldots \otimes \varphi_{i-1} \otimes (v_A^{-1} \triangleright \varphi_i) \otimes \varphi_{i+1} \otimes \ldots \otimes \varphi_g, \\
v_{B(i)}^{-1} \triangleright \left(\varphi_1 \otimes \ldots \otimes \varphi_g\right) &= \varphi_1 \otimes \ldots \otimes \varphi_{i-1} \otimes (v_B^{-1} \triangleright \varphi_i) \otimes \varphi_{i+1} \otimes \ldots \otimes \varphi_g, \\
v_{D_i}^{-1} \triangleright \left(\varphi_1 \otimes \ldots \otimes \varphi_g\right) &= \varphi_1 \otimes \ldots \otimes \varphi_{i-2} \otimes \varphi_{i-1}\!\left(S^{-1}(a_j)a_k?b_k v''^{-1} b_j\right) \otimes \varphi_i\!\left( S^{-1}(a_l)  S^{-1}(v'^{-1}) a_m ? b_m b_l \right)\\
& \:\:\:\:\: \otimes \varphi_{i+1} \otimes \ldots \otimes \varphi_g, \\
v_{E_i}^{-1} \triangleright \left(\varphi_1 \otimes \ldots \otimes \varphi_g\right) &= \varphi_1\!\left(S^{-1}\!\left(v^{(2i-2)-1}\right) ? v^{(2i-1)-1} \right) \otimes \ldots \otimes \varphi_{i-1}\!\left(S^{-1}\!\left( v^{(2)-1} \right) ?  v^{(3)-1} \right)\\
&\:\:\:\:\: \otimes \varphi_i\!\left(S^{-1}(a_j)  S^{-1}\!\left(v^{(1)-1}\right) a_k ? b_k b_j \right) \otimes \varphi_{i+1} \otimes \ldots \otimes \varphi_g,
\end{align*}
with $i \geq 2$ for the two last formulas.
\end{theoreme}

\indent The rest of the section is devoted to the proof of that theorem. First, it is useful to record that
\begin{equation}\label{formuleLambda}
\begin{split}
\Psi_{1,0}^{\otimes g}(\overset{I}{\Lambda}_i) = \Psi_{1,0}^{\otimes g}\!\left( \overset{I}{\underline{C}}{^{(-)}}(1) \ldots \overset{I}{\underline{C}}{^{(-)}}(i-1) \right) &= \overset{I}{S^{-1}(b_j)} \: \widetilde{a_j^{(2i-3)}} a_j^{(2i-2)}  \otimes  \ldots  \otimes  \widetilde{a_j^{(1)}}a_j^{(2)}  \otimes  1^{\otimes g-i+1} \\
\Psi_{1,0}^{\otimes g}\!\left( \overset{I}{\underline{C}}{^{(+)}}(1) \ldots \overset{I}{\underline{C}}{^{(+)}}(i-1) \right) &= \overset{I}{a_j} \: \widetilde{b_j^{(2i-3)}} b_j^{(2i-2)}  \otimes  \ldots  \otimes  \widetilde{b_j^{(1)}}b_j^{(2)}  \otimes  1^{\otimes g-i+1}
\end{split}
\end{equation}
where the matrix $\overset{I}{\Lambda}_k$ is defined in \eqref{matricesAlekseev}. The proof is a simple computation analogous to that of Lemma \ref{expressionM}. Second, recall from the proof of Lemma \ref{lemmevA1Moins1} that 
\begin{equation}\label{VIntegrale}
\mu^l(v)^{-1}\mu^l(vX_i)Y_i = v^{-1}.
\end{equation}
We will write $\mu^l(v)^{-1}\mu^l(v?) = \sum_{I} \mathrm{tr}\!\left(c_I \overset{I}{T}\right)$. Then $v^{-1} = \sum_{I} \mathrm{tr}\!\left(c_I \overset{I}{M}\right)$ under the identification $\mathcal{L}_{0,1}(H) = H$.
\smallskip\\
\indent \textbullet ~ {\em Proof of the formula for the action of $v_{A(i)}^{-1}$.} By definition and by \eqref{formuleLambda}, we have
\begin{align*}
&\Psi_{g,0}(v_{A(i)}^{-1}) = \sum_{I} \mathrm{tr}\!\left(c_I \Psi_{1,0}^{\otimes g}(\overset{I}{\Lambda}_i \, \overset{I}{\underline{A}}(i) \, \overset{I}{\Lambda}_i{^{-1}})\right)\\
& = \sum_I \mathrm{tr}\!\left(c_I \overset{I}{S^{-1}(b_j)} \overset{I}{X_k} \overset{I}{b_l}\right)  \widetilde{a_j^{(2i-3)}}\widetilde{a_l^{(2i-3)}} a_j^{(2i-2)} a_l^{(2i-2)}  \otimes  \ldots  \otimes  \widetilde{a_j^{(1)}}\widetilde{a_l^{(1)}} a_j^{(2)}a_l^{(2)}  \otimes  Y_k  \otimes  1^{\otimes g-i}\\
&=\mu^l(v)^{-1}\mu^l\!\left( v S^{-1}(b_j) X_k b_l\right)  \widetilde{a_j^{(2i-3)}}\widetilde{a_l^{(2i-3)}} a_j^{(2i-2)} a_l^{(2i-2)}  \otimes  \ldots  \otimes  \widetilde{a_j^{(1)}}\widetilde{a_l^{(1)}} a_j^{(2)}a_l^{(2)}  \otimes  Y_k  \otimes  1^{\otimes g-i}\\
&= \mu^l(v)^{-1}\mu^l\!\left(v S^{-1}\!\left(b_jS^{-1}(b_l)\right) X_k\right)  \widetilde{a_j^{(2i-3)}}\widetilde{a_l^{(2i-3)}} a_j^{(2i-2)} a_l^{(2i-2)}  \otimes  \ldots  \otimes  \widetilde{a_j^{(1)}}\widetilde{a_l^{(1)}} a_j^{(2)}a_l^{(2)}  \otimes  Y_k  \otimes  1^{\otimes g-i}\\
&= \mu^l(v)^{-1}\mu^l\!\left( v X_k\right)  1^{\otimes i-1}  \otimes  Y_k  \otimes  1^{\otimes g-i} = 1^{\otimes i-1}  \otimes  v^{-1}  \otimes  1^{\otimes g-i}
\end{align*}
and the formula follows. We used \eqref{quasiCyclic}, the formula $R^{-1} = a_l \otimes S^{-1}(b_l)$ and \eqref{VIntegrale}.
\smallskip\\
\indent \textbullet ~ {\em Proof of the formula for the action of $v_{B(i)}^{-1}$.} This the same proof as for $v_{A(i)}^{-1}$ (the conjugation by $\overset{I}{\Lambda}_i$ vanishes thanks to \eqref{quasiCyclic}).
\smallskip\\
\indent \textbullet ~ {\em Proof of the formula for the action of $v_{D_i}^{-1}$, $i \geq 2$.} We first compute the action of $\overset{I}{A}(i-1)\overset{I}{A}(i)$. We have
\begin{align*}
\Psi_{g,0}\!\left( \overset{I}{A}(i-1)\overset{I}{A}(i) \right) &= \Psi_{1,0}^{\otimes g}\!\left(\overset{I}{\Lambda}_{i-1} \, \overset{I}{\underline{A}}(i-1) \, \overset{I}{\Lambda}{_{i-1}^{-1}} \, \overset{I}{\Lambda}_{i} \, \overset{I}{\underline{A}}(i) \, \overset{I}{\Lambda}{_{i}^{-1}}\right)\\
&= \Psi_{1,0}^{\otimes g}\!\left(\overset{I}{\Lambda}_{i-1} \, \overset{I}{\underline{A}}(i-1) \, \overset{I}{\underline{C}}{^{(-)}}(i-1) \, \overset{I}{\underline{A}}(i) \, \overset{I}{\underline{C}}{^{(-)}}(i-1)^{-1} \, \overset{I}{\Lambda}{_{i-1}^{-1}}\right).
\end{align*}
Hence:
\begin{align*}
&\Psi_{g,0}\!\left( v_{A(i-1)A(i)}^{-1} \right) %&= \sum_{I} \mathrm{tr}\!\left(c_I \Psi_{1,0}^{\otimes g}(\overset{I}{\Lambda}_{i-1} \, \overset{I}{\underline{A}}(i-1) \, \overset{I}{\Lambda}{_{i-1}^{-1}} \, \overset{I}{\Lambda}_{i} \, \overset{I}{\underline{A}}(i) \, \overset{I}{\Lambda}{_{i}^{-1}})\right)
= \sum_{I} \mathrm{tr}\!\left(c_I \Psi_{1,0}^{\otimes g}\!\left(\overset{I}{\Lambda}_{i-1} \, \overset{I}{\underline{A}}(i-1) \, \overset{I}{\underline{C}}{^{(-)}}(i-1) \, \overset{I}{\underline{A}}(i) \, \overset{I}{\underline{C}}{^{(-)}}(i-1)^{-1} \, \overset{I}{\Lambda}{_{i-1}^{-1}}\right)\right)\\
%&= \sum_{I} \mathrm{tr}\!\left(c_I \Psi_{1,0}^{\otimes g} (\overset{I}{\underline{A}}(i-1) \, \overset{I}{C}(i)^{(-)} \, \overset{I}{\underline{A}}(i) \, \overset{I}{C}(i)^{(-)-1})\right)\\
%&= \sum_{I} \mathrm{tr}\!\left(c_I X_j S^{-1}(b_k) X_l b_m \right) 1^{\otimes i-2} \otimes \widetilde{a_k'}\widetilde{a_m'} Y_ja_k'' a_m'' \otimes Y_l \otimes 1^{\otimes g-i}\\
&= \mu^l(v)^{-1}\mu^l\!\left( v S^{-1}(b_j) X_k S^{-1}(b_l) X_m b_n b_o \right)  \widetilde{a_j^{(2i-5)}}\widetilde{a_o^{(2i-5)}} a_j^{(2i-4)} a_o^{(2i-4)} \, \otimes \, \ldots \, \otimes \, \widetilde{a_j^{(1)}}\widetilde{a_o^{(1)}} a_j^{(2)}a_o^{(2)}\\
& \:\:\:\:\:\:\:\:\:\:\:\:\:\:\:\:\:\:\:\:\:\:\:\:\:\:\:\:\:\:\:\:\:\:\:\:\:\:\:\:\:\:\:\:\:\:\:\:\:\:\:\:\:\:\:\:\:\:\:\:\:\:\:\:\:\:\:\:\:\:\:\:\:\:\:\otimes \widetilde{a_l'}\widetilde{a_n'} Y_k a_l'' a_n'' \otimes Y_m \otimes 1^{\otimes g-i}\\
&= \mu^l(v)^{-1}\mu^l\!\left( v X_k S^{-1}(b_l) X_m b_n \right)  1^{\otimes i-2} \otimes \widetilde{a_l'}\widetilde{a_n'} Y_k a_l'' a_n'' \otimes Y_m \otimes 1^{\otimes g-i}\\
&= \mu^l(v)^{-1}\mu^l\!\left( v a_k S^{-1}(b_l) X_m b_n \right)  1^{\otimes i-2} \otimes \widetilde{a_l}\widetilde{a_n'} b_k a_n'' \otimes Y_m \otimes 1^{\otimes g-i}\\
\end{align*}
We used \eqref{quasiCyclic} and the fact that $X_k S^{-1}(b_l) \otimes a_l' \otimes Y_k a_l'' = X_k S^{-1}(b_p)S^{-1}(b_l) \otimes a_l \otimes Y_k a_p = a_kS^{-1}(b_l) \otimes a_l \otimes b_k$. We see that we can assume without loss of generality that $g=2$, $i=2$ since the action is ``local''. Moreover, this can be simplified. Let $F : H^* \to H^*$ be the map defined by
\[ F(\varphi) = \varphi\!\left(a_j ? b_j\right), \:\:\:\:\: F^{-1}(\varphi) = \varphi\!\left(S^{-1}(a_j) ? b_j\right). \]
We compute:
\begin{align*}
&(F^{-1} \otimes \mathrm{id}) \circ \rho\!\left( v_{A(1)A(2)}^{-1} \right) \circ (F \otimes \mathrm{id})(\varphi \otimes \psi)\\
&= \mu^l(v)^{-1}\mu^l\!\left( v a_k b_l X_m b_n \right)   \varphi\!\left( a_j S^{-1}(a_n')a_l S^{-1}(a_o) ? b_o b_k a_n''b_j \right) \otimes \psi(?Y_m)\\
&= \mu^l(v)^{-1}\mu^l\!\left(v a_k b_l X_m b_n \right)   \varphi\!\left( S^{-1}(a_n'') a_j a_l S^{-1}(a_o) ? b_o b_k b_j a_n' \right) \otimes \psi(?Y_m) = (\star).
\end{align*}
We used the formula $R\Delta = \Delta^{\mathrm{op}}R$. Now, we have a Yang-Baxter identity 
\[ a_k b_l \otimes a_j a_l \otimes b_k b_j = R_{13} R_{23} R_{21} = R_{21} R_{23} R_{13} = b_l a_k \otimes a_l a_j \otimes b_j b_k \]
which allows us to continue the computation:
\begin{align*}
(\star) &= \mu^l(v)^{-1}\mu^l\!\left(v b_l a_k X_m b_n \right)   \varphi\!\left( S^{-1}(a_n'') a_l a_j S^{-1}(a_o) ? b_o b_j b_k a_n' \right) \otimes \psi(?Y_m)\\
&= \mu^l(v)^{-1}\mu^l\!\left(v S^{-2}(b_p) b_l a_k X_m b_n \right)   \varphi\!\left( S^{-1}(a_p) a_l ?  b_k a_n \right) \otimes \psi(?Y_m)\\
&= \mu^l(v)^{-1}\mu^l\!\left(v a_k X_m b_n \right)   \varphi\!\left( ?  b_k a_n \right) \otimes \psi(?Y_m)\\
&= \mu^l(v)^{-1}\mu^l\!\left(v X_m \right)   \varphi\!\left( ?  Y_m'' \right) \otimes \psi(? Y_m') = \varphi\!\left( ?  v''^{-1} \right) \otimes \psi(? v'^{-1}).
\end{align*}
We used basic properties of the $R$-matrix and relations \eqref{coproduitRR}, \eqref{VIntegrale}. We have thus shown that
\[ v_{A(1)A(2)}^{-1} \triangleright \varphi \otimes \psi = \varphi\!\left( S^{-1}(a_j)a_k ? b_k v''^{-1} b_j \right) \otimes \psi(?v'^{-1}). \]
Recall that $\overset{I}{D}_2 = \overset{I}{v}{^2} \overset{I}{A}(1)\overset{I}{B}(2)\overset{I}{A}(2){^{-1}} \overset{I}{B}(2){^{-1}}$. Hence $\left(\widetilde{\tau}_{a_2} \widetilde{\tau}_{b_2} \widetilde{\tau}_{a_2} \right)^{-2}(\overset{I}{A}(1)\overset{I}{A}(2)) = \overset{I}{D}_2$. It follows that $\left(\widetilde{\tau}_{a_2} \widetilde{\tau}_{b_2} \widetilde{\tau}_{a_2} \right)^{-2}(v_{A(1)A(2)}^{-1}) = v_{D_2}^{-1}$, and thus by Proposition \ref{propDehnTwist} and Lemma \ref{lemmeActionTresseAB}:
\begin{align*}
v_{D_2}^{-1} \triangleright \varphi \otimes \psi &= \left( v_{A(2)}^{-1} v_{B(2)}^{-1} v_{A(2)}^{-1} \right)^{-2} v_{A(1)A(2)}^{-1} \left( v_{A(2)}^{-1} v_{B(2)}^{-1} v_{A(2)}^{-1} \right)^2 \triangleright \varphi \otimes \psi\\
&= \left( v_{A(2)}^{-1} v_{B(2)}^{-1} v_{A(2)}^{-1} \right)^{-2} \triangleright \varphi\!\left( S^{-1}(a_j)a_k ? b_k v''^{-1} b_j \right) \otimes \psi\!\left( S^{-1}(a_l) g^{-1} v^{-1}  S(v'^{-1}) S(?) b_l \right)\\
&= \varphi\!\left( S^{-1}(a_j)a_k ? b_k v''^{-1} b_j \right) \otimes \psi\!\left( S^{-1}(a_l) g^{-1} v^{-1}  S(v'^{-1}) S\!\left(b_m S^{-1}(?) a_m g^{-1}v\right) b_l \right)\\
&= \varphi\!\left( S^{-1}(a_j)a_k ? b_k v''^{-1} b_j \right) \otimes \psi\!\left( S^{-1}(a_l)  S^{-1}(v'^{-1}) a_m ? b_mb_l \right)
\end{align*}
which is the announced formula.
\smallskip\\
\indent \textbullet ~ {\em Proof of the formula for the action of $v_{E_i}^{-1}$, $i \geq 2$.} We first compute the action of $\overset{I}{C}(1) \ldots \overset{I}{C}(i-1) \overset{I}{A}(i)$. We have
\begin{align*}
\Psi_{g,0}\!\left(\overset{I}{C}(1) \ldots \overset{I}{C}(i-1) \overset{I}{A}(i) \right) &= \Psi_{1,0}^{\otimes g}\!\left( \overset{I}{\underline{C}}{^{(+)}}(1) \ldots \overset{I}{\underline{C}}{^{(+)}}(i-1) \left(\overset{I}{\underline{C}}{^{(-)}}(1) \ldots \overset{I}{\underline{C}}{^{(-)}}(i-1)\right)^{-1} \overset{I} \Lambda_i \overset{I}{\underline{A}}(i) \Lambda_i^{-1} \right)\\
&= \Psi_{1,0}^{\otimes g}\!\left( \overset{I}{\underline{C}}{^{(+)}}(1) \ldots \overset{I}{\underline{C}}{^{(+)}}(i-1) \overset{I}{\underline{A}}(i) \left(\overset{I}{\underline{C}}{^{(-)}}(1) \ldots \overset{I}{\underline{C}}{^{(-)}}(i-1)\right)^{-1} \right)\\
&= \overset{I}{a_j} \overset{I}{X_k} \overset{I}{b_l} \, \widetilde{b_j^{(2i-3)}} \widetilde{a_l^{(2i-3)}} b_j^{(2i-2)} a_l^{(2i-2)}  \otimes  \ldots  \otimes  \widetilde{b_j^{(1)}} \widetilde{a_l^{(1)}} b_j^{(2)} a_l^{(2)}  \otimes  Y_k \otimes 1^{\otimes g-i}\\
&= \overset{I}{X_k} \, \widetilde{Y_k^{(2i-2)}}  Y_k^{(2i-1)} \otimes \ldots \otimes \widetilde{Y_k^{(2)}}  Y_k^{(3)}   \otimes  Y_k^{(1)} \otimes 1^{\otimes g-i}
\end{align*}
thanks to \eqref{formuleLambda} and \eqref{coproduitRR}. Hence, by \eqref{VIntegrale}:
\begin{align*}
\Psi_{g,0}\!\left(v_{C(1) \ldots C(i-1) A(i)}^{-1} \right) &= \mu^l(v)^{-1}\mu^l(v X_k)\, \widetilde{Y_k^{(2i-2)}}  Y_k^{(2i-1)} \otimes \ldots \otimes \widetilde{Y_k^{(2)}}  Y_k^{(3)}   \otimes  Y_k^{(1)} \otimes 1^{\otimes g-i}\\
&= \widetilde{v^{(2i-2)-1}}  v^{(2i-1)-1} \otimes \ldots \otimes \widetilde{v^{(2)-1}}  v^{(3)-1}   \otimes  v^{(1)-1} \otimes 1^{\otimes g-i},
\end{align*}
which means that
\begin{align*}
v_{C(1) \ldots C(i-1) A(i)}^{-1} \triangleright \left(\varphi_1 \otimes \ldots \otimes \varphi_g\right) &= \varphi_1\!\left(S^{-1}\!\left(v^{(2i-2)-1}\right) ? v^{(2i-1)-1} \right) \otimes \ldots \otimes \varphi_{i-1}\!\left(S^{-1}\!\left( v^{(2)-1} \right) ?  v^{(3)-1} \right)\\
&\:\:\:\:\: \otimes \varphi_i\!\left( ? v^{(1)-1} \right) \otimes \varphi_{i+1} \ldots \otimes \varphi_g. 
\end{align*}
Recall that $\overset{I}{E}_i = \overset{I}{v}{^2} \overset{I}{C}(1) \ldots \overset{I}{C}(i-1)\overset{I}{B}(i)\overset{I}{A}(i){^{-1}} \overset{I}{B}(i){^{-1}}$. Hence $\left(\widetilde{\tau}_{a_i} \widetilde{\tau}_{b_i} \widetilde{\tau}_{a_i} \right)^{-2}(\overset{I}{C}(1) \ldots \overset{I}{C}(i-1)\overset{I}{A}(i)) = \overset{I}{E}_i$. As previously, it follows from Proposition \ref{propDehnTwist} and Lemma \ref{lemmeActionTresseAB} that
\begin{align*}
v_{E_i}^{-1} \triangleright \left( \varphi_1 \otimes \ldots \otimes \varphi_g \right) &= \left( v_{A(i)}^{-1} v_{B(i)}^{-1} v_{A(i)}^{-1} \right)^{-2} v_{C(1) \ldots C(i-1)A(i)}^{-1} \left( v_{A(i)}^{-1} v_{B(i)}^{-1} v_{A(i)}^{-1} \right)^2 \triangleright \left( \varphi_1 \otimes \ldots \otimes \varphi_g \right) \\
&= \varphi_1\!\left(S^{-1}\!\left(v^{(2i-2)-1}\right) ? v^{(2i-1)-1} \right) \otimes \ldots \otimes \varphi_{i-1}\!\left(S^{-1}\!\left( v^{(2)-1} \right) ?  v^{(3)-1} \right)\\
&\:\:\:\:\: \otimes \varphi_i\!\left(S^{-1}(a_j)  S^{-1}\!\left(v^{(1)-1}\right) a_k ? b_k b_j \right) \otimes \varphi_{i+1} \ldots \otimes \varphi_g,
\end{align*}
which is the announced formula.
\finDemo

\section{Equivalence with the Lyubashenko representation}
\indent In a series of papers \cite{lyu95a, lyu95b, lyu96}, V. Lyubashenko has constructed projective representations of $\mathrm{MCG}(\Sigma_{g,n})$ by categorical techniques based on the coend of a ribbon category. Our assumptions on $H$ allow to apply his construction to $\mathrm{mod}_l(H)$, the ribbon category of finite-dimensional left $H$-modules. Here we will show that these two representations are equivalent. For the case of the torus $\mathrm{MCG}(\Sigma_{1,0}^{\mathrm{o}})$ and $\mathrm{MCG}(\Sigma_{1,0})$, we have already shown in \cite{Fai18} that the projective representation obtained thanks to $\mathcal{L}_{1,0}(H)$ is equivalent to the Lyubashenko-Majid representation \cite{LM}.

\subsection{The Lyubashenko representation for $\mathrm{mod}_l(H)$}
\indent Let us first quickly recall the Lyubashenko representation in the general framework of a ribbon category $\mathcal{C}$ satisfying some assumptions (see \cite{lyu95b}). 
\smallskip\\
\indent Let $K = \int^X \! X^* \otimes X$ be the coend of the functor $F : \mathcal{C}^{\mathrm{op}} \times \mathcal{C} \to \mathcal{C}$, $F(X,Y) = X^* \otimes Y$ and let $i_X : X^* \otimes X \to K$ be the associated dinatural transformation (see \cite[IX.6]{ML}). Thanks to the universal property of the coend $K$, Lyubashenko defined several morphisms; we will need some of them which we recall now. The first is an algebra structure $K \otimes K \to K$. Consider the following family of morphisms (for each $X,Y \in \mathcal{C}$)
\begin{equation}\label{dinatProduit}
\begin{split}
d_{X,Y} : X^* \otimes X  \otimes Y^* \otimes Y \xrightarrow{\mathrm{id}_{X^*} \otimes c_{X,Y^*} \otimes \mathrm{id}_Y} X^* \otimes Y^*  \otimes X \otimes Y  &\xrightarrow{\mathrm{id}_{X^*} \otimes \mathrm{id}_{Y^*} \otimes c_{X,Y}} X^* \otimes Y^*  \otimes Y \otimes X\\
&\:\:\:\:\:\:\:\:\: \xrightarrow{\sim} (Y \otimes X)^*  \otimes Y \otimes X \xrightarrow{i_{Y \otimes X}} K.
\end{split}
\end{equation}
Since the family $d_{X,Y}$ is dinatural in $X$ and $Y$, it exists a unique $m_K : K \otimes K \to K$ such that $d_{X,Y} = m_K \circ (i_X \otimes i_Y)$, which is in fact an associative product on $K$. Actually, $K$ is endowed with a Hopf algebra structure whose structure morphisms are similarly defined using the universal property, but we do not need this here.
\smallskip\\
\indent Next, consider the following families of morphisms
\begin{equation}\label{dinatRep}
\begin{split}
&\alpha_X : X^* \otimes X \xrightarrow{\theta_{X^*} \otimes \mathrm{id}_X} X^* \otimes X \xrightarrow{i_X} K,\\
&\beta_{X,Y} : X^* \otimes X \otimes Y^* \otimes Y \xrightarrow{\mathrm{id}_{X^*} \otimes \left(c_{Y^*, X} \circ c_{X, Y^*}\right) \otimes \mathrm{id}_Y} X^* \otimes X \otimes Y^* \otimes Y \xrightarrow{i_X \otimes i_Y} K,\\
&\gamma_X^Y : X^* \otimes X \otimes Y \xrightarrow{\mathrm{id}_{X^*} \otimes \left(c_{Y, X} \circ c_{X, Y}\right)} X^* \otimes X \otimes Y \xrightarrow{i_X \otimes \mathrm{id}_Y} K \otimes Y.
\end{split}
\end{equation}
\indent The families $\alpha_X$ and $\gamma_X^Y$ (with $Y$ fixed) are dinatural in $X$, and the family $\beta_{X,Y}$ is dinatural in $X,Y$. Hence by the universal property of $K$, there exists unique morphisms $\mathcal{T} : K \to K, \mathcal{O} : K \otimes K \to K \otimes K, \mathcal{Q}_Y : K \otimes Y \to K \otimes Y$ such that
\begin{equation}\label{morphismesRep}
\alpha_X = \mathcal{T} \circ i_X, \:\:\:\:\: \beta_{X,Y} = \mathcal{O} \circ (i_X \otimes i_Y), \:\:\:\:\: \gamma_X^Y = \mathcal{Q}_Y \circ (i_X \otimes \mathrm{id}_Y).
\end{equation}

\indent Finally, the morphism $\mathcal{S} : K \to K$ is defined by $\mathcal{S} = (\varepsilon_K \otimes \mathrm{id}_K) \circ \mathcal{O} \circ (\mathrm{id}_K \otimes \Lambda_K)$, where $\Lambda_K$ is the two-sided cointegral on $K$. 
\smallskip\\
\indent Let $X$ be any object of $\mathcal{C}$ and $V_X = \mathrm{Hom}_{\mathcal{C}}(X, K^{\otimes g})$. The Lyubashenko representation $Z_X : \mathrm{MCG}(\Sigma_g^{\mathrm{o}}) \to \mathrm{PGL}(V_X)$ \cite[Section 4.4]{lyu95b} takes the following values:
\begin{equation}\label{repLyubashenkoC}
\begin{split}
Z_X(\tau_{a_i}) &= \mathrm{Hom}_{\mathcal{C}}\!\left(X, \mathrm{id}_K^{\otimes g-i} \otimes \mathcal{T} \otimes \mathrm{id}_K^{\otimes i-1}\right),\\
Z_X(\tau_{b_i}) &= \mathrm{Hom}_{\mathcal{C}}\!\left(X, \mathrm{id}_K^{\otimes g-i} \otimes (\mathcal{S}^{-1} \circ \mathcal{T} \circ \mathcal{S}) \otimes \mathrm{id}_K^{\otimes i-1}\right),\\
Z_X(\tau_{d_i}) &= \mathrm{Hom}_{\mathcal{C}}\!\left(X, \mathrm{id}_K^{\otimes g-i} \otimes (\mathcal{O} \circ (\mathcal{T} \otimes \mathcal{T})) \otimes \mathrm{id}_K^{\otimes i-2}\right) \:\:\: \text{for } i \geq 2, \\
Z_X(\tau_{e_i}) &= \mathrm{Hom}_{\mathcal{C}}\!\left(X, \mathrm{id}_K^{\otimes g-i} \otimes \left((\mathcal{T} \otimes \theta_{K^{\otimes i-1}}) \circ \mathcal{Q}_{K^{\otimes i-1}}\right) \right) \:\:\: \text{for } i \geq 2.
\end{split}
\end{equation}
Recall that the curves $a_i, b_i, d_i, e_i$ are represented in Figure \ref{figureCourbesCanoniques}. Since these Dehn twists are a generating set, we have an operator $Z_X(f)$ for all $f \in \mathrm{MCG}(\Sigma_g^{\mathrm{o}})$. If moreover we take $X = \mathbf{1\!\!\!1}$, the unit object of $\mathcal{C}$, then this defines a representation $Z_{\mathbf{1\!\!\!1}} : \mathrm{MCG}(\Sigma_g) \to \mathrm{PGL}(V_{\mathbf{1\!\!\!1}})$ of the mapping class group of $\Sigma_g$.
\smallskip\\
\indent Now, let us explicit the above formulas to the case of $\mathcal{C} = \mathrm{mod}_l(H)$. Recall that the category $\mathrm{mod}_l(H)$ has braiding $c_{I,J} : X \otimes Y \to Y \otimes X$ and twist $\theta_X : X \to X$ given by
\[ c_{X,Y}(x \otimes y) = b_i \cdot y \otimes a_i \cdot x, \:\:\:\: \theta_X(x) = v^{-1} \cdot x \]
and that the action on the dual module $V^*$ is $h \cdot \varphi = \varphi(S(h) \cdot ?)$ for all $\varphi \in V^*, h \in H$, see \cite[Chap. XIII--XIV]{kassel} for more details.
\smallskip\\
\indent It is well-known (and not difficult to see) that $K$ is $H^*$ endowed with the coadjoint action:
\[ \forall \, h \in H, \:\forall \, \varphi \in K, \:\: h \varphi = \varphi(S(h')?h'') \]
and that the dinatural transformation of $K$ is
\[ i_X(\psi \otimes x) = \psi(? \cdot x) \in K. \]
Note that $\psi(? \cdot x)$ is just a matrix coefficient of the module $X$. The dinatural family $d_{X,Y}$ of \eqref{dinatProduit} is
\[ d_{X,Y}(\varphi \otimes x \otimes \psi \otimes y) = \psi(S(b_i)?b_j \cdot y) \varphi(?a_ja_i \cdot x) \]
where in the right of the equality it is the usual product in $H^*$: $\langle fg, h\rangle = f(h')g(h'')$. To compute the product $m_K$ in $K$ explicitly, observe that $i_{H_{\mathrm{reg}}}(\varphi \otimes 1) = \varphi$, where $H_{\mathrm{reg}}$ is the regular representation of $H$. Thus
\begin{align*}
m_K(\varphi \otimes \psi) = m_K \circ (i_{H_{\mathrm{reg}}} \otimes i_{H_{\mathrm{reg}}})\!\left( \varphi \otimes 1 \otimes \psi \otimes 1 \right) = d_{H_{\mathrm{reg}}, H_{\mathrm{reg}}}(\varphi \otimes 1 \otimes \psi \otimes 1)
&= \psi(S(b_i)?b_j) \varphi(?a_ja_i)\\
&= \varphi(a_j ? a_i) \psi(S(b_i)b_j ?)
\end{align*}
where we used $R \Delta = \Delta^{\mathrm{op}} R$ for the last equality. Moreover, the unit element of $K$ is $1_K = \varepsilon$, the counit of $H$. We record the following lemma, already given in \cite{lyu95b}.

\begin{lemme}\label{cointegraleBilatere}
Assume $\mathcal{C} = \mathrm{mod}_l(H)$, and let $\mu^r \in H^*$ be the right integral on $H$ (unique up to scalar). Then $\mu^r$ is the two-sided cointegral in $K$ (unique up to scalar):
\[ \forall\, \varphi \in K, \:\:\: m_K(\mu^r \otimes \varphi) = m_K(\varphi \otimes \mu^r) = \varepsilon_K(\varphi)\mu^r \]
where $\varepsilon_K(\varphi) = \varphi(1)$.
\end{lemme}
\debutDemo
Using \eqref{quasiCyclic}, we get
\[ m_K(\mu^r \otimes \varphi) =  \mu^r(a_j ? a_i) \varphi(S(b_i)b_j ?) = \mu^r(S^2(a_i)a_j ?) \varphi(S(b_i)b_j ?) = \mu^r \varphi = \varphi(1)\mu^r. \]
We used \eqref{quasiCyclic} and the basic properties of $R$ \cite[VIII.2]{kassel}. Similarly:
\[ m_K(\varphi \otimes \mu^r) = \mu^r(S(b_i)?b_j) \varphi(?a_ja_i) = \mu^r(S^2(b_jS^{-1}(b_i))?) \varphi(?a_ja_i) = \mu^r \varphi = \varphi(1) \mu^r. \]
\finDemo

\indent The dinatural families of \eqref{dinatRep} are
\begin{equation*}
\begin{split}
&\alpha_X(\varphi \otimes x) = \varphi(v^{-1}? \cdot x), \\
&\beta_{X,Y}(\varphi \otimes x \otimes \psi \otimes y) = \varphi(?b_j a_i \cdot x) \otimes \psi(S(a_jb_i)? \cdot y) = \varphi(?v'^{-1}v \cdot x) \otimes \psi(S(v''^{-1})v? \cdot y), \\
&\gamma_X^Y(\varphi \otimes x \otimes y) = \varphi(?b_j a_i \cdot x) \otimes a_j b_i \cdot y = \varphi(?v'^{-1}v \cdot x) \otimes v''^{-1}v \cdot y
\end{split}
\end{equation*}
where we used \eqref{ribbon}. It follows that the morphisms defined in \eqref{morphismesRep} are
\begin{equation*}
\begin{split}
&\mathcal{T}(\varphi) = \mathcal{T} \circ i_{H_{\mathrm{reg}}}(\varphi \otimes 1) = \alpha_{H_{\mathrm{reg}}}(\varphi \otimes 1) = \varphi(v^{-1}?), \\
&\mathcal{O}(\varphi \otimes \psi) = \mathcal{O} \circ (i_{H_{\mathrm{reg}}} \otimes i_{H_{\mathrm{reg}}})(\varphi \otimes 1 \otimes \psi \otimes 1) = \beta_{H_{\mathrm{reg}}, H_{\mathrm{reg}}}(\varphi \otimes 1 \otimes \psi \otimes 1) = \varphi(?v'^{-1}v) \otimes \psi(S(v''^{-1})v?), \\
&\mathcal{Q}_Y(\varphi \otimes y) = \mathcal{Q}_Y \circ (i_{H_{\mathrm{reg}}} \otimes \mathrm{id}_Y)(\varphi \otimes 1 \otimes y) = \gamma_{H_{\mathrm{reg}}}^Y(\varphi \otimes 1 \otimes y) = \varphi(?v'^{-1}v) \otimes v''^{-1}v \cdot y.
\end{split}
\end{equation*}
\noindent Note that $(\mathcal{T} \otimes \theta_Y) \circ \mathcal{Q}_Y(\varphi \otimes y) = \varphi(? v'^{-1}) \otimes v''^{-1} \cdot y$ (see \eqref{repLyubashenkoC}). Finally, thanks to Lemma \ref{cointegraleBilatere}, the morphism $\mathcal{S}$ is
\[ \mathcal{S}(\varphi) = \varphi\!\left(v'^{-1}v\right) \mu^r\!\left(S(v''^{-1})v?\right) = \varphi\!\left(S^{-1}(v''^{-1})v\right) \mu^r\!\left(v'^{-1}v?\right) \]
where the second equality is due to the equality $v'^{-1} \otimes S(v''^{-1}) = S^{-1}(v''^{-1}) \otimes v'^{-1}$ (which follows from $S(v^{-1}) = v^{-1}$). Moreover, we will need the following lemma to prove the equivalence of the representations.
\begin{lemme}\label{lemmeST}
Let $\rho$ be the representation of $\mathcal{L}_{1,0}(H)$ on $H^*$, then the following formulas hold:
\begin{equation*}
\begin{split}
&\mathcal{T} = \rho(v_A^{-1}) = (v^{-1})_*, \:\:\:\: \mathcal{S} = \mu^l(v^{-1}) g^{-1}_* \circ \rho(v_A^2v_B) \circ g_*,\\
&\mathcal{S}^{-1} \circ \mathcal{T} \circ \mathcal{S} = (g^{-1} v)_* \circ \rho(v_B^{-1}) \circ (g v^{-1})_*,
\end{split}
\end{equation*}
where $h_*(\varphi) = \varphi(?h)$ for all $h \in H$ and $\varphi \in H^*$.
\end{lemme}
\debutDemo
The formula for $\mathcal{T}$ is obvious. Propositions 4.10 and 5.3 of \cite{Fai18} give $\rho(v_B)$ and then we compute using \eqref{muLmuRgCarre} and \eqref{integraleShifte}:
\begin{align*}
\rho(v_B)(\varphi) &= v_B \triangleright \varphi = \mu^l(v^{-1})^{-1}\left( \mu^l(g^{-1}v^{-1}?)\varphi^v \right)^{v^{-1}} = \mu^l(v^{-1})^{-1}\left( \mu^r(gv^{-1}?)\varphi^v \right)^{v^{-1}}\\
&= \mu^l(v^{-1})^{-1} \mu^r\!\left(v'^{-1}?gv^{-1}\right) \varphi\!\left(S^{-1}(v''^{-1})g^{-1}v\right)\\
&= \mu^l(v^{-1})^{-1} \left( gv^{-2} \right)_*\!\left(\mu^r(v v'^{-1}?)\right) \left\langle g^{-1}_*(\varphi), S^{-1}(v''^{-1})v \right\rangle\\
&=  \mu^l(v^{-1})^{-1} \left( gv^{-2} \right)_* \circ \mathcal{S} \circ g^{-1}_*(\varphi) = \mu^l(v^{-1})^{-1} \rho(v_A^{-2}) \circ g_* \circ \mathcal{S} \circ g^{-1}_*(\varphi)
\end{align*}
where $\varphi^h = \varphi(h?)$ for $h \in H$. The last claimed formula follows from $\mathcal{S} = \mu^l(v^{-1}) (g^{-1}v)_* \circ \rho(v_Av_Bv_A) \circ (gv^{-1})_*$ and the fact that  $v_A, v_B \in \mathcal{L}_{1,0}(H)$ satisfy the braid relation $v_A v_B v_A = v_B v_A v_B$ (see \cite[Prop. 5.5]{Fai18}).
\finDemo

For the representation space, we take $X = H_{\mathrm{reg}}$, so that $V_X = \mathrm{Hom}_H(H_{\mathrm{reg}}, K^{\otimes g}) \cong K^{\otimes g}$. Then by the previous formulas, we get the Lyubashenko projective representation of $\mathrm{MCG}(\Sigma_g^{\mathrm{o}})$ \eqref{repLyubashenkoC} applied to $\mathrm{mod}_l(H)$:
\begin{equation}\label{formulesLyubashenko}
\begin{split}
Z_{H_{\mathrm{reg}}}(\tau_{a_i})(\varphi_1 \otimes \ldots \otimes \varphi_g) &= \varphi_1 \otimes \ldots \otimes \varphi_{g-i+1}(v^{-1}?) \otimes \ldots \otimes \varphi_g, \\
Z_{H_{\mathrm{reg}}}(\tau_{b_i})(\varphi_1 \otimes \ldots \otimes \varphi_g) &= \varphi_1 \otimes \ldots \otimes (g^{-1}v)_* \circ \rho(v_B^{-1}) \circ (gv^{-1})_*(\varphi_{g-i+1}) \otimes \ldots \otimes \varphi_g, \\
Z_{H_{\mathrm{reg}}}(\tau_{d_i})(\varphi_1 \otimes \ldots \otimes \varphi_g) &=  \varphi_1 \otimes \ldots \otimes \varphi_{g-i+1}\!\left(? v'^{-1}\right) \otimes \varphi_{g-i+2}\!\left(S(v''^{-1})?\right) \otimes \ldots \otimes \varphi_g, \\
Z_{H_{\mathrm{reg}}}(\tau_{e_i})(\varphi_1 \otimes \ldots \otimes \varphi_g) &= \varphi_1 \otimes \ldots \otimes \varphi_{g-i} \otimes \varphi_{g-i+1}\!\left(? v^{(1)-1}\right) \otimes \varphi_{g-i+2}\!\left( S(v^{(2)-1}) ? v^{(3)-1} \right) \otimes \ldots \\
&\:\:\:\:\: \otimes \varphi_g\!\left( S(v^{(2i-2)-1}) ? v^{(2i-1)-1} \right),
\end{split}
\end{equation}
with $i \geq 2$ for the two last formulas. If we take $X = \mathbb{C}$, we get 
\[ V_{\mathbb{C}} = \mathrm{Hom}_H(\mathbb{C}, K^{\otimes g}) = (K^{\otimes g})^{\mathrm{inv}} = \left\{ f \in K^{\otimes g} \left| \, \forall\, h \in H, \:\: h \cdot f = \varepsilon(h)f\right.\right\} \]
where by definition of the action of $H$ on $K$, the action of $H$ on $K^{\otimes g}$ is
\begin{equation}\label{actionHKg}
h \cdot \varphi_1 \otimes \ldots \otimes \varphi_g = \varphi_1\!\left(S(h^{(1)}) ? h^{(2)}\right) \otimes \ldots \otimes \varphi_g\!\left(S(h^{(2g-1)}) ? h^{(2g)}\right).
\end{equation}
Then $Z_{\mathbb{C}}$ is a projective representation of $\mathrm{MCG}(\Sigma_g)$ (note that $Z_{\mathbb{C}}$ is just  $Z_{H_{\mathrm{reg}}}$ restricted to $(K^{\otimes g})^{\mathrm{inv}}$).

\medskip

\indent To conclude this section, we explain how to see $\mathcal{L}_{0,1}(H)$ as a coend. Interpreting slightly differently the fusion relation of Definition \ref{defL01}, we can view $\mathcal{L}_{0,1}(H)$ as $H^*$ endowed with a new product. Indeed, we know that $\mathcal{L}_{0,1}(H)$ is generated by matrix coefficients $\overset{I}{M}{^i_j}$ and due to \eqref{isoPsi01}, $\dim(\mathcal{L}_{0,1}(H)) = \dim(H^*)$; hence $\mathcal{L}_{0,1}(H) \cong H^*$ as vector spaces and we identify them: $\overset{I}{M} = \overset{I}{T}$. To avoid confusion, we exceptionally denote by $\star$ (resp. $\ast$) the product of $\mathcal{O}(H)$ (resp. $\mathcal{L}_{0,1}(H)$); both are products on $H^*$ thanks to the identification. Due to \eqref{structureOH} and to obvious commutation relations, the $\mathcal{L}_{0,1}$-fusion relation on $H^*$ is given by
\[ (\overset{J}{a_i})_2 \overset{I}{T}_1 \ast \overset{J}{T}_2 (\overset{I}{b_i})_1 \overset{IJ}{(R')}{^{-1}_{12}}= \overset{I}{T}_1 \star \overset{J}{T}_2. \]
Using the relation $S^{-1}(a_i)a_j \otimes b_j b_i = 1 \otimes 1$ together with obvious commutation relations, we get
\[ \overset{I}{T}_1 \ast \overset{J}{T}_2 = \overset{J}{S^{-1}(a_i)}_2 \overset{I}{T}_1 \star \overset{J}{T}_2 (\overset{I}{b_j})_1 (\overset{J}{a_j})_2 (\overset{I}{b_i})_1 = \overset{I}{T}\!\left( ? b_j b_i \right)_1 \star \overset{J}{T}\!\left(S^{-1}(a_i) ? a_j\right)_2. \]
Since every element of $H^*$ is a linear combination of matrix elements $\overset{I}{T}{^i_j}$ for certain $I,i,j$, the product in $\mathcal{L}_{0,1}(H)$ is
\begin{equation}\label{produitL01Explicite}
\varphi \ast \psi = \varphi\!\left( ? b_j b_i \right) \star \psi\!\left(S^{-1}(a_i) ? a_j\right).
\end{equation}
Moreover, we define a left $H$-module structure on $\mathcal{L}_{0,1}(H)$ by $h \cdot \varphi = \varphi \cdot S^{-1}(h) = \varphi\!\left( S^{-1}(h'') ? h' \right)$ (see \eqref{actionL01}). Since $h \cdot (\varphi \ast \psi) = (h'' \cdot \varphi) \ast (h' \cdot \psi)$, $\mathcal{L}_{0,1}(H)$ is an algebra in $\mathrm{mod}_l(H^{\mathrm{cop}})$, where $H^{\mathrm{cop}}$ is $H$ with opposite coproduct. Moreover, in $H^{\mathrm{cop}}$, we replace $\Delta$ by $\Delta^{\mathrm{op}}$, $R$ by $R'$ and $S$ by $S^{-1}$ so that the formulas for the product and the $H$-action in the coend of $\mathrm{mod}_l(H^{\mathrm{cop}})$ are exactly those of $\mathcal{L}_{0,1}(H)$. We state this as a proposition.
\begin{proposition} It holds:
\[ \mathcal{L}_{0,1}(H) = \int^{X \in \mathrm{mod}_l(H^{\mathrm{cop}})} X^* \otimes X. \]
\end{proposition}

\subsection{Equivalence of the representations}
\indent Recall that we denote by $\rho$ (resp. $\rho_{\mathrm{inv}}$) the representation of $\mathcal{L}_{g,0}(H)$ on $(H^*)^{\otimes g}$ (resp. $\mathrm{Inv}\!\left((H^*)^{\otimes g}\right)$). Also recall the map $F : H^* \to H^*$ 
\[ F(\varphi) = \varphi\!\left(a_i ? b_i\right), \:\:\:\:\: F^{-1}(\varphi) = \varphi\!\left(S^{-1}(a_i) ? b_i\right) \] 
(already used in the proof of Theorem \ref{formulesExplicites}) and let $\sigma : (H^*)^{\otimes g} \to (H^*)^{\otimes g}$ be the permutation
\[ \sigma(\varphi_1 \otimes \varphi_2 \otimes \ldots \otimes \varphi_{g-1} \otimes \varphi_g) = \varphi_g \otimes \varphi_{g-1} \otimes \ldots \otimes \varphi_2 \otimes \varphi_1. \]
It satisfies $\sigma^{-1} = \sigma$.
\begin{theoreme}\label{thmEquivalenceReps}
The representation of Theorem \ref{thmRepMCG} and the Lyubashenko representation of $\mathrm{MCG}(\Sigma_g^{\mathrm{o}})$ and $\mathrm{MCG}(\Sigma_g)$ are equivalent. More precisely:\\
1) The isomorphism of vector spaces
\[ \fonc{(F \circ S)^{\otimes g} \circ \sigma}{K^{\otimes g}}{(H^*)^{\otimes g}}{\varphi_1 \otimes \ldots \otimes \varphi_g}{\varphi_g\!\left(b_i S(?) a_i\right) \otimes \ldots \varphi_1\!\left(b_i S(?) a_i\right)} \]
is an intertwiner between the two representations:% between the Lyubashenko representation and the $\mathcal{L}_{g,0}$-representation of $\mathrm{MCG}(\Sigma_g^{\mathrm{o}})$:
\[ \left[(F \circ S)^{\otimes g} \circ \sigma \right] \circ Z_{H_{\mathrm{reg}}}(f) = \rho(\widehat{f}) \circ \left[(F \circ S)^{\otimes g} \circ \sigma \right]. \]
2) The isomorphism of vector spaces
\[ (F \circ S)^{\otimes g} \circ \sigma : (K^{\otimes g})^{\mathrm{inv}} \rightarrow \mathrm{Inv}\!\left((H^*)^{\otimes g}\right) \]
is an intertwiner between the two representations:% between the Lyubashenko representation and the $\mathcal{L}_{g,0}$-representation of $\mathrm{MCG}(\Sigma_g)$:
\[ \left[(F \circ S)^{\otimes g} \circ \sigma \right] \circ Z_{\mathbb{C}}(f) = \rho_{\mathrm{inv}}(\widehat{f}) \circ \left[(F \circ S)^{\otimes g} \circ \sigma \right]. \]
\end{theoreme}
\debutDemo
1) We show that this isomorphism intertwines the formulas of Theorem \ref{formulesExplicites} and of \eqref{formulesLyubashenko}. Thanks to the properties of $v$ \eqref{ribbon}, it is clear that $(F \circ S)^{\otimes g} \circ \sigma \circ Z(\tau_{a_i}) = \rho(v_{A(i)}^{-1}) \circ (F \circ S)^{\otimes g} \circ \sigma$. Next, thanks to \eqref{actionvAvB}, \eqref{muLmuRgCarre} and \eqref{integraleShifte}, we have
\[ v_B^{-1} \triangleright \varphi = \mu^l(v)^{-1} \mu^r\!\left(gv^{-1} v'?\right) \varphi\!\left( vS^{-1}(gv'') \right). \]
Hence, for $\varphi \in H^*$,
\begin{align*}
\rho(v_B^{-1}) \circ (F \circ S)(\varphi) &= \mu^l(v)^{-1} \mu^r\!\left(gv^{-1}v'?\right) \varphi\!\left( v b_i g v'' a_i \right) = \mu^l(v)^{-1} \mu^r\!\left(g\overline{Y}_j?\right) \varphi\!\left( v^2 b_i g \overline{X}_j a_i \right)\\
&= \mu^l(v)^{-1} \mu^r\!\left(g S(a_j)S^{-1}(b_k)?\right) \varphi\!\left( v^2 g S^{-2}(b_i) b_j a_k a_i \right) = (\star)
\end{align*}
with $\overline{X}_i \otimes \overline{Y}_i = (RR')^{-1}$. We have a Yang-Baxter relation
\begin{align*}
S(a_j)S^{-1}(b_k) \otimes S^{-2}(b_i)b_j \otimes a_k a_i &= a_jS^{-1}(b_k) \otimes S^{-1}\!\left(b_j S^{-1}(b_i)\right) \otimes a_k a_i = (\mathrm{id} \otimes S^{-1} \otimes \mathrm{id})(R_{12} R_{31}^{-1} R_{32}^{-1})\\
&= (\mathrm{id} \otimes S^{-1} \otimes \mathrm{id})( R_{32}^{-1} R_{31}^{-1} R_{12}) = S^{-1}(b_k)S(a_j) \otimes b_j S^{-2}(b_i) \otimes a_i a_k
\end{align*}
which allows us to continue the computation:
\begin{align*}
(\star) &= \mu^l(v)^{-1} \mu^r\!\left(g S^{-1}(b_k)S(a_j)?\right) \varphi\!\left( v^2 g b_j S^{-2}(b_i) a_i a_k \right) = \mu^l(v)^{-1} \mu^r\!\left(g S^{-1}(b_k)S(a_j)?\right) \varphi\!\left( v S^2(b_j) a_k \right)\\
&= \mu^l(v)^{-1} \mu^r\!\left(g S^{-1}(a_j b_k)?\right) \varphi\!\left( v b_j a_k \right) = \mu^l(v)^{-1} \mu^r\!\left(g v S^{-1}(v''^{-1})?\right) \varphi\!\left( v^2 v'^{-1} \right).
\end{align*}
We used \eqref{u} and \eqref{ribbon}. On the other hand, we compute
\begin{align*}
(F \circ S) \circ Z(\tau_b)(\varphi) &= (F \circ S) \circ (\mathcal{S}^{-1} \circ \mathcal{T} \circ \mathcal{S})(\varphi) = (F \circ S) \circ (g^{-1} v)_* \circ \rho(v_B^{-1}) \circ (g v^{-1})_*(\varphi)\\
&= F \circ S\!\left( \mu^l(v)^{-1} \mu^r\!\left( v'?\right) \varphi\!\left( S^{-1}(v'') \right) \right) = \mu^l(v)^{-1} \mu^r\!\left( v' b_i S(?) a_i\right) \varphi\!\left( S^{-1}(v'') \right)\\
&= \mu^l(v)^{-1} \mu^r\!\left( v S^2(a_i) \overline{Y}_j b_i S(?)\right) \varphi\!\left( v S^{-1}(\overline{X}_i) \right)\\
&= \mu^l(v)^{-1} \mu^r\!\left( v S^2(a_i) S(a_j) S^{-1}(b_k) b_i S(?)\right) \varphi\!\left( v S^{-1}(b_ja_k) \right) = (\star\star).
\end{align*}
We used Lemma \ref{lemmeST}, \eqref{quasiCyclic} and \eqref{ribbon}. As previously, we have a Yang-Baxter relation
\[ S^2(a_i) S(a_j) \otimes S^{-1}(b_k) b_i \otimes b_ja_k = S(a_j)S^2(a_i) \otimes b_i S^{-1}(b_k) \otimes a_k b_j \]
which allows us to continue the computation:
\begin{align*}
(\star \star) &= \mu^l(v)^{-1} \mu^r\!\left( v S(a_j) S^2(a_i) b_i S^{-1}(b_k) S(?)\right) \varphi\!\left( v S^{-1}(a_k b_j) \right)\\
& = \mu^l(v)^{-1} \mu^r\!\left(  S(a_j) g S^{-1}(b_k) S(?)\right) \varphi\!\left( v S^{-1}(a_k b_j) \right) = \mu^l(v)^{-1} \mu^r\!\left( g a_j b_k S(?)\right) \varphi\!\left( v b_j a_k \right)\\
&= \mu^l(v)^{-1} \mu^r\!\left( g v v''^{-1} S(?)\right) \varphi\!\left( v^2 v'^{-1} \right) = \mu^l(v)^{-1} \mu^r \circ S\!\left( ? S^{-1}(v''^{-1})vg^{-1}\right) \varphi\!\left( v^2 v'^{-1} \right)\\
&= \mu^l(v)^{-1} \mu^l\!\left( ? S^{-1}(v''^{-1})vg^{-1}\right) \varphi\!\left( v^2 v'^{-1} \right) = \mu^l(v)^{-1} \mu^r\!\left( gv S^{-1}(v''^{-1})? \right) \varphi\!\left( v^2 v'^{-1} \right).
\end{align*}
We used \eqref{u} to simplify $S^2(a_i) b_i = S(S^{-1}(b_i)S(a_i)) = gv^{-1}$ and the properties of $\mu^l$ and $\mu^r$ recorded in section \ref{preliminaries}. Hence, it holds $\rho(v_B^{-1}) \circ (F \circ S) = (F \circ S) \circ Z(\tau_b)$, which clearly implies that $\rho(v_{B(i)}^{-1}) \circ (F \circ S)^{\otimes g} \circ \sigma = (F \circ S)^{\otimes g} \circ \sigma \circ Z(\tau_{b_i})$. Let us now proceed with $\tau_{d_i}$ ($i \geq 2$):
\begin{align*}
 &(F \circ S)^{\otimes g} \circ \sigma \circ Z(\tau_{d_i}) \circ \sigma \circ (S^{-1} \circ F^{-1})^{\otimes g}\!\left(\varphi_1 \otimes \ldots \varphi_g\right)\\
&=(F \circ S)^{\otimes g} \circ \sigma \circ Z(\tau_{d_i})\!\left( \varphi_g\!\left(S^{-1}(a_j) S^{-1}(?) b_j\right) \otimes \ldots \otimes \varphi_1\!\left(S^{-1}(a_j) S^{-1}(?) b_j\right) \right)\\
&= (F \circ S)^{\otimes g} \circ \sigma \!\left( \varphi_g\!\left(S^{-1}(a_j) S^{-1}(?) b_j\right) \otimes \ldots \otimes \varphi_i\!\left(S^{-1}(a_j) S^{-1}(v'^{-1})S^{-1}(?) b_j\right) \otimes  \varphi_{i-1}\!\left(S^{-1}(a_j) S^{-1}(?) v''^{-1} b_j\right) \right.\\
&\:\:\:\:\:\:\:\:\:\:\:\:\:\:\:\:\:\:\:\:\:\:\:\:\:\:\:\:\:\:\:\:\: \left. \otimes \ldots \otimes \varphi_1\!\left(S^{-1}(a_j) S^{-1}(?) b_j\right) \right)\\
&=  \varphi_1 \otimes \ldots \otimes \varphi_{i-1}\!\left(S^{-1}(a_j)a_k ? b_k v''^{-1} b_j \right) \otimes \varphi_i\!\left(S^{-1}(a_j) S^{-1}(v'^{-1})a_k ? b_k b_j \right) \otimes \ldots \otimes \varphi_g\\
&= \rho(v_{D_i}^{-1})\!\left(\varphi_1 \otimes \ldots \varphi_g\right).
\end{align*}
Finally, for $\tau_{e_i}$ ($i \geq 2$):
\begin{align*}
 &(F \circ S)^{\otimes g} \circ \sigma \circ Z(\tau_{e_i}) \circ \sigma \circ (S^{-1} \circ F^{-1})^{\otimes g}\!\left(\varphi_1 \otimes \ldots \varphi_g\right)\\
&= (F \circ S)^{\otimes g} \circ \sigma \circ Z(\tau_{e_i})\!\left( \varphi_g\!\left(S^{-1}(a_j) S^{-1}(?) b_j \right) \otimes \ldots \otimes \varphi_1\!\left(S^{-1}(a_j) S^{-1}(?) b_j\right) \right)\\
&= (F \circ S)^{\otimes g} \circ \sigma\!\left( \varphi_g\!\left(S^{-1}(a_j) S^{-1}(?) b_j\right) \otimes \ldots \otimes \varphi_{i+1}\!\left(S^{-1}(a_j) S^{-1}(?) b_j\right) \otimes \varphi_{i}\!\left(S^{-1}(a_j)S^{-1}(v^{(1)-1}) S^{-1}(?) b_j\right) \right. \\
& \:\:\:\:\:\:\:\:\:\:\:\:\:\:\:\:\: \left.  \otimes \, \varphi_{i-1}\!\left( S^{-1}(a_j) S^{-1}(v^{(3)-1}) S^{-1}(?) v^{(2)-1} b_j \right) \otimes \ldots \otimes \varphi_1\!\left( S^{-1}(a_j) S^{-1}(v^{(2i-1)-1})  S^{-1}(?) v^{(2i-2)-1} b_j \right) \right)\\
&= \varphi_1\!\left( S^{-1}(a_j) S^{-1}(v^{(2i-1)-1}) a_k ? b_k v^{(2i-2)-1} b_j \right) \otimes \ldots \otimes \varphi_{i-1}\!\left( S^{-1}(a_j) S^{-1}(v^{(3)-1}) a_k ? b_k v^{(2)-1} b_j \right) \\
& \:\:\:\: \otimes \varphi_{i}\!\left(S^{-1}(a_j)S^{-1}(v^{(1)-1}) a_k ? b_k b_j\right) \otimes \varphi_{i+1} \otimes \ldots \otimes \varphi_g\\
&= \varphi_1\!\left( S^{-1}(v^{(2i-2)-1})  ? v^{(2i-1)-1} \right) \otimes \ldots \otimes \varphi_{i-1}\!\left( S^{-1}(v^{(2)-1}) ? v^{(3)-1} \right) \otimes \varphi_{i}\!\left(S^{-1}(a_j)S^{-1}(v^{(1)-1}) a_k ? b_k b_j\right) \otimes \varphi_{i+1}\\
&\:\:\:\: \otimes \ldots \otimes \varphi_g\\
&= \rho(v_{E_i}^{-1})\!\left(\varphi_1 \otimes \ldots \varphi_g\right).
\end{align*}
We used $\Delta^{\mathrm{op}}R = R\Delta$ for the last equality.\\
\noindent 2) It is not difficult to see that $(F \circ S)^{\otimes g} \circ \sigma : K^{\otimes g} \to (H^*)^{\otimes g}$ is a morphism of $H$-modules, where $K^{\otimes g}$ is endowed with the action \eqref{actionHKg} and $(H^*)^{\otimes g}$ is endowed with the action \eqref{actionH} (with $n=0$). Hence, the restriction of $(F \circ S)^{\otimes g} \circ \sigma$ to $(K^{\otimes g})^{\mathrm{inv}}$ indeed takes values in $\mathrm{Inv}\!\left((H^*)^{\otimes g}\right)$. Since 
\[ Z_{\mathbb{C}}(f) = \left(Z_{H_{\mathrm{reg}}}(f)\right)_{\vert \, (K^{\otimes g})^{\mathrm{inv}}} \:\:\: \text{ and } \:\:\: \rho_{\mathrm{inv}}(\widehat{f}) = \rho(\widehat{f})_{\vert \, \mathrm{Inv}\!\left((H^*)^{\otimes g}\right)}, \]
the result follows from the first part of the theorem.
\finDemo


\begin{thebibliography}{30}
\setlength{\itemsep}{0\baselineskip}
{\small
\bibitem[AGPS18]{AGPS} N. Aghaei, A.M. Gainutdinov, M. Pawelkiewicz, V. Schomerus, {\em Combinatorial Quantisation of $\mathrm{GL}(1|1)$ Chern-Simons Theory I: The Torus}, arXiv:1811.09123.
\bibitem[Ale94]{alekseev} A.Y. Alekseev, {\em Integrability in the Hamiltonian Chern-Simons theory}, Algebra i Analiz 6, Issue 2 (1994), 53--66.
\bibitem[AGS95]{AGS} A.Y. Alekseev, H. Grosse, V. Schomerus, {\em Combinatorial quantization of the Hamiltonian Chern-Simons theory, I}, Commun. Math. Phys. 172, Issue 2 (1995), 317--358.
\bibitem[AGS96]{AGS2} A.Y. Alekseev, H. Grosse, V. Schomerus, {\em Combinatorial quantization of the Hamiltonian Chern-Simons theory, II}, Commun. Math. Phys. 174, Issue 3 (1996), 561--604.
\bibitem[AS96a]{AS} A.Y. Alekseev, V. Schomerus, {\em Representation theory of Chern-Simons observables}, Duke Math. J. 85, no. 2 (1996), 447--510. 
\bibitem[AS96b]{AS2} A.Y. Alekseev, V. Schomerus, {\em Quantum Moduli Spaces of Flat Connections}, arXiv:q-alg/9612037.
\bibitem[Ari10]{arike} Y. Arike, {\em A construction of symmetric linear functions on the restricted quantum group $\overline{U}\!_q(\mathfrak{sl}(2))$}, Osaka J. Math. 47 (2010), 535-557.
\bibitem[BaR]{BaR} S. Baseilhac, P. Roche, {\em Non-restricted moduli algebras of punctured spheres}, in preparation.
\bibitem[BZBJ18]{BZBJ} D. Ben-Zvi, A. Brochier, D. Jordan, {\em Integrating quantum groups over surfaces}, Journal of Topology 11, Issue 4 (2018),  874--917.
%\bibitem[BJ17]{BJ} A. Brochier, D. Jordan, {\em Fourier transform for quantum D-modules via the punctured torus mapping class group}, Quantum Topol. 8 (2017), 361-379.
\bibitem[BR95]{BR} E. Buffenoir, P. Roche, {\em Two dimensional lattice gauge theory based on a quantum group}, Commun. Math. Phys. 170, Issue 3 (1995), 669--698.
\bibitem[BR96]{BR2} E. Buffenoir, P. Roche, {\em Link invariants and combinatorial quantization of Hamiltonian Chern Simons theory}, Commun. Math. Phys. 181, Issue 2 (1996), 331--365.
\bibitem[BNR02]{BNR} E. Buffenoir, K. Noui, P. Roche, {\em Hamiltonian quantization of Chern-Simons theory with $\mathrm{SL}(2, \mathbb{C})$ group}, Class. Quantum Grav. 19, no. 19 (2002), 4953--5015.
\bibitem[BFKB98a]{BFKB} D. Bullock, C. Frohman, J. Kania-Bartoszynska, {\em Skein Quantization and Lattice Gauge Field Theory}, Chaos, Solitons \& Fractals 9, Issues 4--5 (1998), 811--824.
\bibitem[BFKB98b]{BFKB2} D. Bullock, C. Frohman, J. Kania-Bartoszynska, {\em Topological interpretations of lattice gauge field theory}, Commun. Math. Phys. 198, Issue 1 (1998), 47--81.
%\bibitem[CMR17]{CMR} L. Chekhov, M. Mazzocco, V. Rubtsov, {\em Algebras of quantum monodromy data and decorated character varieties}, arXiv:1705.01447.
\bibitem[CR62]{CR} C. Curtis, I. Reiner, {\em Representation Theory of Finite Groups and Associative Algebras}, Interscience Publishers (1962).
\bibitem[Fai18a]{F} M. Faitg, {\em A note on symmetric linear forms and traces on the restricted quantum group $\overline{U}\!_q(\mathfrak{sl}(2))$}, arXiv:1801.07524.
\bibitem[Fai18b]{Fai18} M. Faitg, {\em Modular group representations in combinatorial quantization with non-semisimple Hopf algebras}, arXiv:1805.00924.
%\bibitem[Fai2]{theseF} M. Faitg, Ph.D. thesis, in preparation.
\bibitem[FM12]{FM} B. Farb, D. Margalit, {\em A Primer on Mapping Class Groups}, Princeton University Press (2012).
\bibitem[FR98]{FockRosly} V.V. Fock, A.A. Rosly, {\em Poisson structure on moduli of flat connections on Riemann surfaces and $r$-matrix}, arXiv:math/9802054.
\bibitem[FSS12]{FSS1} J. Fuchs, C. Schweigert, C. Stigner, {\em Modular invariant Frobenius algebras from ribbon Hopf algebra automorphisms}, J. Algebra 363 (2012), 29--72
\bibitem[FSS14]{FSS2} J. Fuchs, C. Schweigert, C. Stigner, {\em Higher genus mapping class group invariants from factorizable Hopf algebras}, Adv. Math. 250 (2014), 285--319.
\bibitem[FGST06]{FGST} B.L. Feigin, A.M. Gainutdinov, A.M. Semikhatov, I.Y. Tipunin, {\em Modular group representations and fusion in logarithmic conformal field theories and in the quantum group center}, Commun. Math. Phys. 265, Issue 1 (2006), 47-93.
\bibitem[GT09]{GT} A.M. Gainutdinov, I.Y. Tipunin, {\em Radford, Drinfeld and Cardy boundary states in the $(1, p)$ logarithmic conformal field models}, J. Phys. A: Math. Theor 42, 315207 (2009).
%\bibitem[Iba16]{ibanez} E. Ibanez, {\em Evaluable Jones-Wenzl idempotents at root of unity and modular representation on the center of $\overline{U}\!_qsl(2)$}, PhD Thesis (in French), arXiv:1604.03681.
\bibitem[Kas95]{kassel} C. Kassel, {\em Quantum Groups}, Graduate texts in Mathematics 155, Springer (1995).
%\bibitem[KS11]{KS} H. Kondo, Y. Saito, {\em Indecomposable decomposition of tensor products of modules over the restricted quantum group associated to $\mathfrak{sl}_2$}, Journal of Algebra 330, Issue 1 (2011), 103-129.
%\bibitem[LM94]{LM} V. Lyubashenko, S. Majid, {\em Braided Groups and Quantum Fourier Transform}, Journal of Algebra 166, Issue 3 (1994), 506-528.
\bibitem[Lab13]{labourie} F. Labourie, {\em Lectures on Representations of Surface Groups}, EMS publishing house, Zurich Lectures in Advanced Mathematics, 2013.
\bibitem[Lyu95a]{lyu95a} V. Lyubashenko, {\em Modular transformations for tensor categories}, J. Pure Appl.Alg. 98, Issue 3 (1995) 279--327.
\bibitem[Lyu95b]{lyu95b} V. Lyubashenko, {\em Invariants of 3-manifolds and projective representations of mapping class groups via quantum groups at roots of unity}, Commun. Math. Phys. 172, Issue 3 (1995), 467--516.
\bibitem[Lyu96]{lyu96} V. Lyubashenko, {\em Ribbon abelian categories as modular categories}, J. Knot Theory and its Ramif. 5, no. 03 (1996), 311--403.
\bibitem[LM94]{LM} V. Lyubashenko, S. Majid, {\em Braided Groups and Quantum Fourier Transform}, Journal of Algebra 166, Issue 3 (1994), 506-528.
\bibitem[ML98]{ML} S. Mac Lane, {\em Categories for the working mathematician}, 2nd edition, Graduate texts in Mathematics 5, Springer (1998).
\bibitem[Maj93]{majid93} S. Majid, {\em Braided groups}, J. Pure Appl. Algebra 86, Issue 2 (1993), 187--221.
\bibitem[Maj95]{majid} S. Majid, {\em Foundations of Quantum Group Theory}, Cambridge University Press, 1995.
\bibitem[Mon93]{Mon} S. Montgomery, {\em Hopf algebras and their action on rings}, CBMS Regional Conference Series in Mathematics no. 82, American Mathematical Society.
\bibitem[MW15]{MW} C. Meusburger, D.K. Wise, {\em Hopf algebra gauge theory on a ribbon graph}, arXiv:1512.03966.
%\bibitem[Rad94]{radford} D.E. Radford, {\em The trace function and Hopf algebras}, Journal of Algebra 163 (1994), 583-622.
\bibitem[Sch98]{S} V. Schomerus, {\em Deformed Gauge Symmetry in Local Quantum Physics}, Habilitation thesis, Hamburg (1998).
%\bibitem[Sut94]{suter} R. Suter, {\em Modules over $\bar U_q(\mathfrak{sl}(2))$}, Commun. Math. Phys. 163, Issue 2 (1994), 359-393.
\bibitem[Waj83]{wajnryb} B. Wajnryb, {\em A simple presentation for the mapping class group of an orientable surface}, Israel J. Math. 45, Issue 2--3 (1983), 157--174.
\bibitem[Wit91]{W}, E. Witten, {\em On Quantum Gauge Theories in Two Dimensions}, Commun. Math. Phys. 141, Issue 1 (1991), 153--209.
}
\end{thebibliography}
\end{document}